\documentclass{article}

\usepackage{subfig}
\usepackage{graphicx}
\usepackage{adjustbox} 
\usepackage[dvipsnames]{xcolor}
\usepackage{tabularx}
\usepackage{amsmath,amssymb}
\usepackage{amsthm}
\usepackage{mathrsfs} 
\usepackage{gensymb}
\usepackage[left=2cm, right=2cm, top=2.5cm, bottom=2.5cm]{geometry}
\usepackage{cases}
\usepackage{mathtools} 
\usepackage{gensymb} 
\usepackage{caption} 
\usepackage{bm} 
\usepackage[colorlinks=true, urlcolor=black, linkcolor=black, citecolor=black]{hyperref}
\usepackage{cleveref}
\usepackage{dsfont}
\usepackage{numprint} 
\usepackage{marginnote}
\usepackage{mathtools}
\usepackage{multirow} 
\usepackage{multicol} 
\usepackage{numprint}
\usepackage{longtable}
\usepackage{lineno} 
\usepackage{yfonts}
\usepackage{eufrak}
\usepackage{ulem}
\usepackage{siunitx}
\usepackage{rotating}
\usepackage{algorithm}%
\usepackage{algorithmicx}%
\usepackage{algpseudocode}%
\usepackage{titlesec}
\usepackage{xparse}
\usepackage{comment}

\newcommand{\sib}[1]{[\si{#1}]}

\newcommand{\partialder}[2]{\frac{\partial#1}{\partial#2}}

    \NewDocumentCommand{\secondpartialder}{m m o}{\IfNoValueTF{#3}{
            \frac{\partial^2 #1}{\partial #2^2}
        }{
            \frac{\partial^2 #1}{\partial #2 \partial #3}
        }
    }

\newcommand{\divergence}[0]{\nabla \cdot}

\newcommand{\laplacian}[0]{\nabla^2}

\newcommand{\loss}[0]{\mathscr{L}}
\newcommand{\on}[0]{\ \mathrm{on} \ }

\makeatletter
\newcommand{\slantedparallel}{\mathrel{\mathpalette\new@parallel\relax}}
\newcommand{\new@parallel}[2]{
	\begingroup
	\sbox\z@{$#1T$}
	\resizebox{!}{\ht\z@}{\raisebox{\depth}{$\m@th#1/\mkern-5mu/$}}%
	\endgroup
}
\makeatother

\newcommand{\eps}[0]{\varepsilon}

\newcommand{\real}[0]{\text{Re}}
\newcommand{\imag}[0]{\text{Im}}
\newcommand{\RR}[0]{\mathbb{R}}
\newcommand{\CC}[0]{\mathbb{C}}

\newcommand{\bfu}[0]{\boldsymbol{u}}
\newcommand{\bfsigma}[0]{\boldsymbol{\sigma}}

\newcommand{\bfn}[0]{\boldsymbol{n}}
\newcommand{\bfP}[0]{\boldsymbol{P}}
\newcommand{\bfPhi}[0]{\boldsymbol{\Phi}}
\newcommand{\bfx}[0]{\boldsymbol{x}}

\newcommand{\boundarylapl}[1][]{\mathcal{B}_{{\scriptscriptstyle L}\if\relax\detokenize{#1}\relax\else,#1\fi}}
\newcommand{\boundarymech}[1][]{\mathcal{B}_{{\scriptscriptstyle M}\if\relax\detokenize{#1}\relax\else,#1\fi}
}

\newtheorem{problem}{Problem}

\theoremstyle{definition}              

\newtheorem{remark}{Remark}[section]

\crefname{theorem}{theorem}{theorems}
\Crefname{theorem}{Theorem}{Theorems}

\crefname{lemma}{lemma}{lemmas}
\Crefname{lemma}{Lemma}{Lemmas}

\crefname{proposition}{proposition}{propositions}
\Crefname{proposition}{Proposition}{Propositions}

\crefname{problem}{problem}{problems}
\Crefname{problem}{Problem}{Problems}

\crefname{definition}{definition}{definitions}
\Crefname{definition}{Definition}{Definitions}

\crefname{remark}{remark}{remarks}
\Crefname{remark}{Remark}{Remarks}

\crefname{observation}{observation}{observations}
\Crefname{observation}{Observation}{Observations}

\crefname{example}{example}{examples}
\Crefname{example}{Example}{Examples}

\definecolor{raspberry}{rgb}{1,0,0.2}
\definecolor{darkred}{rgb}{0.5,0.1,0.1}
\definecolor{darkgreen}{rgb}{0.1,0.6,0.2}

\usepackage[numbers,square]{natbib}

\title{A holomorphic neural network framework for 3D boundary value problems governed by harmonic potentials
} 

\author{Enrico Ballini\textsuperscript{1}, Allan Peter Engsig-Karup\textsuperscript{2}, Tito Andriollo\textsuperscript{1}  \\
\small \textsuperscript{1}Department of Mechanical and Production Engineering, Aarhus University, Aarhus, Denmark \\
\small\textsuperscript{2}Department of Applied Mathematics and Computer Science, Technical University of Denmark, Kongens Lyngby, Denmark}

\date{}

\begin{document}

\maketitle



\begin{abstract}
We present a neural-network-based framework for the solution of three-dimensional boundary value problems where the solution is expressible in terms of harmonic potentials. The approach leverages the Whittaker integral formula, which allows representing the solution through functions that are holomorphic with respect to a suitable complex variable. These functions are subsequently approximated using holomorphic neural networks, which guaranty fulfillment of the holomorphicity requirement. A key feature of the proposed formulation is that the governing partial differential equations (PDEs) are satisfied exactly by construction. Therefore, in contrast to standard physics-informed neural networks, no residual minimization of PDEs is required in the interior of the domain, and training is based exclusively on boundary collocation points. The method is validated against three-dimensional Laplace and linear elasticity problems, where, in the latter case, displacement and stress fields are expressed via the Papkovich-Neuber potentials. The numerical results show an accurate approximation of both scalar and vector fields, with errors remaining controlled throughout the domain. Overall, the work demonstrates that the incorporation of analytical structures into neural network architectures provides a natural and effective framework for the meshless approximation of three-dimensional boundary value problems while preserving the underlying properties of the governing equations. 
\end{abstract}

\section{Introduction}

Partial differential equations play a central role in the mathematical description of physical phenomena and engineering processes. Among the many classes of boundary value problems, \textit{harmonic} problems, namely those whose solution can be expressed via harmonic functions, occupy a particularly important position, as they arise in a wide range of contexts. Examples include potential flow in fluid mechanics \cite{landau1987, Engsig-Karup2013}, diffusion and transport phenomena \cite{Quarteroni1994}, elastic deformation of solid materials \cite{Sadd2021}, mechanical bending of thin plates \cite{meleshko_selected_2003}, and even stochastic processes \cite{Redner2001}. \\
The widespread occurrence of harmonic problems, together with the increasing demand for accurate and computationally efficient simulations, motivates the continued development of numerical methods for their solution. In this regard, several methods have been proposed for the numerical solution of PDEs, including the finite element method (FEM) \cite{Quarteroni1994}, the spectral element method \cite{Engsig-Karup2016, Xu2018}, isogeometric analysis \cite{Hughes2005} and the boundary element method \cite{Gu2001}. Neural networks have also been used to solve differential equations \cite{Lee1990, Psichogios1992, Dissanayake1994}. Early approaches combined neural networks with collocation methods to incorporate residuals of differential equations in the objective function, e.g. \cite{Lagaris1998,Lagaris2000} and later \cite{Rudd2013,Berg2018}, providing a basis on which trainable parameters were optimized to satisfy both the governing PDE and the associated boundary conditions. Along these lines, the work of \citet{Raissi2019} revisited the original ideas and utilized modern deep neural network software that includes support for automatic differentiation and introduced the so-called physics-informed neural networks (PINNs). In this framework, neural networks are employed as global ansatz functions for the approximation of the PDE solution, and are trained by minimizing a loss function that incorporates both the PDE residuals and the boundary condition mismatch. Also, it is possible to incorporate a loss function for data mismatch, which makes it a general hybrid modeling framework particularly well-suited to address a broad range of problems.
For example, PINNs have demonstrated remarkable flexibility across a wide range of PDE classes and formulations \cite{Toscano2025}, including fluid flow problems \cite{Wassing2025}, advection--diffusion--reaction systems \cite{Lekaba2026, Huang2024}, harmonic problems \cite{Vaishampayan2024}, and solid mechanics applications \cite{Herrmann2024, Hu2024a}. Among the latter, linear elasticity has attracted particular attention due to its broad range of engineering applications \cite{Hu2024a}, with existing studies addressing both forward  \cite{Henkes2022, Vahab2022, Rezaei2022, Roy2023, Kadlag2026} and inverse problems \cite{Xu2023, Zhang2022a}.

Within the PINN framework, fully-connected real-valued neural networks have proven effective for approximating PDE solutions, as supported by several approximation results, e.g., \cite{Park2020, Cai2022, Li2023, Hanin2017}. However, the convergence of PINN approximations to the exact solution may be slow for certain classes of problems, which can significantly affect the overall efficiency of the method \cite{Mishra2022}. Consequently, problem-specific formulations and architectures are often required to improve the performance of neural-network-based solvers, as research in  several areas indicates \cite{Fukushima1980,Vaswani2017}. \\
One common strategy is to design loss functions that more closely reflect the underlying physical structure of the problem. For example, \citet{Samaniego2020} proposed the so-called deep energy method, which takes advantage of the principle of minimum potential energy to optimize the neural network. Building on this idea, subsequent developments have sought to more closely couple the discretization strategy with the governing equation, for example, by approximating potential functions that inherently satisfy the underlying physical principles \cite{Wang2024a}. \\
Another important research direction concerns the reduction of the highest derivative order appearing in the loss function. This can be achieved by employing weak formulations of the governing equations \cite{Rezaei2022, Tresckow2022} or through variational formulations, as proposed by \citet{Kharazmi2019, Kharazmi2021}. Similarly, \citet{Sun2023} reduced the derivative order by defining the loss function exclusively in terms of the residuals of boundary integral equations. \\
Problem-specific mathematical structures can also be incorporated directly into the neural network architecture, in order to enforce intrinsic properties of the PDE solution. Examples include the work of \citet{Ghosh2023}, where \textit{harmonic} neural networks were proposed to solve two-dimensional Laplace problems, and that of \citet{Calafa2024}, where \textit{holomorphic} neural networks were introduced to solve boundary value problems where the solution can be represented via holomorphic functions. Although problem-specific, this line of research appears particularly promising because it exploits well-established analytical representation frameworks to construct approximations that satisfy the governing PDEs \textit{a priori}. Compared to the standard PINN approach, this results in more efficient training, since evaluations are only required at the domain boundary, and in lower memory requirements due to the reduced number of training points \cite{Calafa2024, Calafa2025a, Calafa2025,Zhou2026, Zhou2026a}. On the other hand, the approaches developed thus far have been limited to two-dimensional settings, which restricts their large-scale applicability.

\paragraph{Contributions and outline.}
Motivated by the approaches proposed in \cite{Ghosh2023, Calafa2024}, we introduce a novel neural-network-based framework for solving harmonic boundary value problems in the general three-dimensional setting. The method is based on the Whittaker integral representation  \cite{Whittaker1903}, which allows generation of three-dimensional harmonic functions from complex-valued functions of two variables -- one complex and one real. To approximate these generator functions, a dedicated complex-valued neural network architecture is developed that is holomorphic in one argument only. \\
The key feature of the proposed approach is that the governing PDEs are satisfied exactly, independently of the training accuracy that results from approximating the solution of these PDEs by a neural network. Indeed, the networks are trained exclusively through the mismatch between the predicted and prescribed values of the boundary conditions, without requiring evaluation of the PDE residual at interior points. Moreover, the method avoids the computation of higher-order derivatives, resulting in improved numerical efficiency. \\
To assess the strengths and limitations of the proposed methodology, standard Laplace and linear elasticity boundary value problems are considered. Both problems are first recast in terms of suitable holomorphic functions and then approximated using the proposed neural network architecture.

The paper is structured as follows. In \Cref{sec:models}, we present the governing equations of the considered boundary value problems. In \Cref{sec:nn}, we describe the adopted numerical strategies together with the proposed neural network architecture. Finally, in \Cref{sec:tests}, we present several test cases to assess the performance of the method in different scenarios.

\section{Mathematical models}\label{sec:models}

We present in this section the mathematical models of the problems of interest, providing their formulation in terms of complex-valued functions. \\
To ease the readability of the following sections, we recall here the definition of a holomorphic function and its properties. Consider an open domain $\Omega_\text{2D} \subset \CC$. Using the standard identification $\CC \simeq \RR^2$, each point $\zeta \in \Omega_\text{2D}$ is written as $\zeta = \xi + i\eta$, with $(\xi,\eta) \in \RR^2$. Let $g:\Omega_\text{2D}\to\CC$ be given by
$g(\zeta)=u(\xi,\eta)+iv(\xi,\eta)$, where $u,v:\RR^2\to\RR$. 
The function $g$ is said to be \textit{holomorphic} in $\Omega_\text{2D}$ if the limit  
\begin{equation*}
	g' = \lim_{\zeta \to \zeta_0}\frac{g(\zeta) - g(\zeta_0)}{\zeta-\zeta_0}
\end{equation*}
exists for any point $\zeta_0 \in \Omega_\text{2D}$, see \cite{Stein2003, Ahlfors2021}. 
A direct consequence of holomorphicity is that the Cauchy--Riemann equations are satisfied: 
\begin{align}\label{eq:chauchy_riemann}
\begin{aligned}
    &\partialder{u}{\xi} - \partialder{v}{\eta} = 0, \\
    &\partialder{u}{\eta} + \partialder{v}{\xi} = 0.
\end{aligned}
\end{align}
By computing the mixed derivatives of the terms in \eqref{eq:chauchy_riemann}, we obtain that $u$ and $v$ are harmonic functions as specified in the following equations:
\begin{align}\label{eq:har_components}
\begin{aligned}
	&\secondpartialder{u}{\xi} + \secondpartialder{u}{\eta} = 0, \\
	&\secondpartialder{v}{\xi} + \secondpartialder{v}{\eta} = 0.
\end{aligned}
\end{align}

\subsection{3D Laplace problem}\label{sec:3d_laplace}
In this section, we address the solution of the Laplace problem by deriving a formulation that is well suited for discretization using the neural networks presented in \Cref{sec:nn}. \\
Let $\Omega \subset \RR^3$ be a simply-connected open domain with boundary $\partial\Omega$, decomposed into two disjoint subsets $\partial\Omega_D$ and $\partial\Omega_N$ such that $\partial\Omega = \partial\Omega_D \cup \partial\Omega_N$. We denote by $\bfn$ the outward unit normal vector on $\partial\Omega$.
We consider the Laplace equation in $\Omega$ supplemented with boundary conditions:
\begin{subequations}\label{eq:lapl_sys}
\begin{align}
    &\laplacian V = 0,
    && \text{in } \Omega, \label{eq:lapl} \\
    &\boundarylapl(V) = 0,
    && \text{on } \partial\Omega,
\end{align}
\end{subequations}
where $\laplacian \cdot = \secondpartialder{\cdot}{x} +\secondpartialder{\cdot}{y} +\secondpartialder{\cdot}{z}$ is the Laplace operator; $V:\RR^3\to\RR$ denotes the scalar potential; $\boundarylapl(V)$ represents the prescribed boundary operator. Typical examples include
\begin{align*}
    &V = V_0,
    && \text{on } \partial\Omega_D, \\
    &\nabla V \cdot \bfn = g_0,
    && \text{on } \partial\Omega_N,
\end{align*}
where $V_0$ and $g_0$ are given functions. \\
A general representation formula for solutions of \eqref{eq:lapl} was derived by \citet{Whittaker1903}:
\begin{equation}\label{eq:general_sol_whittaker}
    V(x,y,z) = \int_0^{2\pi}
    f(\zeta(x,y,z,\theta),\theta)\, d\theta,
\end{equation}
where $\theta \in \RR$; the coordinate $\zeta:\RR^4\to\CC$ is a non-injective complex valued function defined as
\footnote{
In the original formulation found in \cite{Whittaker1903}, the potential is written as
$V = \int_0^{2\pi}f(r,\theta)\, d\theta$,
with $r(x,y,z, \theta) = z + i(\cos\theta\,x + \sin\theta\,y)$.
Noticing that, without loss of generality, $r = i\overline{\zeta}$, with $\zeta(x,y,z,\theta) = \cos\theta\,x + \sin\theta\,y + iz$, we may write
$f(r,\theta) = f(i\overline{\zeta},\theta)$.
Related generalizations are discussed in
\cite{Ketchum1928,Piltner1987,Plaksa2014}.
}
\begin{align}\label{eq:zeta_abc}
\begin{aligned}
    &\zeta = ax + by + cz, \\
    &(a,b,c) \text{ is a permutation of }
    (\cos(\theta), \sin(\theta), i).
\end{aligned}
\end{align}
Here, the function $f:\CC\times\RR\to\RR$ is harmonic in $x,y,z$ through $\zeta$
\begin{equation}\label{eq:f_har}
    \laplacian f(\zeta(x,y,z,\theta), \theta) = 0.
\end{equation}
We now proceed by determining how to replace $f$ with a holomorphic function. 
For the sake of exposition, we consider here $a=\cos(\theta)$, $b=\sin(\theta)$, $c=i$ and, using the identification $\CC\simeq\RR^2$, we write $\zeta = \xi + i\eta$ and $f = f(\xi, \eta)$. Applying the Laplace operator to $f$ yields
\begin{align*}
    \laplacian f =
    a^2 \secondpartialder{f}{\xi} + b^2 \secondpartialder{f}{\xi} +
    \secondpartialder{f}{\eta} =
    \secondpartialder{f}{\xi}
    +
    \secondpartialder{f}{\eta},
\end{align*}
since $a^2+b^2=1$. Then, using \eqref{eq:f_har}, we obtain
\begin{equation*}
    \secondpartialder{f}{\xi}
    +
    \secondpartialder{f}{\eta}
    =
    0,
\end{equation*}
which shows that $f$ is harmonic in the two-dimensional complex plane. 
As recalled in \eqref{eq:har_components}, both the real and imaginary parts of a holomorphic function are harmonic. Moreover, on simply connected domains, every harmonic function admits a representation as the real or imaginary part of a holomorphic function \cite{Ahlfors2021}.
Therefore, without loss of generality, we may write
\begin{equation*}
    f(\zeta, \theta) = \real(\tilde f(\zeta, \theta)),
\end{equation*}
where $\tilde f$ is holomorphic in $\zeta$. \\
It follows that the Laplace problem \eqref{eq:lapl_sys} can be reformulated in the following form.

\begin{problem}[Holomorphic Laplace formulation]\label{pb:lapl}
Find a holomorphic function $\tilde f(\zeta,\theta)$ such that
\begin{equation*}
    \boundarylapl(V)=0,
    \qquad
    \text{on } \partial\Omega,
\end{equation*}
where
\begin{equation}\label{eq:general_sol_whittaker_hol}
    V
    =
    \int_0^{2\pi}
    \real\bigl(\tilde f(\zeta,\theta)\bigr)\, d\theta.
\end{equation}
\end{problem}
We notice that the problem of finding a harmonic function has been transferred to finding a holomorphic function.
This observation is particularly relevant when combined with the following property of holomorphic functions: given $h$ and $g$ both holomorphic, their composition $h \circ g$ is holomorphic \cite{Stein2003}. This property enables us to approximate $\tilde f$ as a composition of ``simple'' holomorphic functions, $\tilde f \approx g\circ h \circ \ldots$, allowing the use of fully connected neural networks, presented in the following \Cref{sec:nn}, as they are a composition of affine and nonlinear functions.

\begin{remark}[Separated variable form is not complete]\label{sec:separated_variables}
    Assume, by contradiction, that $f$ in \eqref{eq:general_sol_whittaker} can be written in separated variable form $f(\zeta, \theta) = g(\zeta)h(\theta)$. Moreover, let consider a particular $V$ in the form  $V = V_1+V_2$. From \eqref{eq:general_sol_whittaker} and from the separated variable assumption, it follows that $V_1$ can be expressed as
    $V_1 = \int_0^{2\pi}g(\zeta_1)h_1(\theta)\, d\theta$ and similarly for $V_2$. By linearity, we obtain $V = \int g_1 h_1 + \int g_2 h_2 \neq \int gh$ as initially assumed. This remark shows that a generic non-linear dependence on both $\zeta$ and $\theta$ is required, justifying the need of a semi-holomorphic neural networks, see \Cref{sec:nn}. 
\end{remark}

\subsection{3D Linear elasticity problem}\label{sec:3d_elasticity}
In this section, we address the solution of the standard boundary value problem of linear elasticity. We consider a isotropic and homogeneous solid that occupies a simply-connected domain, denoted by $\Omega$, with Lipschitz continuous boundary. The mechanical behavior of the solid is linear elastic, characterized by shear modulus $\mu$ and Lamé first parameter $\lambda$. Neglecting body forces, the linear elastic problem is defined as \cite{Sadd2021}
\begin{subequations}\label{eq:mech}
\begin{align}
    &\mu \laplacian \bfu + (\lambda + \mu)\nabla(\divergence \bfu) = 0, && \text{in } \Omega,  \label{eq:Navier} \\
    &\boundarymech(\bfu,\bfsigma) = 0, && \text{on } \partial\Omega, 
\end{align}
\end{subequations}
where $\bfu$ is the unknown displacement vector and $\boundarymech$ is an operator defining proper boundary conditions in terms of $\bfu$ and of the stress tensor $\bfsigma$. This last quantity is related to $\bfu$ via the expression
\begin{equation}\label{eq:stress_tensor}
    \bfsigma = \mu\left(\nabla \bfu + (\nabla \bfu)^\top\right) + \lambda \text{Tr}(\nabla \bfu) \mathbf{I}
\end{equation}
where $\text{Tr}(\cdot)$ denotes the trace operator and $\mathbf{I}$ is the identity tensor. A typical choice for the boundary conditions is
\begin{align*}
    &\bfu = \bfu_0, &\on\partial\Omega_D, \\
    &\bfsigma\cdot \bfn = \bm{t}_0, &\on\partial\Omega_N, 
\end{align*}
where $\partial\Omega_D$ and $\partial\Omega_N$ are disjoint subsets of $\partial\Omega$ and $\bm{u}_0$ and $\bm{t}_0$ are prescribed displacement and traction vector. 
%
%
%
Let $\bfP = (P_x, P_y, P_z) : \RR^3 \to \RR^3$ and $B: \RR^3 \to \RR$ be such that $\laplacian \bfP = 0$ and $\laplacian B = 0$; these functions are referred to as \textit{potentials}.
Under the aforementioned hypothesis, the solution of \eqref{eq:mech} can be expressed using the Papkovich--Neuber representation \cite{Sadd2021}:
\begin{equation}
    \bfu = \frac{1}{2\mu} \big(\bfP - \nabla\left(B  + \kappa(P_x x + P_y y + P_z z) \right) \big) , \quad \text{in } \overline{\Omega}.  \label{eq:mech_PN_u} 
\end{equation}
where $\kappa = (\lambda+\mu)/(2\lambda + 4\mu)$. Hence, the system \eqref{eq:mech} can be replaced by the following:
\begin{subequations}\label{eq:mech_PN}
\begin{align}
    &\laplacian \bfP = 0, && \text{in } \Omega,  \label{eq:mech_PN_P} \\
    &\laplacian B = 0, && \text{in } \Omega,  \label{eq:mech_PN_B} \\
    &\boundarymech(\bfu(\bfP, B), \bfsigma(\bfP, B)) = 0, && \text{on } \partial\Omega.
\end{align}
\end{subequations}
The primary variables in \eqref{eq:mech_PN} are $\bfP$ and $B$. The number of unknown functions is thus four, the three components of $\bfP$ plus $B$, instead of the original three components of $\bfu$. As a consequence, the solution may not be unique, as discussed in e.g. \cite{TranCong1989}.
Without loss of generality, $\bfP$ and $B$ can be written using the integral representation \eqref{eq:general_sol_whittaker_hol}; therefore, we set 
\begin{align}\label{eq:P_B_Phi_chi}
\begin{aligned}
    &\bfP = \int_0^{2\pi} \real(\bfPhi(\zeta, \theta)) \; d\theta, \\
    &B = \int_0^{2\pi} \real(\chi(\zeta, \theta)) \; d\theta,
\end{aligned}
\end{align}
where $\bfPhi = (\phi_x, \phi_y, \phi_z) : (\CC\times [0, 2\pi))^3 \to \RR^3$ and $\chi: \CC\times [0, 2\pi) \to \RR$ are holomorphic functions in the first argument.
Since the governing PDEs in \eqref{eq:mech_PN} are identically satisfied by the properties of $\bfPhi$ and $\chi$, the linear elasticity problem can be formulated in the following way: 
\begin{problem}[Holomorphic linear elasticity]\label{pb:linear_elasticity}
Find $\bfPhi(\zeta, \theta)$ and $\chi(\zeta, \theta)$, holomorphic functions in $\zeta$, such that  
\begin{align*}
\begin{aligned}
    &\boundarymech(\bfu(\bfPhi, \chi), \bfsigma(\bfPhi, \chi)) = 0, && \text{on } \partial\Omega,
\end{aligned}
\end{align*}
where, by replacing \eqref{eq:P_B_Phi_chi} into \eqref{eq:mech_PN_u}, the displacement components can be expressed as
%
%
%
\begin{align}\label{eq:displ_phi_chi}
	\begin{aligned}
		&u_x = \frac{1}{2\mu}\int_0^{2\pi} \real((1-\kappa)\phi_x - a\chi' - a\kappa( 
		x \phi_x' + y \phi_y' + z \phi_z')), \\
		&u_y = \frac{1}{2\mu}\int_0^{2\pi} \real((1-\kappa)\phi_y - b\chi' - b\kappa(x \phi_x' + y \phi_y' + z \phi_z')), \\
		&u_z = \frac{1}{2\mu}\int_0^{2\pi} \real((1-\kappa)\phi_z - c\chi' - c\kappa(x \phi_x' + y \phi_y' + z \phi_z')), \\
	\end{aligned}
\end{align}
and the stress tensor, $\bfsigma$, is defined by the relationship \eqref{eq:stress_tensor}. The coefficients $a,b,c$ are defined in \eqref{eq:zeta_abc}. 
\end{problem}
%
%
In this work, the unknown functions $\bfPhi$ and $\chi$ are approximated by neural networks as described in the following \Cref{sec:nn}.

\section{Numerical strategies}\label{sec:nn}

In this section, we first review the notion of holomorphic neural network, originally introduced by \citet{Calafa2024}. We then present an architecture that preserves the holomorphicity of the network output with respect to a prescribed input argument, while relaxing the corresponding constraint for the remaining arguments. \\
We consider a $L$-layer neural network. Since the composition of holomorphic functions is itself holomorphic, the holomorphic property can be enforced in a neural network in a layer-wise manner. Consequently, we  focus on a layer $\ell \in 1, \ldots, L-1$, whose width is denoted by $n^\ell$. The relationship between the input $z^{\ell} \in \CC^{n^\ell}$ and the output $z^{\ell+1}\in\CC^{n^\ell+1}$ of the layer $\ell$ is given by
\begin{align*}
	z^{\ell+1} = \rho^\ell ( A^{\ell}(z^{\ell})),
\end{align*}
where $A^{\ell}(z^{\ell}) = W^{\ell} z^{\ell} + b^{\ell}$, the coefficients $W^{\ell}\in\CC^{n^{\ell+1}\times n^\ell}$ and $b^{\ell}\in\mathbb{R}^{n^{\ell+1}}$ are the trainable weights of the unknown affine transformation, and $\rho^{\ell}$ is a non-linear holomorphic activation function. 
It follows that the holomorphic multi-layered neural network can be written as the following composition
\begin{equation*}
	z^{L+1}  = (A^{L+1} \circ \rho^{L} \circ \ldots \circ \rho^{1} \circ A^{1})(\zeta),
\end{equation*}
where $\zeta \in \CC$ and $z^{(L+1)} \in \CC$ are the input and the output of the neural network, respectively. When non-polynomial activation functions that are holomorphic on the whole complex plane are used, the holomorphic neural network possesses the universal approximation property on simply connected domains \cite{Calafa2024,Geuchen2025}.

\begin{figure}[h]
    \centering
    \includegraphics[width=0.25\linewidth]{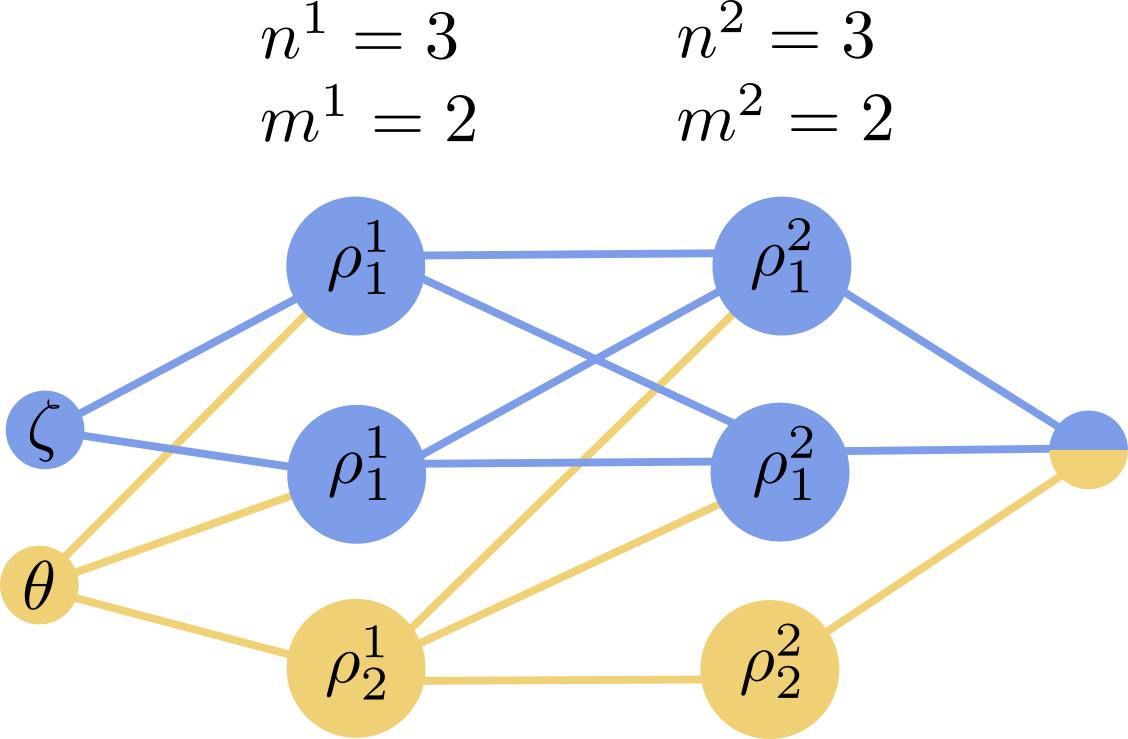}
    \caption{Example of a semi-holomorphic neural network with 2 hidden layers.}
    \label{fig:semi_hol_nn}
\end{figure}

We focus now in specifying the requirements for obtaining a network that approximates a function of the form $f(\zeta, \theta):\CC \times \RR \to \CC$ that is holomorphic with respect to the first argument $\zeta$. 
The use of a single network which considers general nonlinear interactions of the inputs $\zeta$ and $\theta$ is required, as discussed in \Cref{sec:separated_variables}.
To this end, let us consider a split value $m^\ell$ such that $1\leq m^\ell <n^\ell$, and let us consider $z^\ell \in \CC^{n^\ell - m^\ell}$, and $x^\ell \in \CC^{m^\ell}$. Then, we consider a neural network whose layer $\ell$, with input $[{z^\ell}^\top, {x^\ell}^\top]^\top \in \CC^{n_{\ell}}$ and output $[{z^{\ell+1}}^\top, {x^{\ell+1}}^\top]^\top \in \CC^{n^{\ell+1}}$, is constructed in the following form:
\begin{align*}
	\begin{bmatrix}
		z^{\ell+1} \\
		x^{\ell+1}
	\end{bmatrix} = 
	\begin{bmatrix}
		\rho_1^{\ell} \\
		\rho_2^{\ell} 
	\end{bmatrix}\left(
	\begin{bmatrix}
		W_{11}^{\ell} &W_{12}^{\ell} \\
		0 & W_{22}^{\ell}
	\end{bmatrix} 
	\begin{bmatrix}
		z^{\ell} \\
		x^{\ell}
	\end{bmatrix} 
	+ 
	\begin{bmatrix}
		b_1^{\ell} \\ 
		b_2^{\ell}    
	\end{bmatrix}
	\right),
\end{align*}
where $\rho_1^{\ell}$ is a holomorphic activation function, $\rho_2^{\ell}$ is a non-holomorphic activation function, 
$W_{11}^\ell\in\RR^{m^{\ell+1}\times m^\ell}$,
$W_{12}^\ell\in\RR^{m^{\ell+1}\times p^\ell}$,
$W_{21}^\ell\in\RR^{p^{\ell+1}\times p^\ell}$,
where $p^\ell$ is such that $m^\ell + p^\ell = n^\ell$.
We call $\varsigma^{\ell}  = [\rho_1^{\ell}, \rho_2^{\ell}]^\top$, and $\displaystyle H^{\ell} (\cdot) = \begin{bmatrix}
	W_{11}^{\ell} &W_{12}^{\ell} \\
	0 & W_{22}^{\ell}
\end{bmatrix} (\cdot)
+ 
\begin{bmatrix}
	b_1^{\ell} \\ 
	b_2^{\ell}    
\end{bmatrix}$.
It follows that the multi-layered neural network can be written as the following composition:
\begin{equation}\label{eq:semi_hol_nn}
	z^{L+1} = (H^{L+1} \circ \varsigma^{L} \circ \ldots \circ \varsigma^1 \circ H^1) ([\zeta^\top, \theta^\top]^\top).
\end{equation}
It can be readily verified that the architecture defined by \eqref{eq:semi_hol_nn}, which is schematically shown in \Cref{fig:semi_hol_nn}, has the holomorphic property in $\zeta$ and not necessarily in $\theta$, as desired. 
\textit{Semi-holomorphic} neural networks constructed in this way are used to approximate the function $\tilde f$ in \Cref{pb:lapl} and the potentials $\phi_x, \phi_y, \phi_z$ and $\chi$ in \Cref{pb:linear_elasticity}. We denote these neural networks as $N_V$, $N_{\phi_x}, N_{\phi_y}, N_{\phi_z}$, and $N_\chi$ respectively.

\paragraph{Training.}
To simplify the notation, we denote by $\mathcal{W}$ the union of all trainable parameters, and by $|\mathcal{W}|$ their total number.
The trainable parameters are initialized using the strategy proposed in \cite{Calafa2024} and their optimal values, denoted by $\mathcal{W}_\text{opt}$, are obtained by solving the following optimization problem:
\begin{align*}
	\mathcal{W}_\text{opt} = \arg \min_{\mathcal{W} \in \Xi} \loss(\mathcal{W}),
\end{align*}
where $\Xi$ denotes the admissible search space for $\mathcal{W}$ and $\loss$ is the loss function. 
The loss function must measure a suitable error on the boundary conditions. Considering the $\ell_2$-norm, we can write
\begin{align*}    
    & \loss =  \frac{1}{n_\text{train}}\sum_{i=1}^{n_\text{train}} \loss_p(\bfx_i), \\
    &\loss_p(\bfx) =\|\mathcal{B}_j|_{\bfx}\|^2, \quad j = {\scriptstyle{L, M}}, \ \bfx \in \partial\Omega, 
\end{align*}
where $\loss_p$ is the per-sample loss, $\bfx_i$ are the training point, $n_\text{train}$ the total number of trainig points. 
Depending on the specific problem at hand, examples of per-sample losses are: 
$\loss_p = \|\overline{V}-V_0\|^2_{\partial\Omega_D} + \alpha_1\|\nabla \overline{V}\cdot \bfn - g_0\|^2_{\partial\Omega_N}$ for \Cref{pb:lapl}, and 
$\loss_p = \|\overline{\bfu} - \bfu_0\|^2_{\partial\Omega_D} + \alpha_2\|\overline{\bfsigma}\cdot\bfn - \bm{t}_0\|^2_{\partial\Omega_N}$ for \Cref{pb:linear_elasticity},
where $\overline{\cdot}$ are the quantities derived from the neural networks and $\alpha_1,\alpha_2 \in \mathbb{R}_+$ are weighting factors. \\
We highlight that the holomorphic property holds for any $\mathcal{W} \in \Xi$, which means that this is not a property obtained after training convergence, but rather an intrinsic property of the presented neural network.
\begin{remark}[Comparison with PINN]
    In the standard PINN approach \cite{Raissi2019}, the loss function is defined as the sum of the errors associated with the residual of the governing PDEs and the corresponding boundary conditions, namely $\loss = \alpha \|\mathcal{D}\|^2 +  \|\mathcal{B}\|^2$, where $\mathcal{D}$ denotes the differential operator and $\alpha\in\mathbb{R}_+$. In the present approach, the term $\|\mathcal{D}\|^2$ vanishes identically by construction, therefore the loss function reduces to $\loss = \|\mathcal{B}\|^2$. 
    Moreover, since the operator $\mathcal{D}$ involves derivatives of higher order than those appearing in $\mathcal{B}$ \cite{Reddy1998}, the proposed training strategy avoids the computation of derivatives of order higher than that strictly required in the computation of the boundary conditions. \\
    Moreover, let $d$ be the characteristic inter-point distance, assumed approximately constant throughout the domain. Also, let $n_{\text{train},\Omega}$ and $n_{\text{train}, \partial\Omega}$ denote the total number of training points in $\Omega$ and $\partial\Omega$, respectively. These quantities scale as $n_{\text{train}, \Omega} \sim \frac{|\Omega|}{d^3}$ and $n_{\text{train}, \partial\Omega} \sim \frac{|\partial\Omega|}{d^2}$. In the standard PINN approach, a total of $n_{\text{train},\Omega} + n_{\text{train},\partial\Omega}$ points must be considered, whereas in the proposed method only $n_{\text{train},\partial\Omega}$ is required. As a consequence, the proposed method exhibits a more favourable scaling in $d$ compared to standard PINNs.
\end{remark}
\paragraph{Integration.}
The integrals in \eqref{eq:general_sol_whittaker} and \eqref{eq:displ_phi_chi} are approximated through a quadrature rule of the form $\int_0^{2\pi} f(\theta)\; d\theta \approx \sum_i w_i f(\theta_i)$, where $w_i$ are proper scalar weights. We denote by $\mathcal{Q}(\cdot)$ such approximation of $\int_0^{2\pi}\cdot \, d\theta$. 
Moreover, we denote by $\overline{V}:\RR^3\to\RR$ the quantity
\begin{align}\label{eq:V_discrete}
	\overline{V} = \mathcal{Q}(\real(N_V)),
\end{align}
which serves as the neural network approximation of $V$ via a discrete integral. Similarly, we denote by $\overline{\bfu}$ and $\overline{\bfsigma}$ the approximated displacement and stress.  
\begin{remark}[Null error on PDE]\label{th:null_residual}
    We observe that the error introduced by the discretization of the integral does not affect the satisfaction of the governing PDEs. Indeed, $\laplacian \overline{V} = \laplacian\mathcal{Q}(\real(N_V)) = \mathcal{Q}(\laplacian\real(N_V)) = 0$. Similarly for $\bfP$ and $B$. This implies that \eqref{eq:lapl}, \eqref{eq:mech_PN_P}, \eqref{eq:mech_PN_B} are satisfied exactly, independently of the accuracy of the integration method and the training of the network. Therefore, the only source of error introduced by the proposed method is the approximation of the boundary conditions.
\end{remark}
\begin{remark}[Max error on boundary]\label{th:error_boundary} For \Cref{pb:lapl}, the discretization error in $\Omega$ is bounded by the discretization error on $\partial\Omega$. Indeed, let $\Delta V$ denote the discrepancy between the exact $V$ and the approximated one, so that $\overline{V} = V + \Delta V$. We have $\laplacian \overline{V} = \laplacian (V+\Delta V) = \laplacian \Delta V = 0$. The maximum principle \cite{Evans2022} applies to $\Delta V$, meaning that the maximum value of the discrepancy between the exact and approximated solution is attained on the boundary. It follows that the error $\|\Delta V\|$  has its maximum on the boundary.  
\end{remark}

\section{Test cases}\label{sec:tests}

We present in this section a series of test cases to assess the effectiveness and limitations of the proposed method. Two cases, in \Cref{sec:test_lapl}, refer to the Laplace problem, presented in \Cref{sec:3d_laplace}, while the remaining two, in \Cref{sec:test_mech}, deal with the linear elastic problem, presented in \Cref{sec:3d_elasticity}.

To train neural networks, points are required on the boundary of the domain geometry. These can be determined using any suitable meshing strategy. In this work, for simplicity, such  boundary points are computed as the centroids of a simplex grid. Nevertheless, we emphasize that the method is meshless, and the use of a mesh is solely the result of a practically convenient procedure for defining a proper representation of the geometry and its boundary points. A similar remark applies to the figures presented in this section: the meshed appearance of some geometries is merely the outcome of a convenient procedure for illustrating three-dimensional structures. \\
The points on the boundary are partitioned into a training set and a test set: the former contains $n_\text{train}$ points, while the latter contains $n_\text{test}$ points. For each case, the train--test split is defined by $n_\text{train} = 0.9(n_\text{train} + n_\text{test})$. \\
When multiple networks are used, as in the elasticity-related test cases, the trainable parameters of the four networks are combined into a single set, $\mathcal{W}_\text{tot} = \mathcal{W}_{\phi_x} \cup \mathcal{W}_{\phi_y} \cup \mathcal{W}_{\phi_z} \cup \mathcal{W}_{\chi}$, where $\mathcal{W}_{\phi_x}, \ldots$ denote the parameters of $N_{\phi_x}, \ldots$. The set $\mathcal{W}_\text{tot}$ thus becomes the collection of trainable weights to be optimized concurrently using an appropriate strategy. \\
In each test case, the neural networks are trained using the ADAM algorithm, with the parameters suggested in \cite{Kingma2014}. Preliminary studies reveal no meaningful improvements from adopting different strategies. \\
To assess the accuracy of the proposed method, we make use of the following norms, defined analogously for scalars $a$,  vectors with components $b_i$, and matrices with components $c_{ij}$:
\begin{align*}
	\|a\|_{\Omega} \coloneq \sqrt{\int_\Omega a^2}, \quad \|b\|_{\Omega} \coloneq \sqrt{\int_\Omega \sum_i b_i^2}, \quad 
	\|c\|_{\Omega} \coloneq \sqrt{\int_\Omega \sum_{ij} c_{ij}^2}.
\end{align*}
The proper evaluation and comparison of the method with existing approaches is challenging due to the lack of appropriate metrics for comparing different methods, see, e.g., \cite{McGreivy2024}. In this work, we restrict ourselves to assessing the accuracy of the method against a reference solution, either analytical or obtained with FEM on a sufficiently fine grid, and to highlighting its main properties.\\ 
%
To ease the readability, we restrict the discussion to the main features of each test case. Further implementation details are provided in \Cref{sec:hyperparameters_details}. \\
All numerical values of the physical quantities reported in the following test cases are considered dimensionless.

\subsection{Laplace problem}\label{sec:test_lapl}
In this section, we show two test cases regarding the Laplace problem.
To assess the accuracy of the results, we introduce the following relative errors:
\begin{subequations}\label{eq:errs_lapl}
\begin{align}
	&\varepsilon_V = \frac{\sqrt{|\Omega|}\,|\overline{V}-V_\text{ref}|}{\|V_\text{ref}\|_{\Omega}}, \label{eq:errs_lapl_local}\\
    &\overline{\varepsilon}_V = \frac{\|\overline{V} - V_\text{ref}\|_\Omega}{\|V_\text{ref}\|_\Omega}, \label{eq:errs_lapl_global} \\
    &\varepsilon_{\nabla V} = \frac{\sqrt{|\Omega|}\,|\nabla\overline{V}-\nabla V_\text{ref}|}{\|\nabla V_\text{ref}\|_{\Omega}}, \label{eq:errs_grad_lapl_local}\\
    &\overline{\varepsilon}_{\nabla V} = \frac{\|\nabla\overline{V} - \nabla V_\text{ref}\|_\Omega}{\|\nabla V_\text{ref}\|_\Omega}, \label{eq:errs_grad_lapl_global}
\end{align}
\end{subequations}
where the reference solution, $V_\text{ref}$, is obtained either analytically or via FEM.

\subsubsection{Manufactured solution on complex geometry}\label{sec:test_lapl_complex}
\begin{figure}[H]
    \centering
    \hspace*{-27mm}\raisebox{-2mm}{\includegraphics[width=0.35\linewidth]{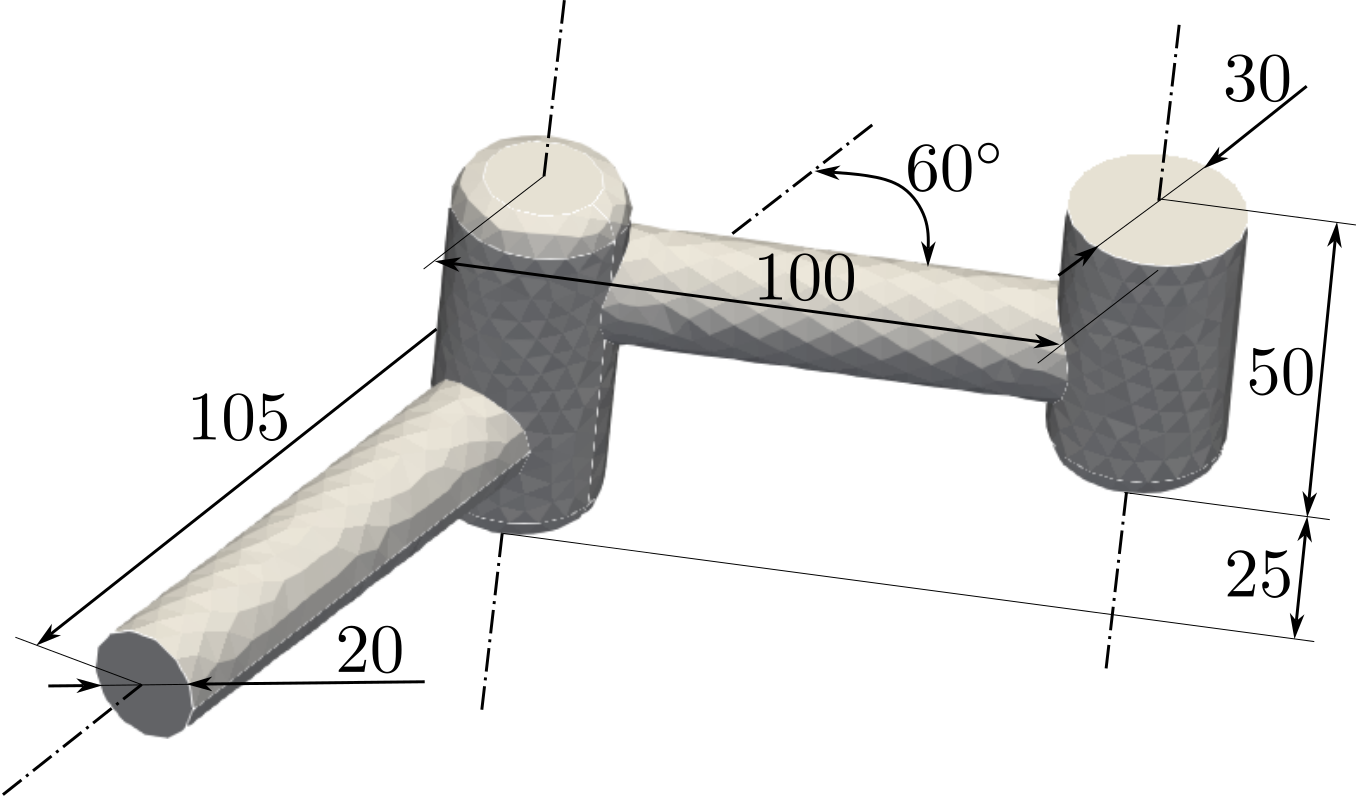}}
    \includegraphics[width=0.3\linewidth]{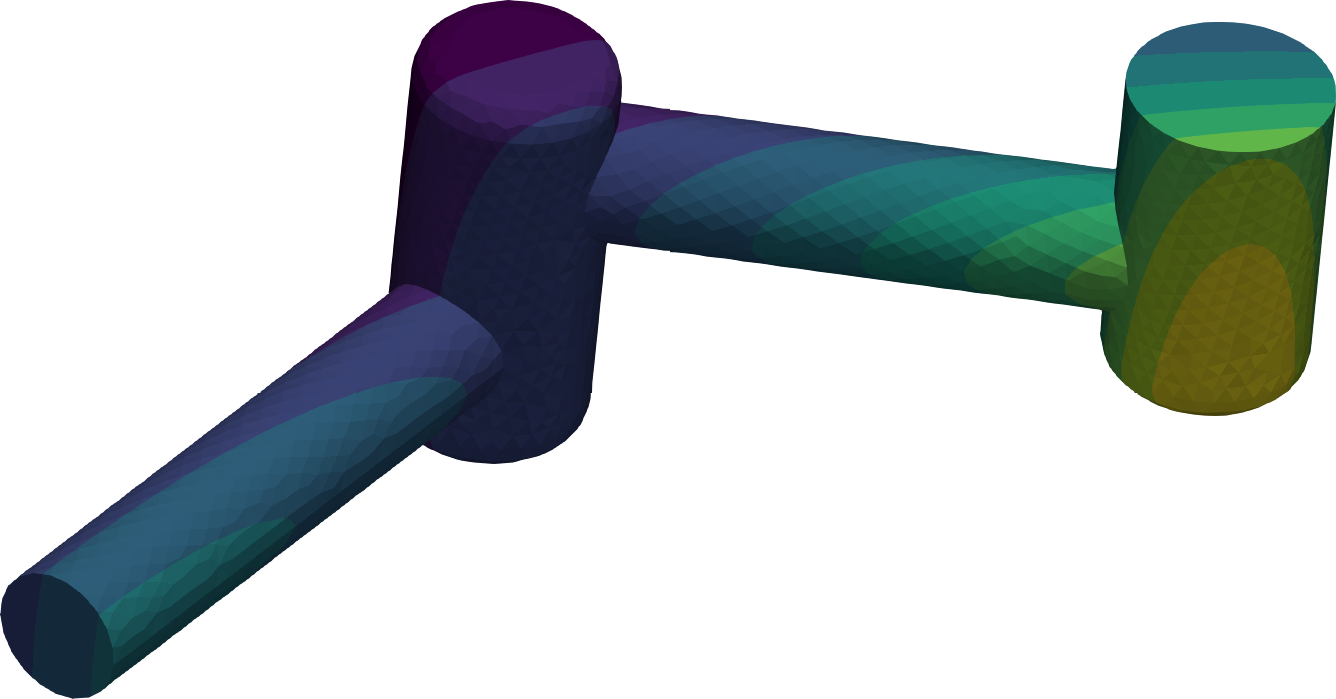} 
    \includegraphics[width=0.136\linewidth]{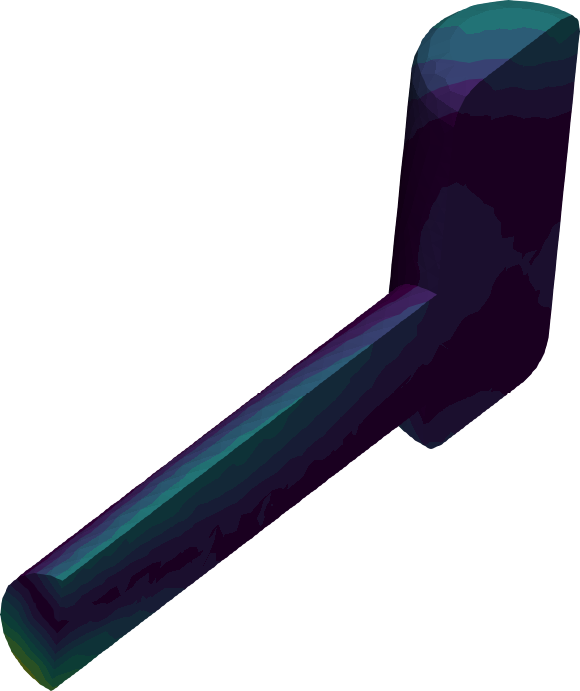}
    \vspace*{6mm} \\
    %
    \includegraphics[width=0.3\linewidth]{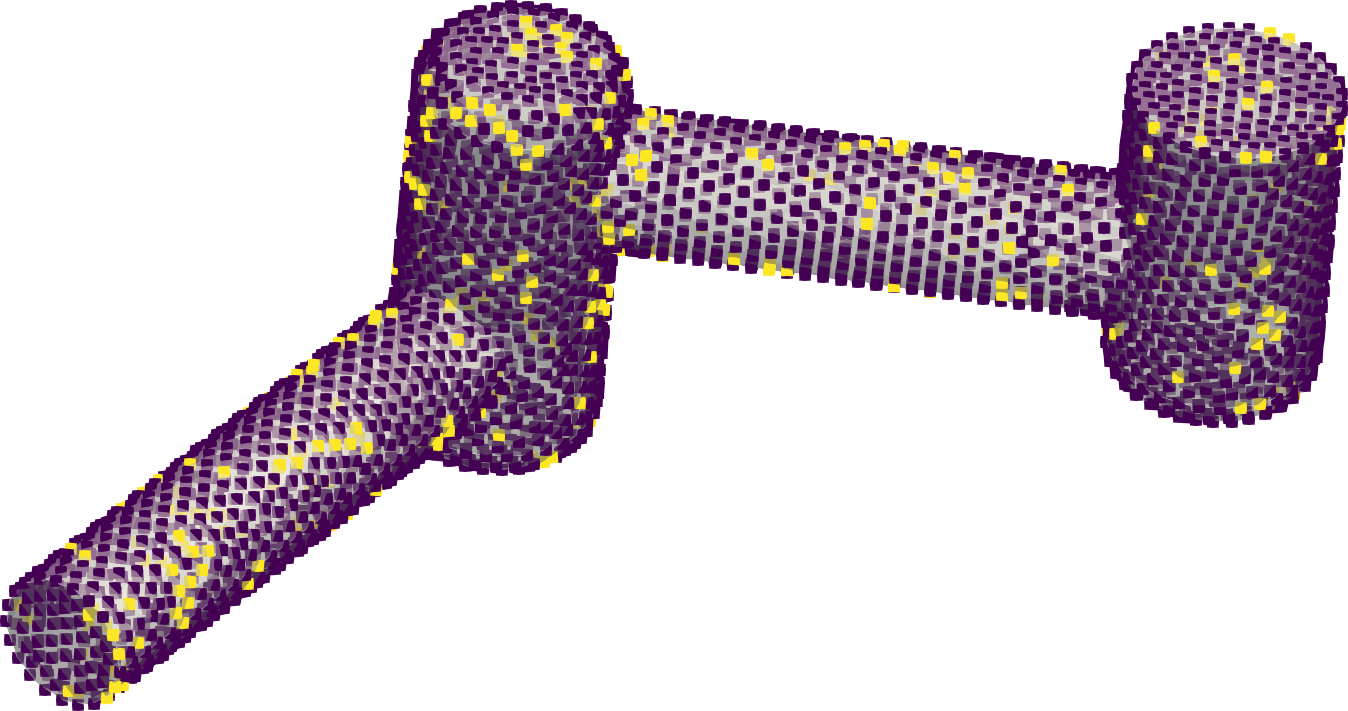}
    \includegraphics[width=0.3\linewidth]{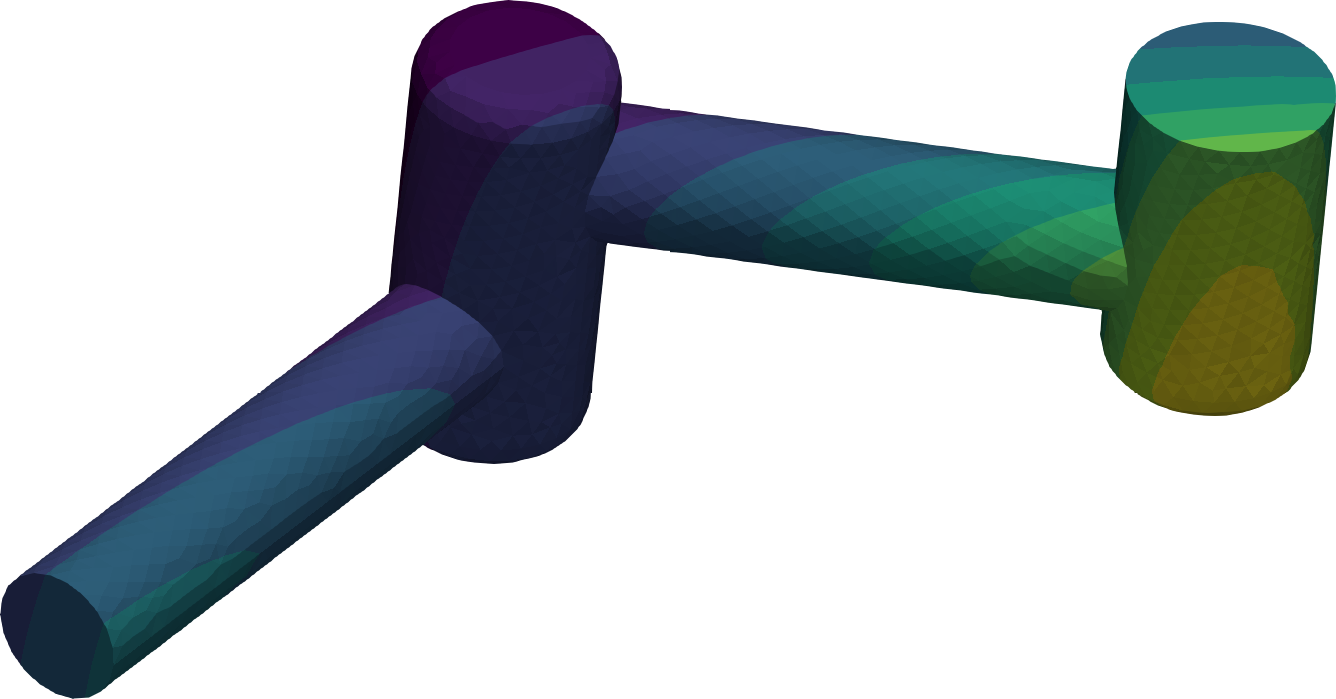} 
    \includegraphics[width=0.3\linewidth]{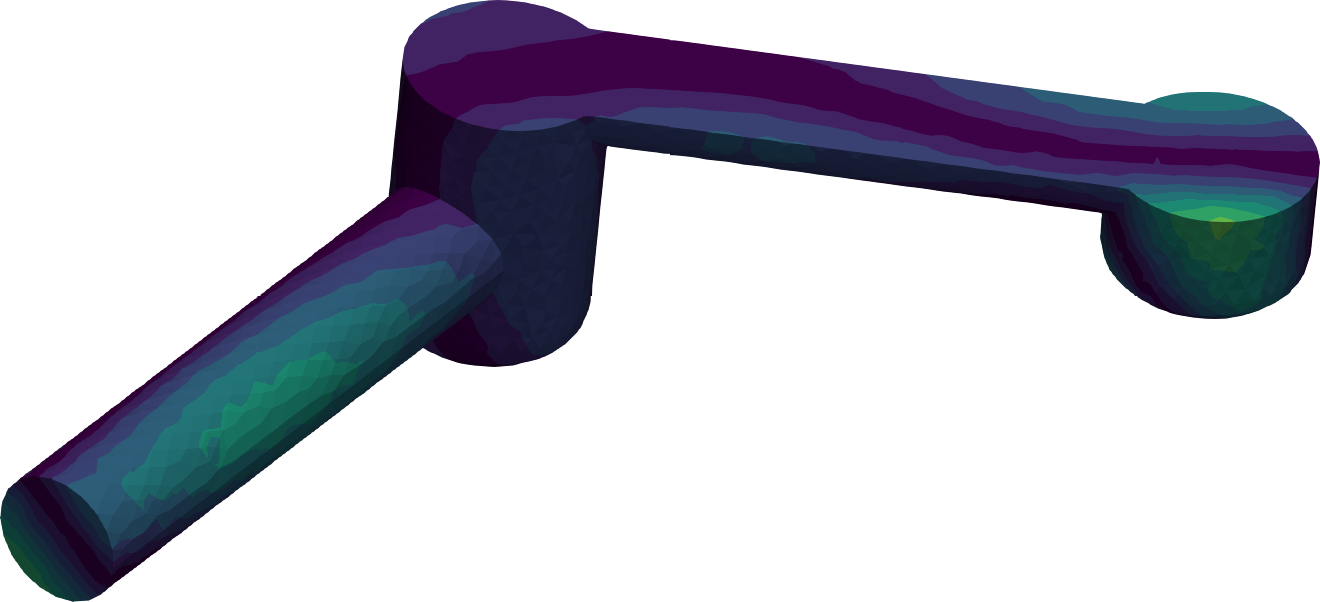} \\
    \subfloat[Domain and dataset points]{\includegraphics[width=0.2\linewidth]{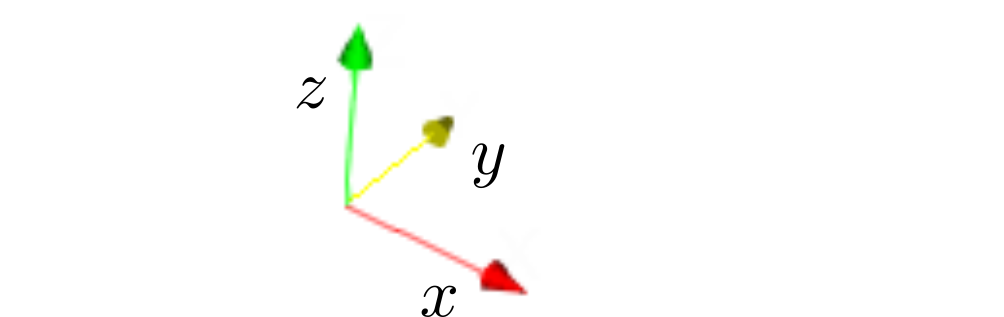}}\hspace{10mm}
    \subfloat[$V_\text{ref}$ and $\overline{V}$]{\includegraphics[width=0.25\linewidth]{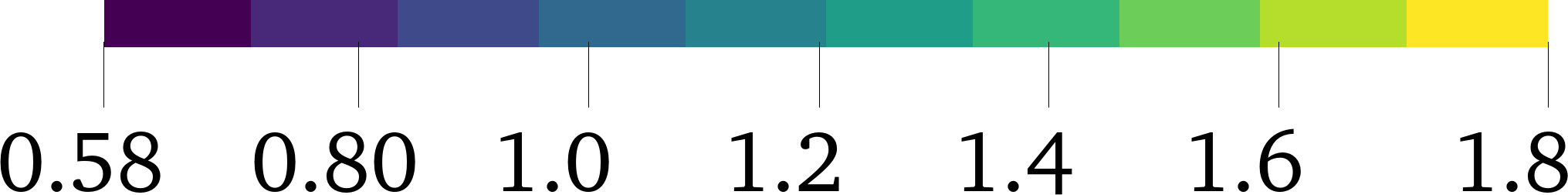}}\hspace{10mm}
    \subfloat[Relative error, $\eps_V$]{\includegraphics[width=0.25\linewidth]{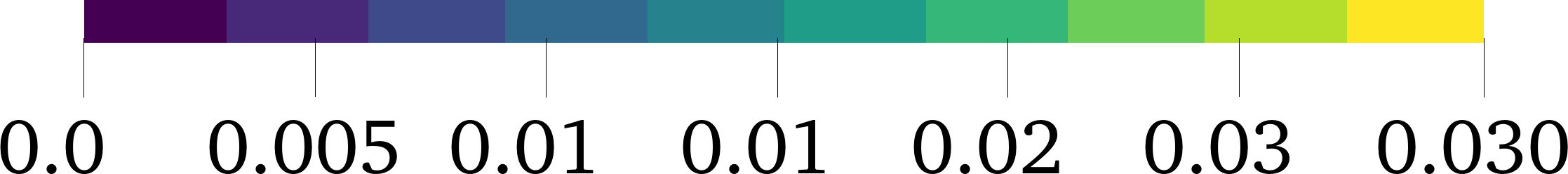}}
    \caption{Domain and numerical results. (a) Domain and dimensions; the training points (dark color) and test points (light color) are displayed below.
    	(b) Analytical harmonic function $V_\text{ref}$ on top and its neural network approximation below, $\overline{V}$. The two solutions are graphically similar.
    	(c) Relative error \eqref{eq:errs_lapl_local} on both the boundary and interior points, visualized through two slices of the geometry. The error remains low across the entire domain. }
    \label{fig:test_lapl_supp}
\end{figure}

This first test case involves a domain with complex geometry, which is illustrated in \Cref{fig:test_lapl_supp}, and the use of an analytical reference solution:
\begin{equation*}
	V_\text{ref} = \exp(\sqrt{13}\, x) \cos(2 y) \cos(3 z), \quad \text{in } \overline{\Omega}.
\end{equation*}
We define the boundary operator as
\begin{equation*}
	\boundarylapl(V) \coloneq V - V_\text{ref}.
\end{equation*}
Dirichlet boundary conditions are then imposed as
\begin{equation*}
	\boundarylapl(V) = 0 \quad \on \partial\Omega.
\end{equation*}
This problem may be interpreted, for example, as a heat conduction process in a mechanical device with prescribed temperature on the boundary. We employ a semi-holomorphic neural network whose layer sizes are $(2, 16, 16, 16, 16, 1)$, where the first and last entries denote the input and output dimensions, respectively. The value of $m^\ell$, which determines the split between the holomorphic and non-holomorphic components, is set to $m^\ell = n^\ell /2, \  \ell=1,\ldots,5$. The exponential activation function is used at each layer of the holomorphic part, while a complex PReLU function, that we denote by cPReLU, is used for the non-holomorphic part. The latter is defined as:
\begin{align*}
	\text{cPReLU}(z) = \text{PReLU}(\real(z)) + i\text{PReLU}(\imag(z)), &&\text{where } \ \text{PReLU}(x) = \begin{cases}
		x  &\text{if } x > 0,\\
		0.1 \, x  &\text{if } x\leq 0.
	\end{cases}
\end{align*}
In this case, we set $\zeta = \sin(\theta)y+\cos(\theta)z + i x$ and we adopt the composite midpoint rule \cite{Burden2011} with $64$ intervals equally spaced in $[0, 2\pi)$ to discretize the integral in \eqref{eq:general_sol_whittaker_hol}. 
The data used to train the neural network consist of a discrete set of coordinate points on the boundary, shown in panel (a) of \Cref{fig:test_lapl_supp}, with $n_\text{train} = \numprint{2232}$ and $n_\text{test} = 248$. \\
The network is trained by minimizing the following loss function, $\loss$, associated with Dirichlet boundary conditions:
\begin{align*}
\begin{aligned} 
	&\loss = \frac{1}{n_\text{train}}\sum_{i=1}^{n_\text{train}}\loss_p(\bfx_i), \\      
    &\loss_p(\bfx) = \left\| \boundarylapl(\overline{V}(\bfx)) \right\|^2, & \bfx \in \partial\Omega,
\end{aligned}
\end{align*}
where $\overline{V}$ is defined in \eqref{eq:V_discrete}.
Training lasts $\numprint{3000}$ epochs, until the loss function reaches a fairly stationary low value. Further details can be found in \Cref{sec:appendix_support}. \\
The results are reported in \Cref{fig:test_lapl_supp}. Panel (b) shows good agreement between the exact solution and its approximation.
This qualitative observation is quantitatively confirmed by the relative error displayed in panel (c) of the same figure. Indeed, its local maximum value is about $3\%$. Moreover, we observe that the error peaks on the boundary, as noted in \Cref{th:error_boundary}.
To assess the robustness of the results, we report in \Cref{tab:ablation} a sensitivity study in which the number of quadrature intervals and the layer width are varied independently, while all remaining hyperparameters, e.g., the number of layers and the training settings, are kept fixed. We observe that the mean relative error $\overline{\eps}_{\nabla V}$ remains low across all tested configurations, increasing only slightly for the coarser settings. This suggests that the method is not strongly sensitive to these particular choices of hyperparameters within the explored ranges.
\begin{table}[h]
\centering
\begin{tabular}{lllll}
    \# Intervals & Layer widths & $|\mathcal{W}|$ & $\overline{\eps}_V$ & \\
    \cline{1-4}
    16 & 16 & 881 &  $9.67\times10^{-3}$ & \\
    32 & 16 & 881 &  $8.91\times10^{-3}$ & \\
    64$^\dagger$ & 16 & 881 & $8.36\times10^{-3}$ & \\
    128 & 16 & 881 & $7.76\times10^{-3}$ & \\
    64 & 4 &  77 & $1.11\times10^{-2}$ & \\
    64 & 8 &  249 & $9.00\times10^{-3}$ & \\
    64 &  32 & \numprint{3297} & $8.98\times10^{-3}$ &\\
    $^\dagger$Shown in \Cref{fig:test_lapl_supp}.
\end{tabular}
\caption{Sensitivity study: mean relative error $\overline{\eps}_V$ for varying numbers of 
quadrature intervals and layer widths. The total number of trainable parameters, $|\mathcal{W}|$, is reported for reference. Results remain moderately stable across all tested configurations.}
\label{tab:ablation}
\end{table}

\subsubsection{Neumann problem on a quarter torus}\label{sec:test_lapl_fluid}
\begin{figure}[H]
    \centering
    \subfloat[Domain and b.c.]{\includegraphics[width=0.40\linewidth]{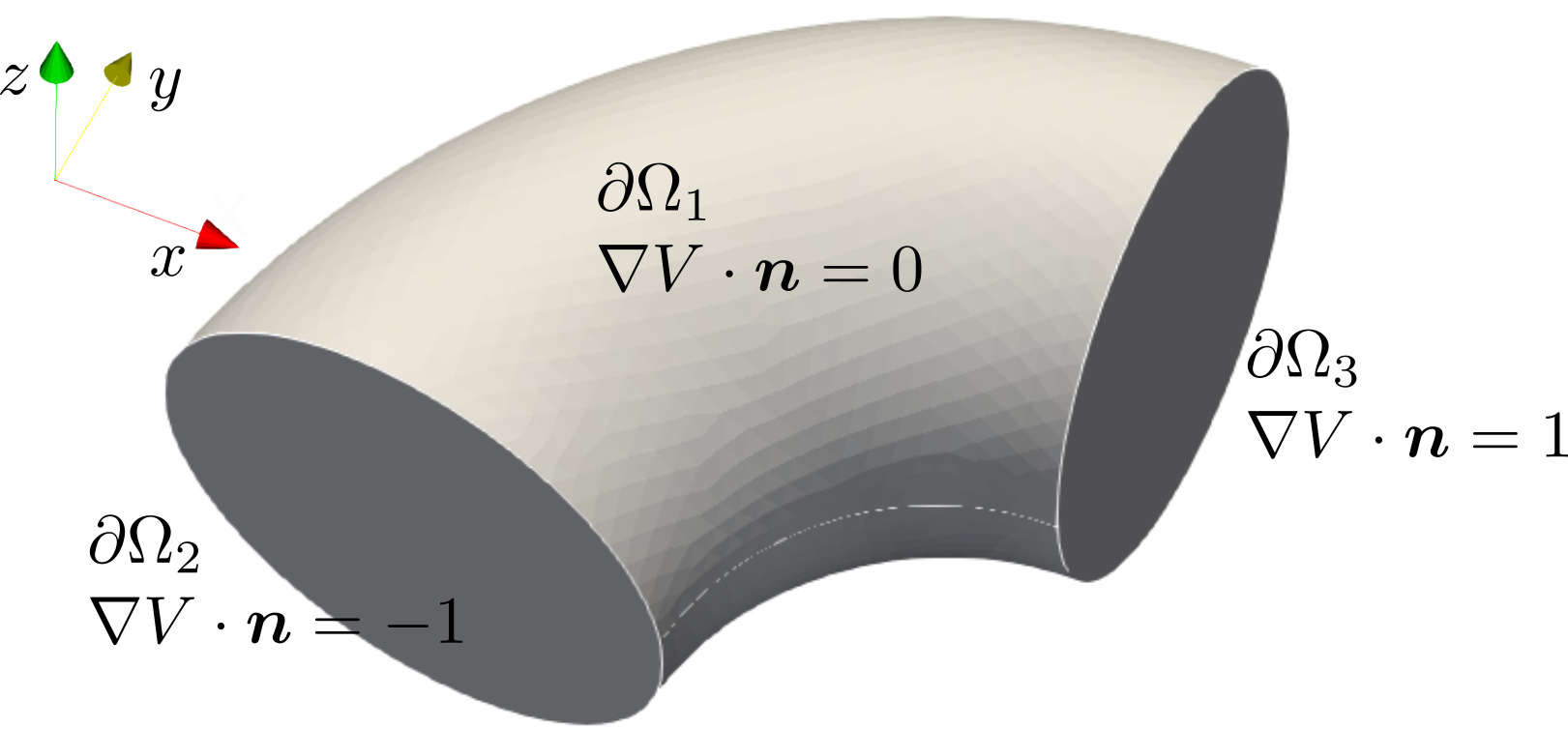}}\hspace{10mm}
    \subfloat[Velocity - FEM]{\includegraphics[width=0.28\linewidth]{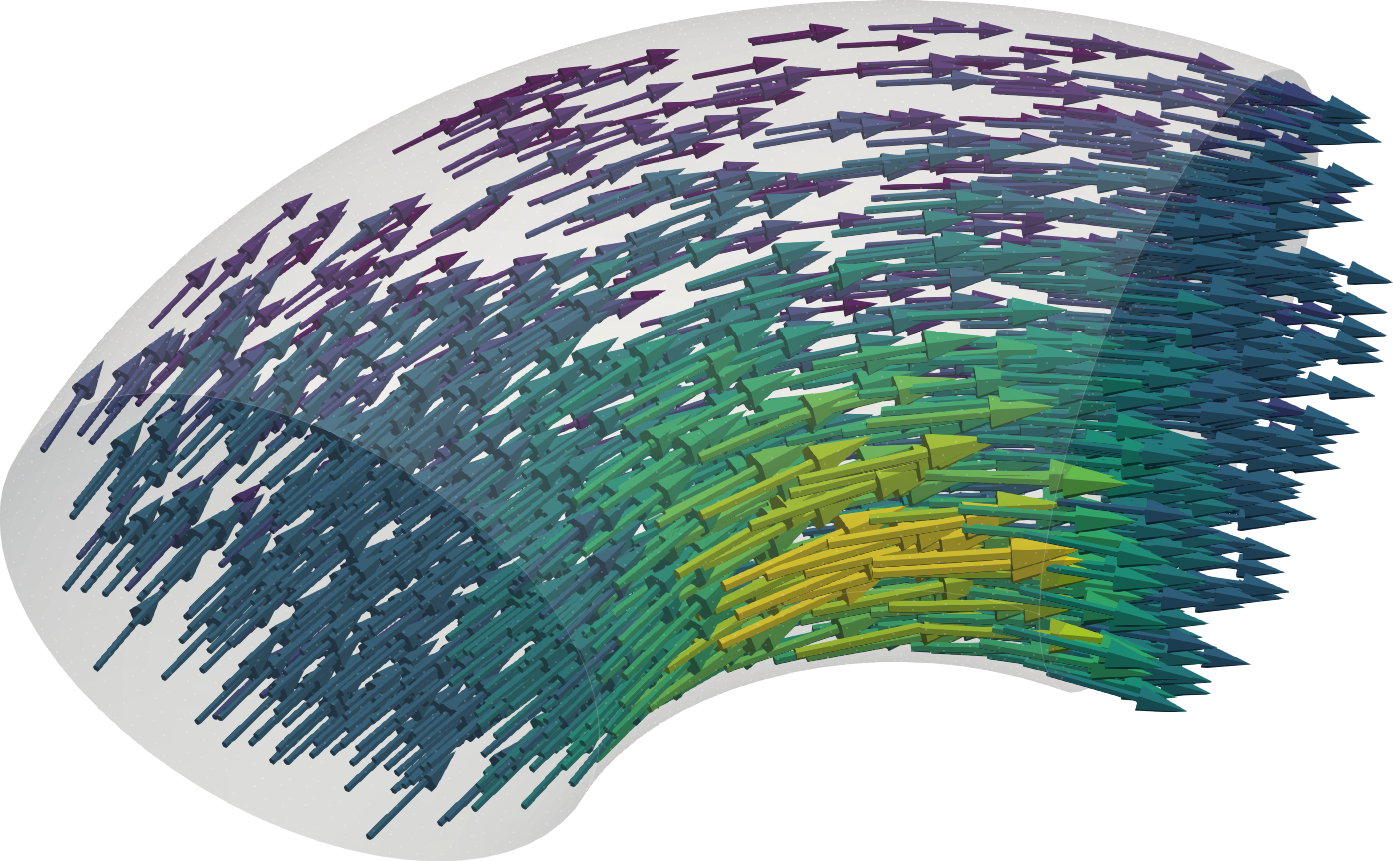}\hspace{3mm} \includegraphics[width=0.1\linewidth]{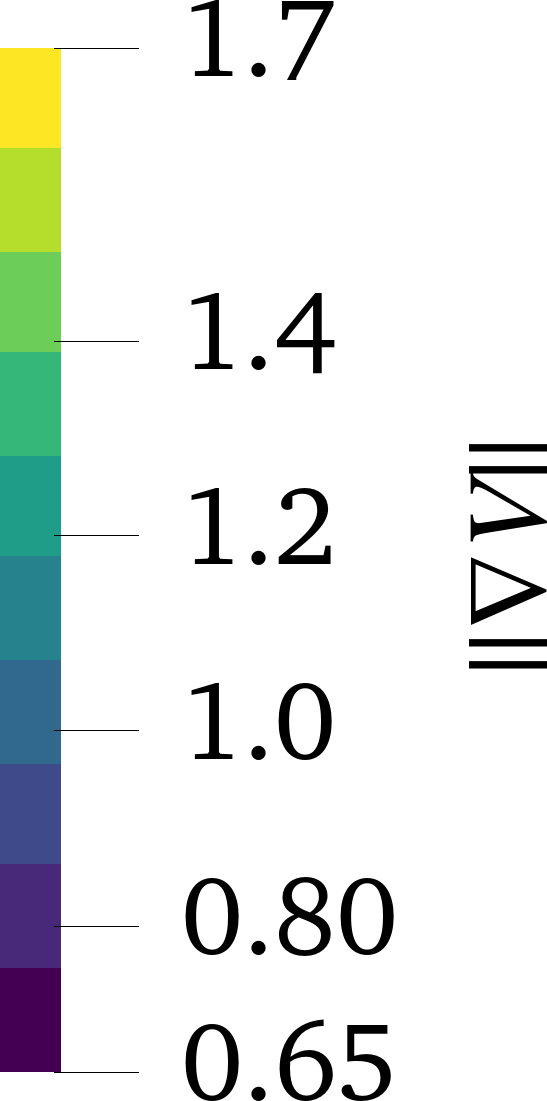}} \\
    %
    \subfloat[Velocity - PINN]{\includegraphics[width=0.28\linewidth]{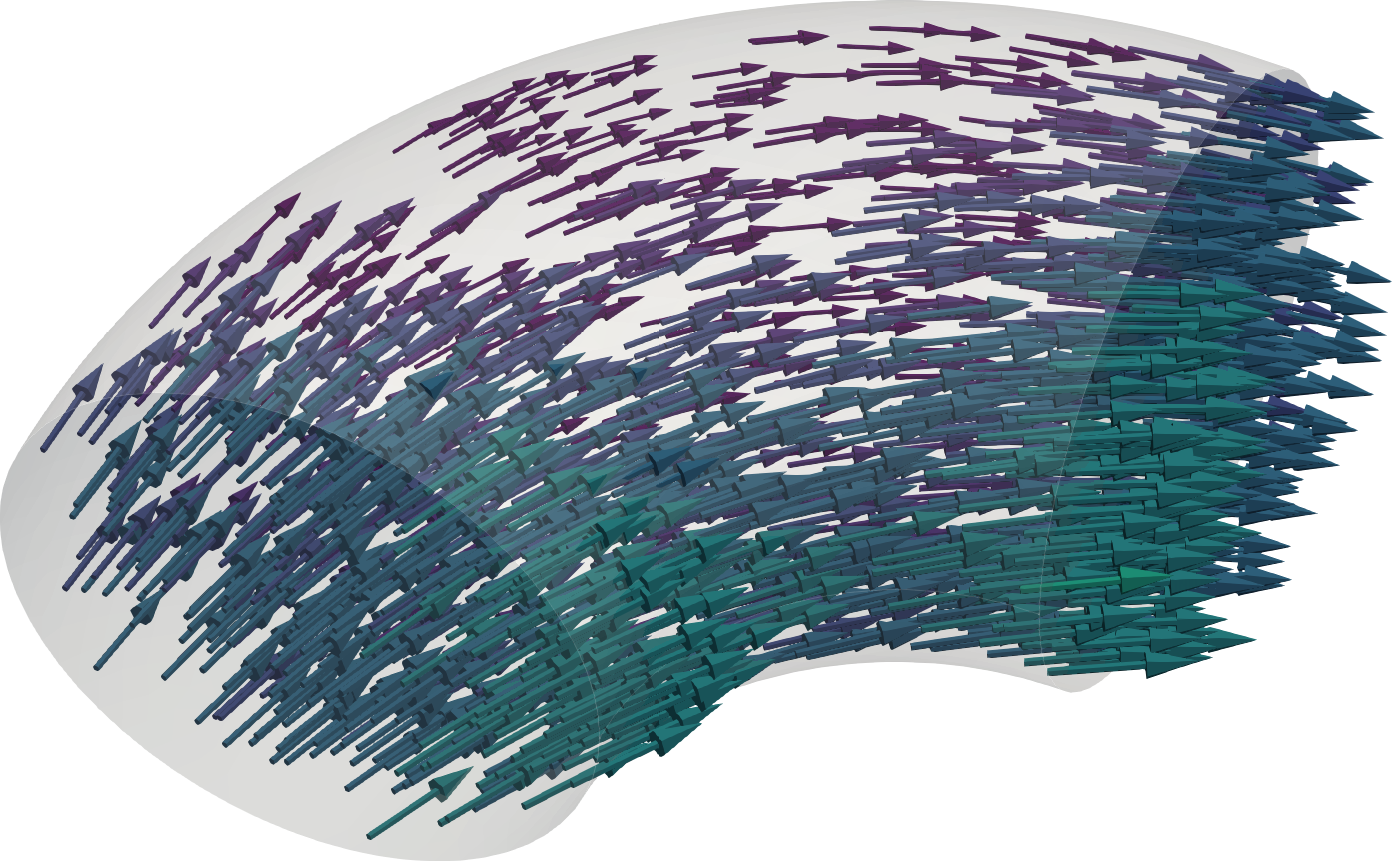}}\hspace{20mm} 
    \subfloat[Velocity - HOL]{\includegraphics[width=0.28\linewidth]{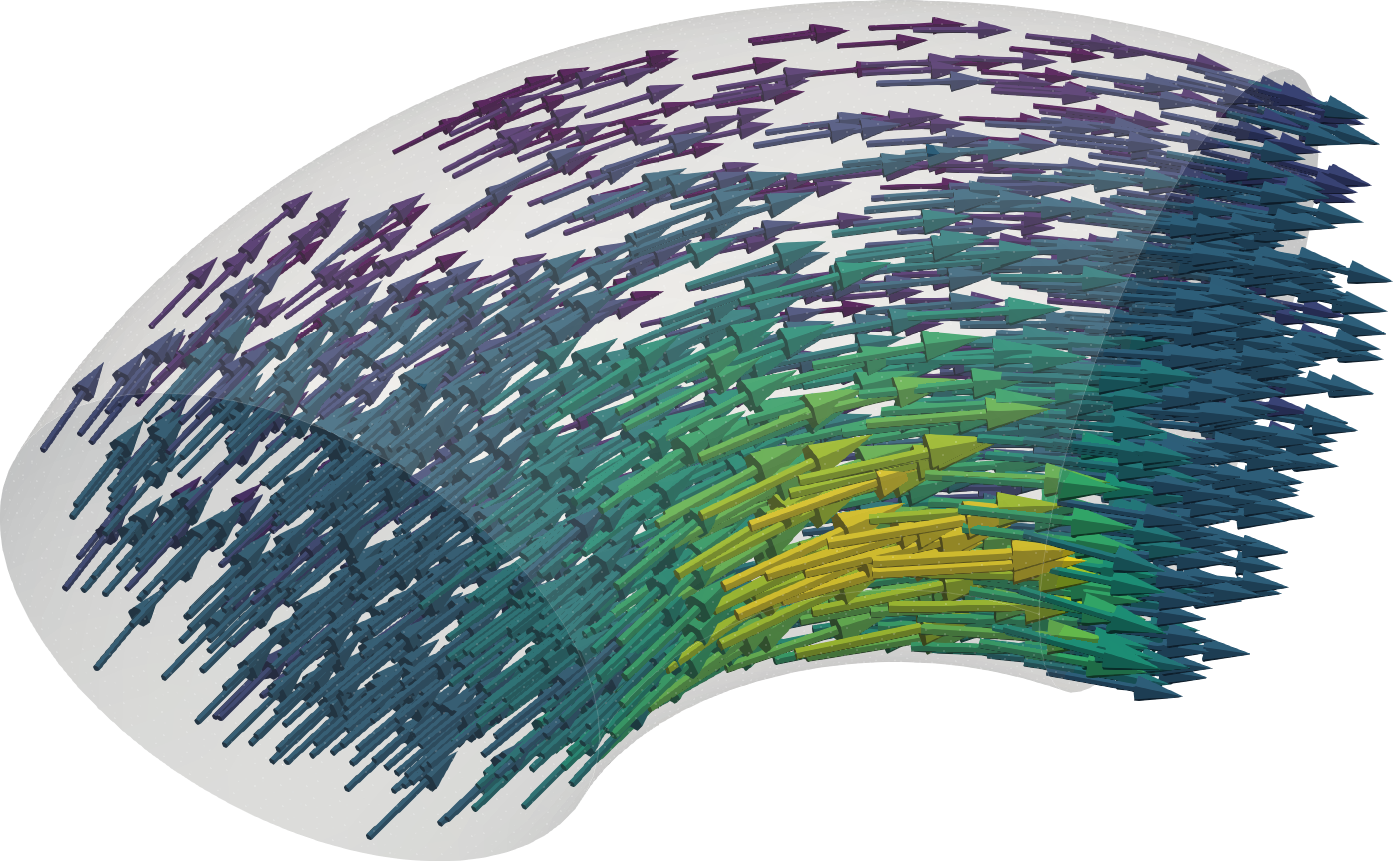}\hspace{3mm}
    \includegraphics[width=0.1\linewidth]{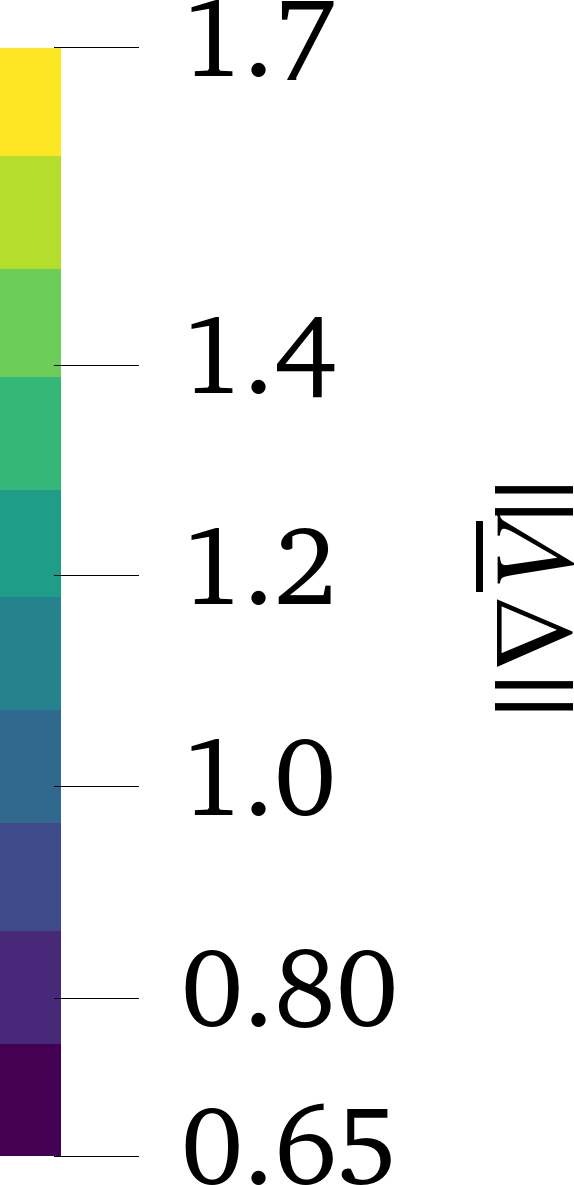}} \\
    %
    \subfloat[Relative error - PINN]{\includegraphics[width=0.28\linewidth]{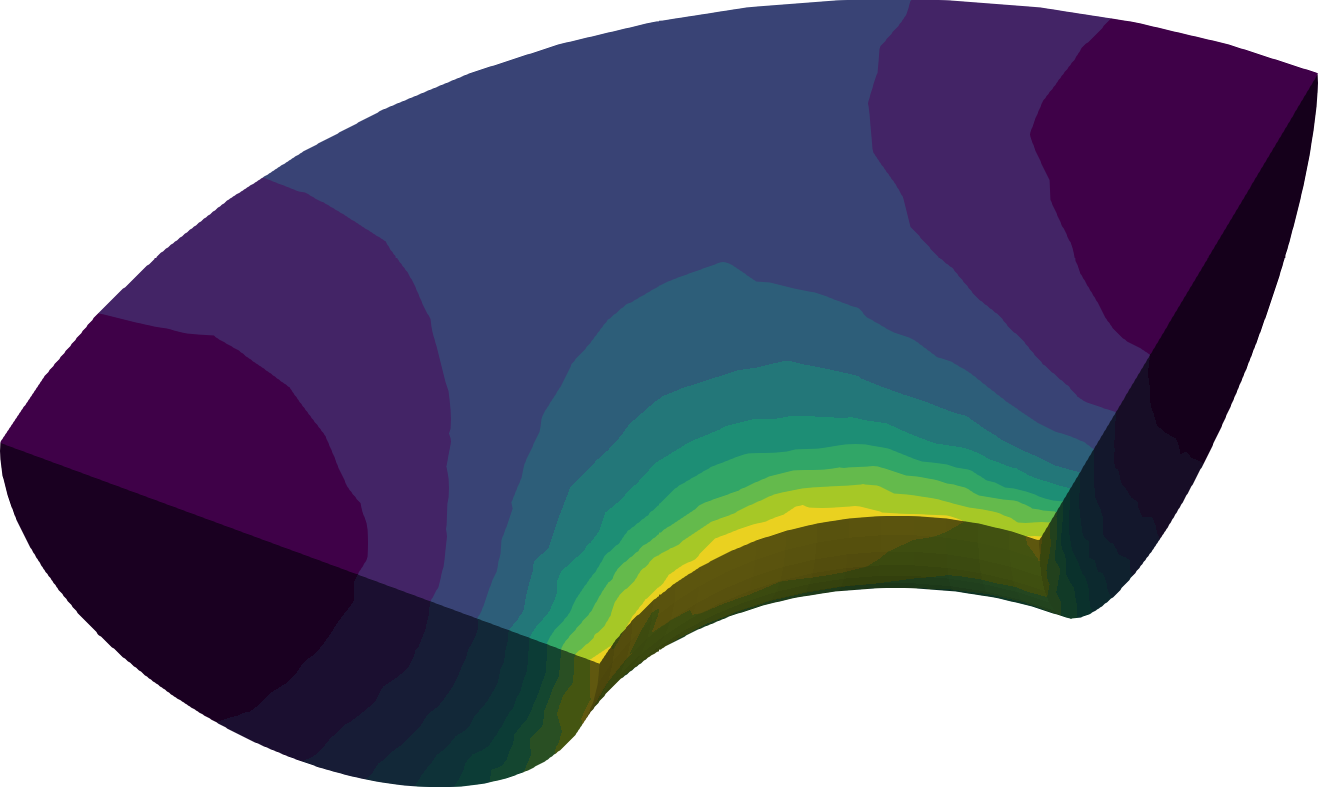}}\hspace{20mm}
    \subfloat[Relative error - HOL]{\includegraphics[width=0.28\linewidth]{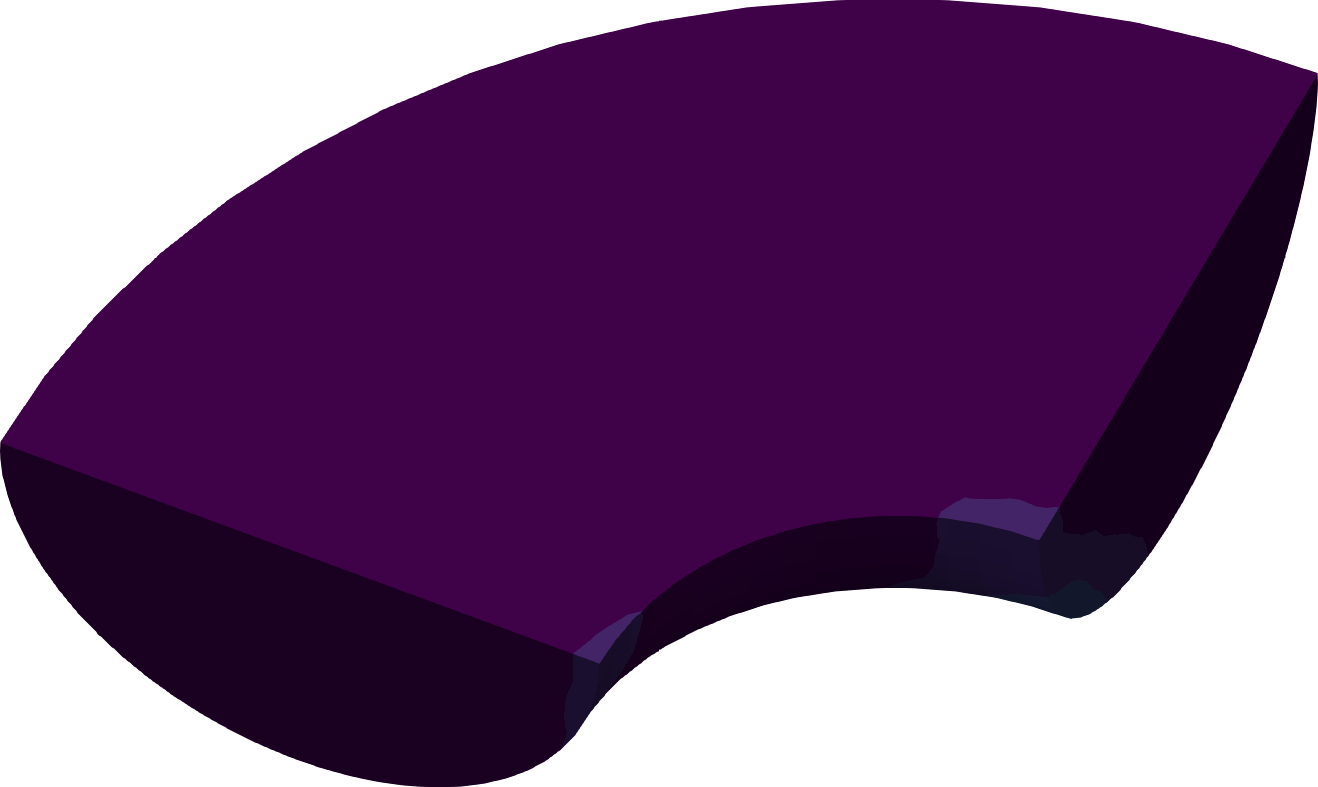}\hspace{3mm}\includegraphics[width=0.1\linewidth]{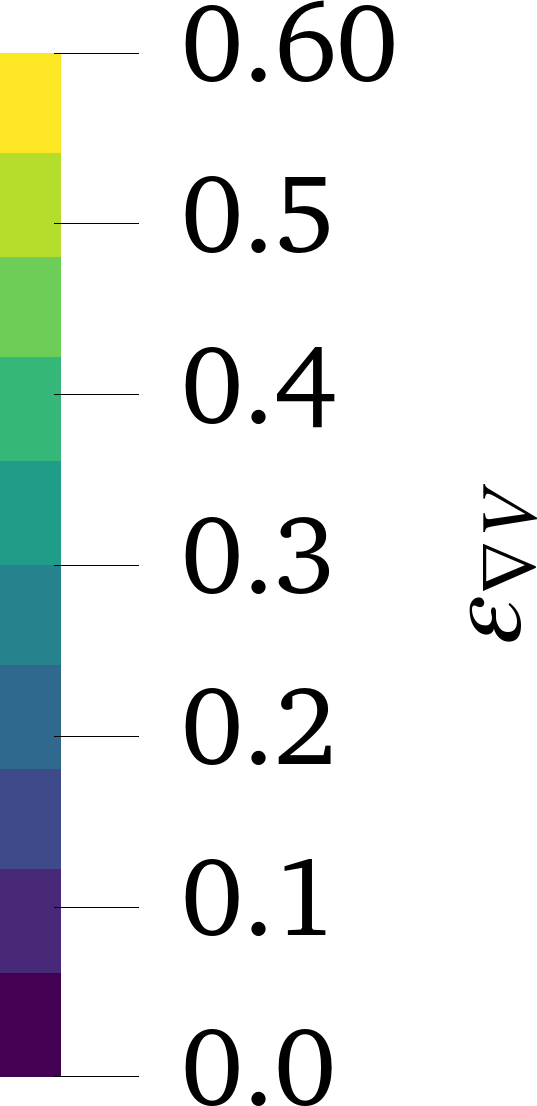}}\\
    %
    \subfloat[Mass conservation - PINN]{\includegraphics[width=0.28\linewidth]{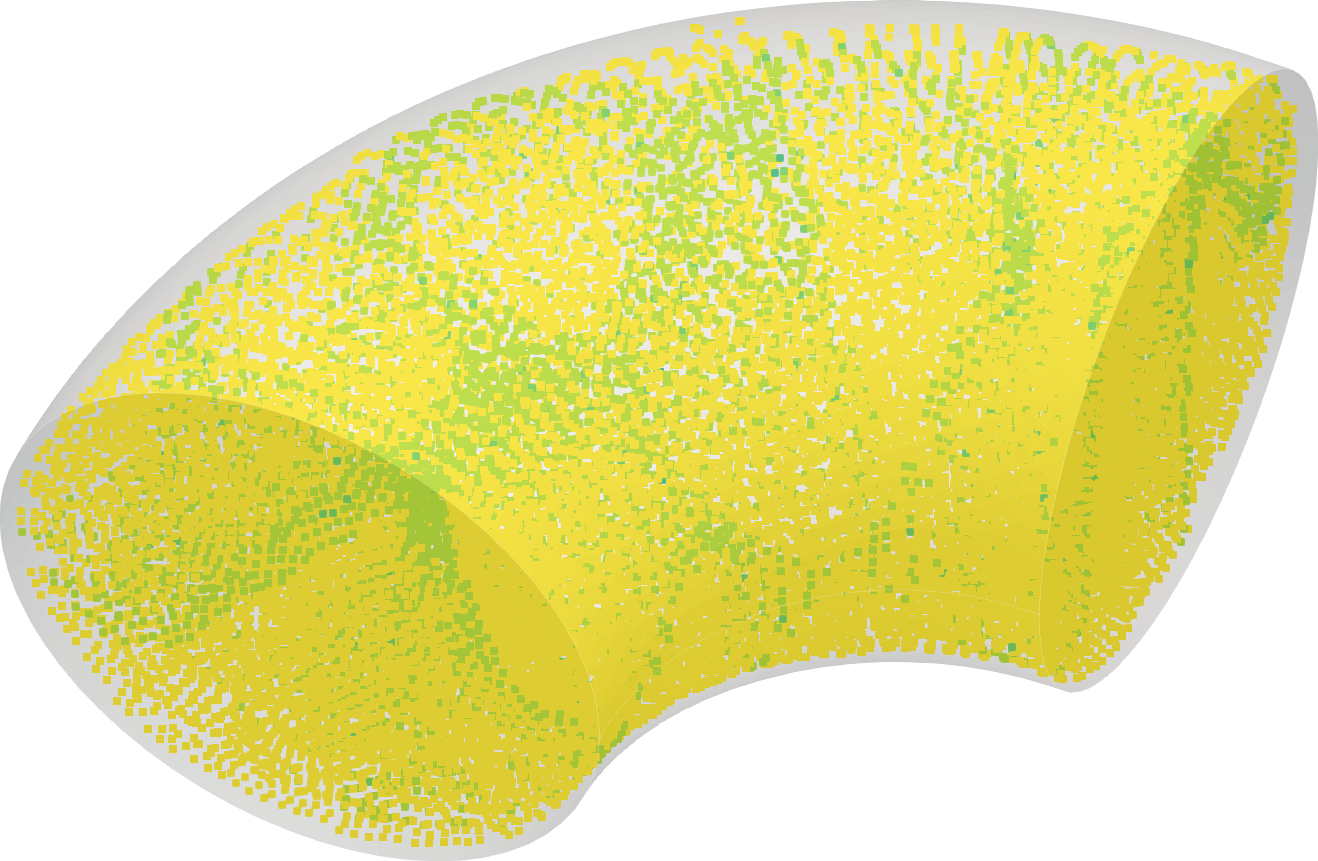}}\hspace{20mm}
    \subfloat[Mass conservation - HOL]{\includegraphics[width=0.28\linewidth]{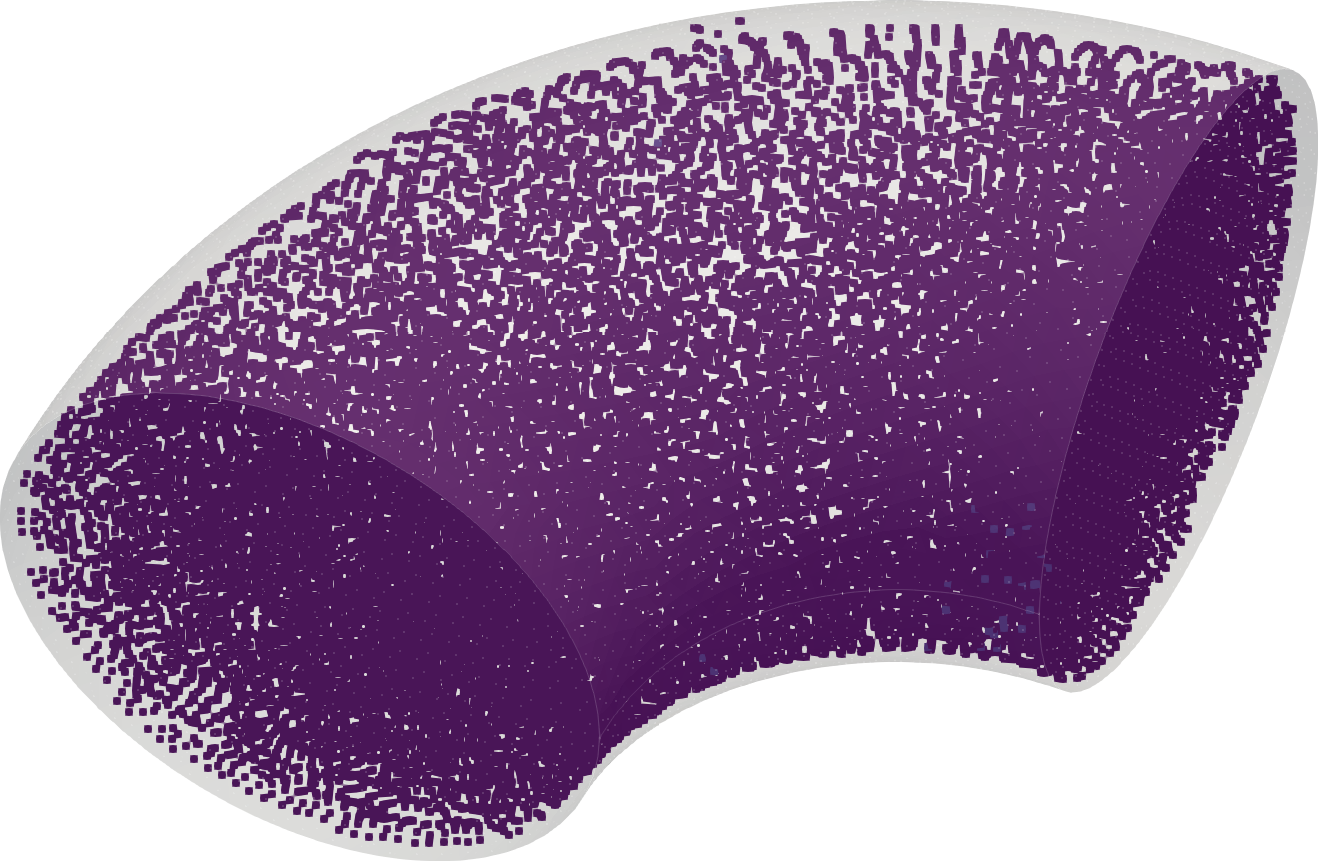}\hspace{3mm}\includegraphics[width=0.1\linewidth]{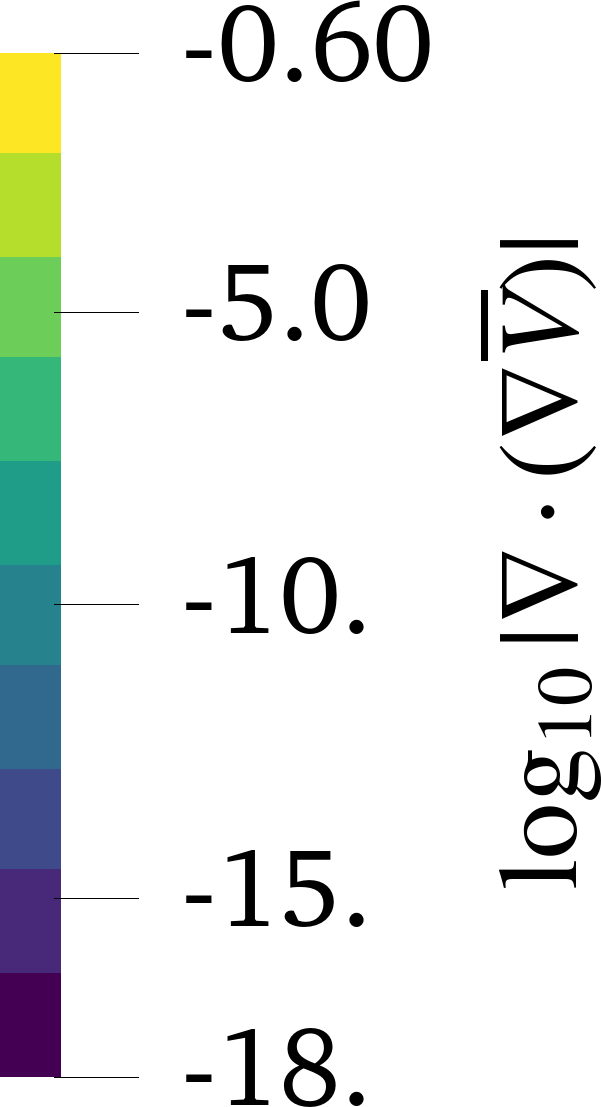}}\\
    \caption{Domain and numerical results. (a) The boundary conditions are imposed only to $\nabla V$. (a)-(c) Visualization of $\nabla \overline{V}$, which can be interpreted as the velocity field of an incompressible and inviscid flow. (e)-(f) Relative error defined in \eqref{eq:errs_grad_lapl_local}.  The proposed method (HOL) shows an error below $5\%$. (g)-(h) The term $\divergence \nabla \overline{V}$ is displayed; it can be interpreted as the divergence of the velocity field, which vanishes by incompressibility. The proposed method satisfies this property up to machine precision.}
    \label{fig:test_lapl_flow}
\end{figure}

The main features of this example are the absence of Dirichlet boundary conditions, implying that $V$ is defined up to a constant, and the lack of an analytical solution for $V$. \\
The domain represents a quarter of a torus with inner radius $0.5$ and outer radius $1.5$, see panel (a) of \Cref{fig:test_lapl_flow}. We define the following boundary operators:
\begin{align*}
\begin{aligned}
	&\boundarylapl[1](V) \coloneq \nabla V \cdot \bfn,  && \on \partial{\Omega_1}, \\
	&\boundarylapl[2](V) \coloneq  \nabla V \cdot \bfn  + 1, && \on \partial{\Omega_2}, \\
	&\boundarylapl[3](V) \coloneq  \nabla V \cdot \bfn - 1, && \on \partial{\Omega_3},
\end{aligned}
\end{align*}
where $\partial\Omega_1$, $\partial\Omega_2$, and $\partial\Omega_3$ are shown in \Cref{fig:test_lapl_flow}, panel (a). Consequently, we set the following boundary conditions
\begin{align*}
	\boundarylapl[1](V) = 0, && \boundarylapl[2](V) = 0, &&& \boundarylapl[3](V) = 0.
\end{align*}
This test can be interpreted as an inviscid, irrotational fluid flow, where $V$ denotes the velocity potential and $\nabla V$ the velocity field \cite{landau1987}. The boundary operators correspond to imposing in the normal direction at the different boundaries no net flow, constant inflow, and constant outflow, respectively. \\
The layer widths of the neural networks, $N_V$, including input and output, are $(2, 16, 16, 16, 1)$, with the splitting $m^\ell =  n^\ell / 2$ for $\ell = 1,\ldots,4$. We employ the exponential function and cPReLU as activation functions.
We set $\zeta = \sin(\theta)y+\cos(\theta)z + i x$ and employ the composite midpoint rule with $32$ intervals uniformly spaced in $[0, 2\pi)$. 
Let us define the mask function $\omega$ as follows:
\begin{equation}\label{eq:selection_fcn}
	\omega_i(\bfx) = 
	\begin{cases*}
		1 \quad \text{if } \bfx \in \partial\Omega_i, \\
		0 \quad \text{otherwise}.	
	\end{cases*}	
\end{equation}
The loss function is given by
\begin{subequations}
\begin{align}
	&\loss = \frac{1}{n_\text{train}}\sum_{i=1}^{n_\text{train}} \loss_p(\bfx_i) + \loss_c(\bfx_i), \\
    &\loss_p(\bfx) = \alpha_1 \, \omega_1(\bfx) \| \boundarylapl[1](\overline{V}(\bfx)) \|^2 + \alpha_2 \, \omega_2(\bfx)\| \boundarylapl[2](\overline{V}(\bfx)) \|^2 + \alpha_3 \, \omega_3(\bfx)\| \boundarylapl[3](\overline{V}(\bfx)) \|^2, \label{eq:pointwise_loss} \\
    &\loss_c = 10^{-4}\, \overline{V}^2, \label{eq:loss_stab}
\end{align}
\end{subequations}
where $n_\text{train} = \numprint{3128}$ while $n_\text{text} = 348$ points form the test dataset; the weights $\alpha_i$, $i = 1,2,3$, are manually tuned and set to $1$, $10$, and $10$, respectively, in order to balance the accuracy across the domain of the potential $\overline{V}$ defined in \eqref{eq:V_discrete}; $\loss_c$ is a small penalty term added to the loss function to promote potentials with a small average over the training dataset. This regularization is beneficial because, under the prescribed boundary conditions, the problem is ill-posed as the solution is not unique. However, the primary quantity of interest is the gradient of $V$, rather than $V$ itself. Therefore, to stabilize the training procedure, we introduce the penalty term mentioned earlier. \\
Training lasts $\numprint{2000}$ epochs, until the loss function reaches a low stationary value. More details can be found in \Cref{sec:appendix_flow}.

The reference solution is obtained with the FEM software Abaqus using a mesh of 10-node, second order tetrahedral elements with mean size 0.05. The solution is represented in panel (b) of \Cref{fig:test_lapl_flow}, which provides a graphical representation of the flow field. The vectors follow the geometric boundaries and increase in magnitude near the inner radius of the torus, as expected.
In panel (d) of the same figure, we report the gradient of the solution, $\nabla \overline{V}$, obtained with the proposed method. Compared with the reference solution, a good visual agreement can be observed.
This qualitative result is further confirmed quantitatively by the relative error $\eps_{\nabla V}$, displayed in panel (f) of \Cref{fig:test_lapl_flow}. The error remains below $5\%$ throughout the domain.

\paragraph{Evaluation against PINN.}
To better assess the properties of the proposed method, hereafter referred to as HOL (from holomorphic), we present a comparison with a standard PINN formulation \cite{Raissi2019}. Specifically, we adopt a PINN approach based on a single neural network to approximate the mapping $(x,y,z) \to V$. Also, we design the network to have both a comparable number of layers and total trainable parameters to $N_V$, which is equal to $\numprint{609}$. Accordingly, the network has layers of width $(3, 16, 16, 16, 1)$, yielding a total number of trainable parameters equal to $\numprint{625}$. 
To ensure smoothness of the output, the Softsign activation function is employed. 
The PINN loss function is defined as
\begin{align*}
	&\loss_\text{PINN} = \sum_{i=1}^{n_\text{train}} \loss_p(\bfx_i) + \loss_c(\bfx_i)  + \sum_{i=1}^{n_{\text{train},\Omega}} \loss_\text{in}(\bfx_i), \\
    & \loss_\text{in}(\bfx) = \|\laplacian V(\bfx)\|^2,
\end{align*}
where $\loss_p$ and $\loss_c$ ar defined in \eqref{eq:pointwise_loss} and \eqref{eq:loss_stab} respectively; $n_{\text{train},\Omega}$ denotes the number of training points in the domain, equal to $\numprint{44514}$. Training lasts $\numprint{2000}$ epochs, until the loss function reaches a fairly stationary low value. The training times for both PINN and HOL methods are reported in \Cref{tab:pinn_hol_lapl}. We observe that, for the same number of epochs, equal to 2000, the ratio between the two training times is approximately 5 in favour of the HOL method.
Further details can be found in \Cref{sec:appendix_flow}.
In panel (c) of \Cref{fig:test_lapl_flow}, we report the gradient of the solution, representing the velocity field, obtained with the PINN method. The corresponding relative error with respect to the reference solution is shown in panel (e).
The peak local error $\eps_{\nabla V}$ associated with the PINN approach is higher than that of the proposed method, reported in panel (f) of the same figure.
%
%
The same trend is observed for the global error $\overline{\eps}_{\nabla V}$, which is larger for the PINN method ($\overline{\eps}_{\nabla V} = 0.1856$ against $\overline{\eps}_{\nabla V} = 0.0166$ for the proposed method), as reported in \Cref{tab:pinn_hol_lapl}. \\
An additional relevant feature is illustrated in panels (g) and (h) of \Cref{fig:test_lapl_flow}, which report the values of $\laplacian V$ in $\Omega$, expressed as $\divergence \nabla V$ to emphasize its connection with the physical incompressibility constraint, namely the vanishing divergence of the velocity field. Owing to the properties of the proposed method, this constraint is satisfied up to machine precision. In contrast, for the PINN approach, this quantity is influenced by the training procedure and does not attain comparably small values.
\begin{table}[h]
\centering
\begin{tabular}{lllll}
     & $|\mathcal{W}|$  & Tot. trainig points & Training \sib{\min} & $\overline{\eps}_{\nabla V}$ \\
     \hline
 PINN &  625 & \numprint{47642} & 43 & 0.1856 \\
 HOL  & 609  & 3128 & 8 & 0.0166 \\ 
\end{tabular}
\caption{Comparison between the physics-informed neural network (PINN) approach and the proposed method (HOL), in terms of number of trainable parameters, number of training points, training time, and global relative error defined in \eqref{eq:errs_grad_lapl_global}.}
\label{tab:pinn_hol_lapl}
\end{table}

\subsection{Linear elasticity problem}\label{sec:test_mech}
In this section, we show two test cases regarding the problem described in \Cref{sec:3d_elasticity}.
To assess the accuracy of the results, we introduce the following relative errors
\begin{subequations}\label{eq:errs_mech}
\begin{align}
    &\varepsilon_{\bfu}(x,y,z) = \frac{\sqrt{|\Omega|} \, \| \overline{\bfu} - \bfu_\text{ref} \|_2}{ \|\bfu_\text{ref}\|_{\Omega}}, \label{eq:errs_mech_local_u} \\
    &\varepsilon_{\bfsigma}(x,y,z) = \frac{\sqrt{|\Omega|} \, \|\overline{\bfsigma} - \bfsigma_\text{ref} \|_2}{\|\bfsigma_\text{ref}\|_{\Omega}}, \label{eq:errs_mech_local_sigma} \\
    &\overline{\varepsilon}_{\bfu} = \frac{\| \overline{\bfu} - \bfu_\text{ref} \|_{\Omega}}{ \|\bfu_\text{ref}\|_{\Omega}}, \label{eq:errs_mech_global_u}\\
    &\overline{\varepsilon}_{\bfsigma} = \frac{\| \overline{\bfsigma} - \bfsigma_\text{ref} \|_{\Omega}}{\|\bfsigma_\text{ref}\|_{\Omega}}, \label{eq:errs_mech_global_sigma}
\end{align}
\end{subequations}
where the reference quantities, denoted by $\cdot_\text{ref}$, are obtained either analytically or via FEM.

\subsubsection{Elementary cases under uniform stress}\label{sec:test_mech_uniform}
\begin{figure}[H]
    \textbf{Uniaxial stretch} \hspace*{3cm} \textbf{Biaxial stretch} \hspace*{3cm} \textbf{Simple shear} \\
    %
    \vspace*{-3mm}
    \adjustbox{valign=t}{\includegraphics[scale=1]{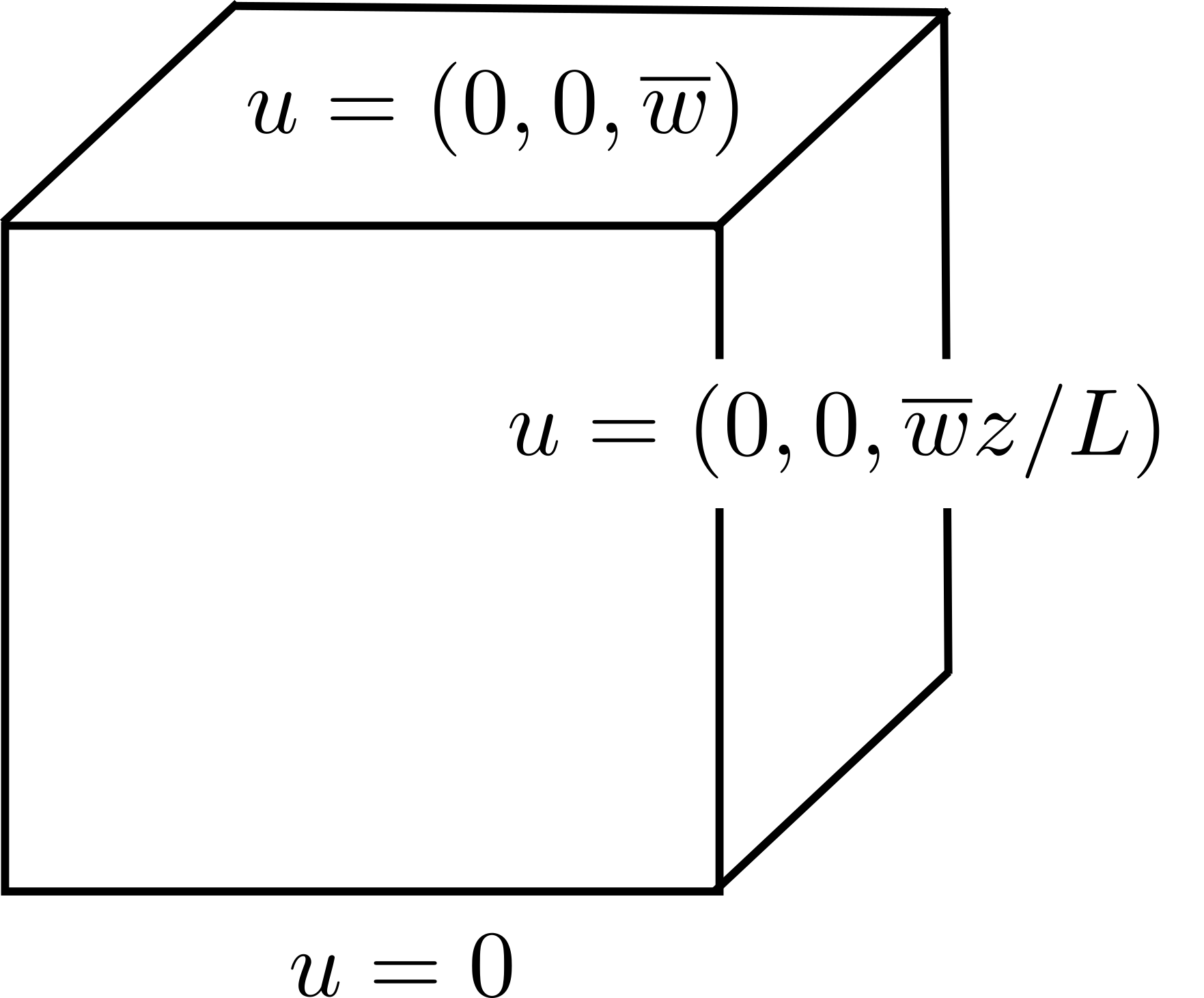}}\hspace{5mm}
    \adjustbox{valign=t}{\includegraphics[scale=1]{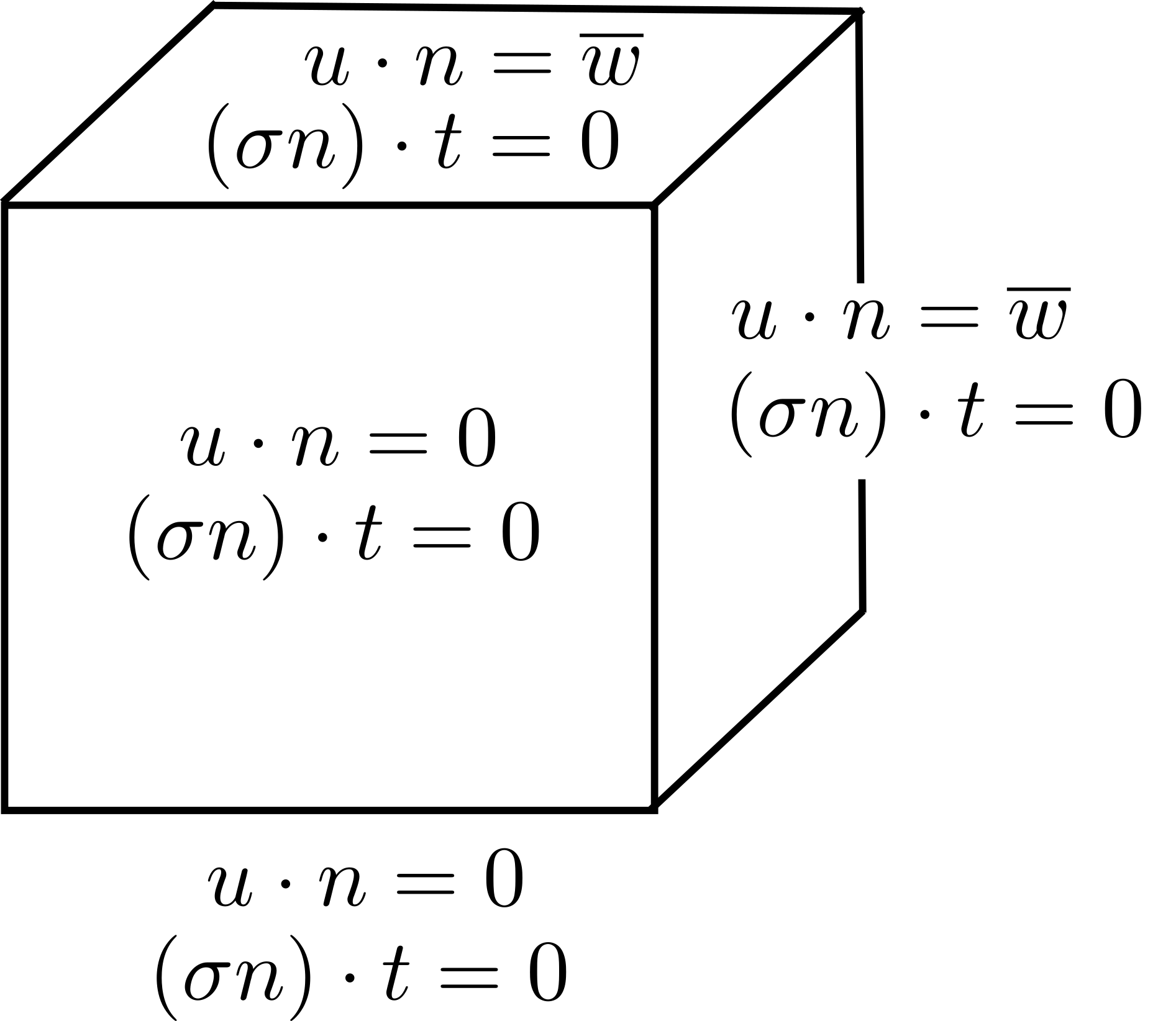}}\hspace{5mm}
    \adjustbox{valign=t}{\includegraphics[scale=1]{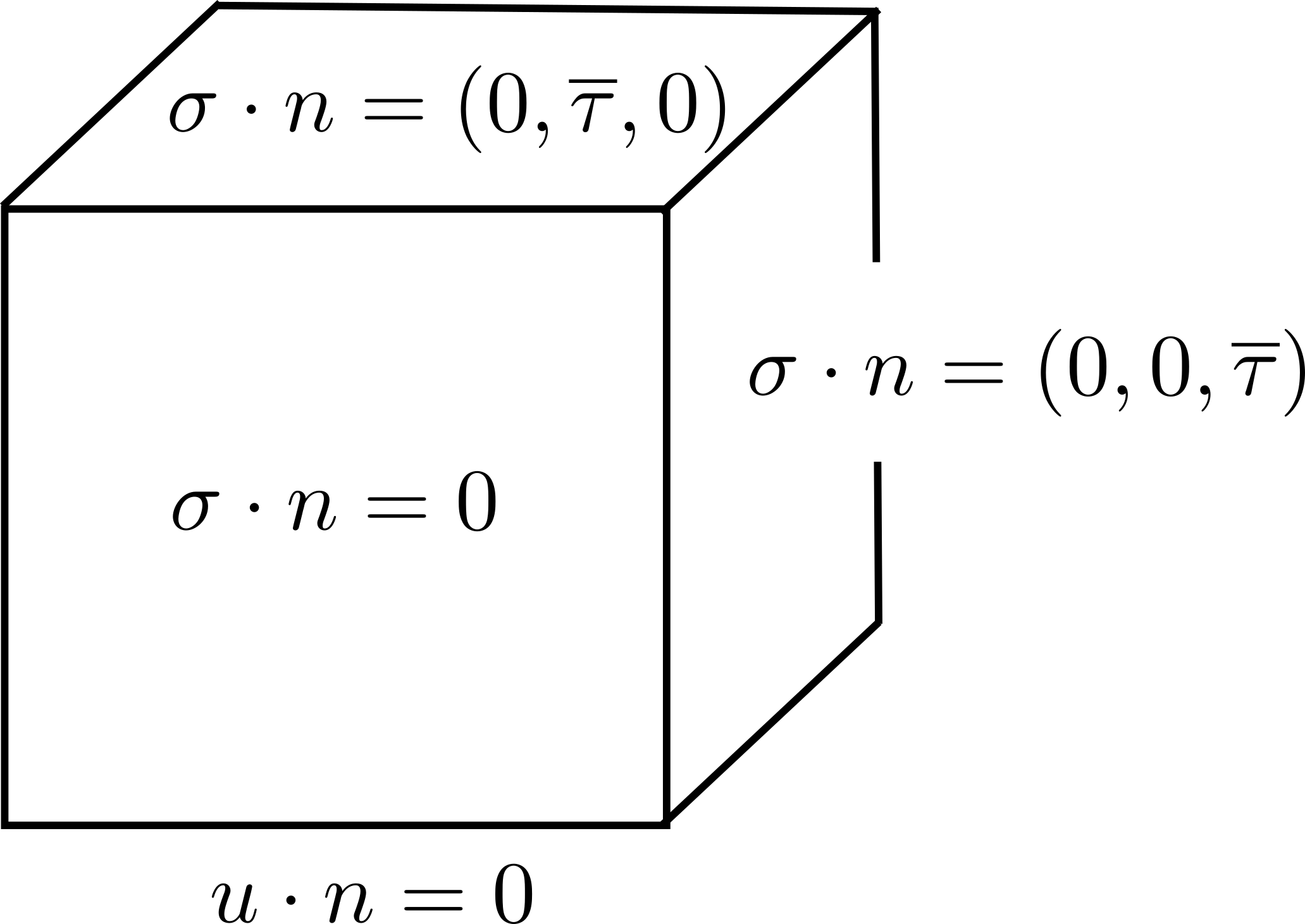}}\\
    \vspace*{3mm}
    \includegraphics[scale=0.8]{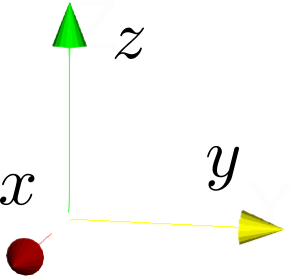} \\
    \vspace*{6mm}
    \rotatebox{90}{\hspace*{-2.5cm}Deformed domain}\hspace{2mm}
    \adjustbox{valign=t}{\includegraphics[scale=0.15]{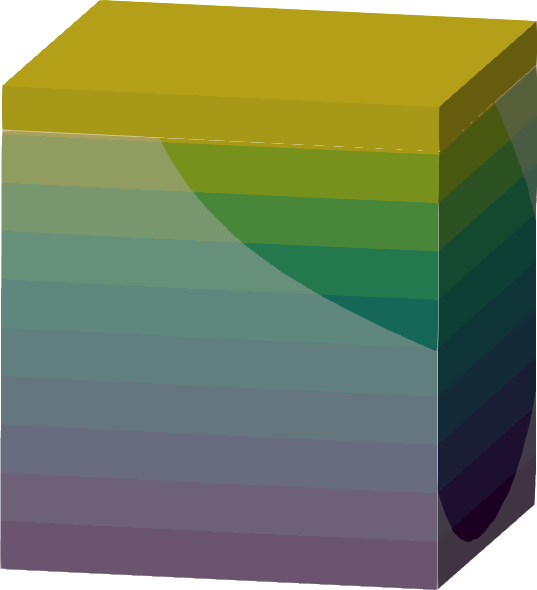}}\hspace{1mm}
    \adjustbox{valign=t}{\includegraphics[scale=0.13]{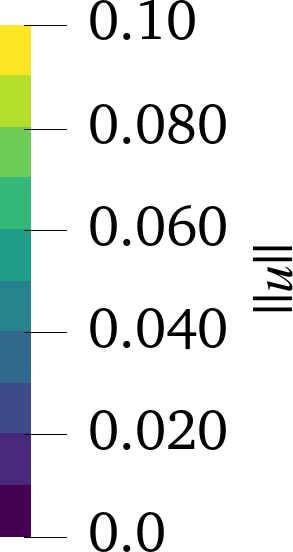}}\hspace{5mm}
    \adjustbox{valign=t}{\includegraphics[scale=0.15]{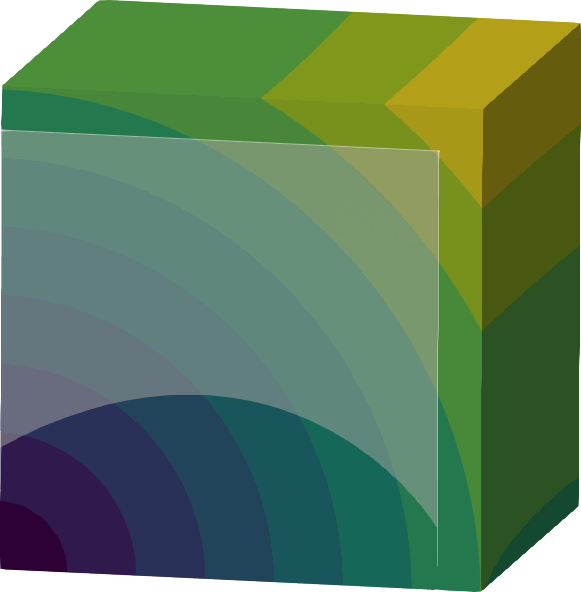}}\hspace{1mm}
    \adjustbox{valign=t}{\includegraphics[scale=0.13]{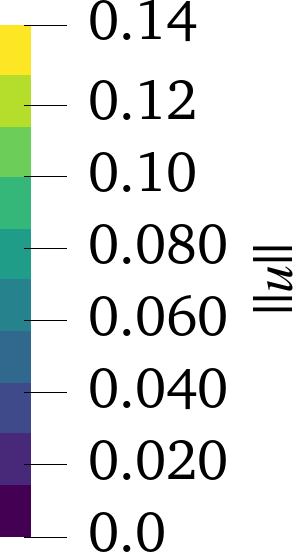}}\hspace{5mm}
    \adjustbox{valign=t}{\includegraphics[scale=0.15]{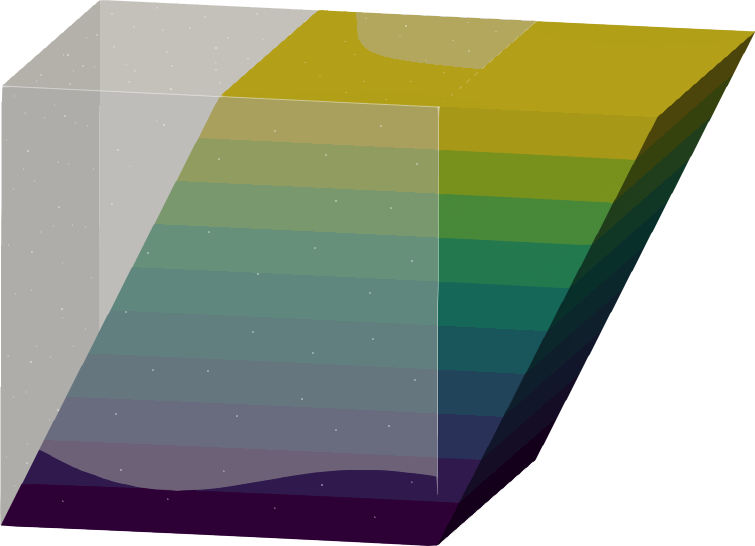}}\hspace{1mm}
    \adjustbox{valign=t}{\includegraphics[scale=0.13]{fig/test_1_deformed_label.png}} \\
    \vspace*{6mm}  
    \rotatebox{90}{\hspace*{-2cm}Stress error}\hspace{2mm}
    \adjustbox{valign=t}{\includegraphics[scale=0.15]{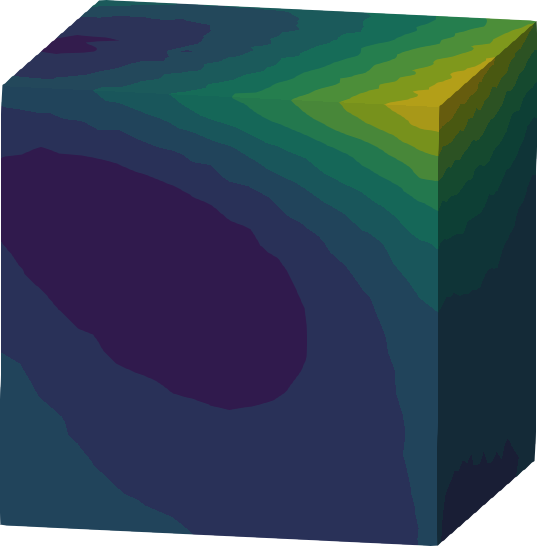}}\hspace{1mm}
    \adjustbox{valign=t}{\includegraphics[scale=0.13]{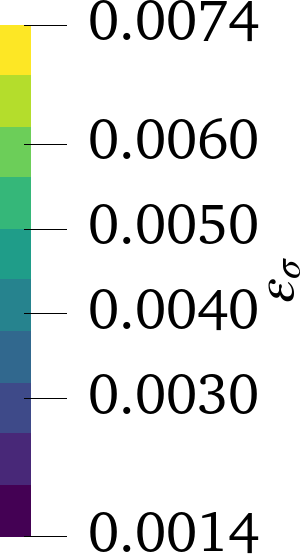}}\hspace{5mm}
    \adjustbox{valign=t}{\includegraphics[scale=0.15]{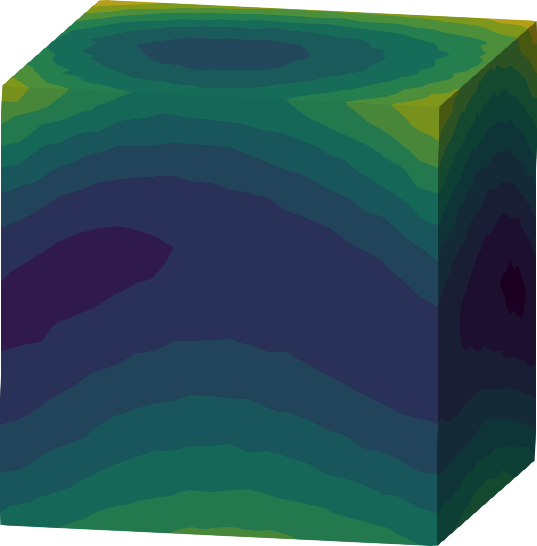}}\hspace{1mm}
    \adjustbox{valign=t}{\includegraphics[scale=0.13]{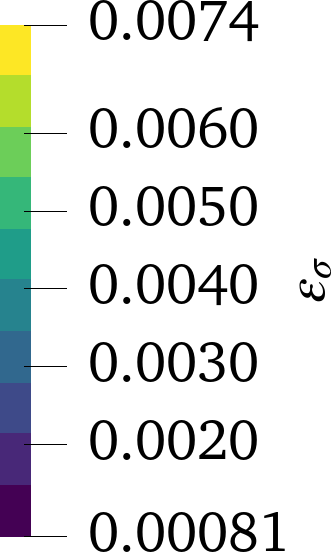}}\hspace{5mm}
    \adjustbox{valign=t}{\includegraphics[scale=0.15]{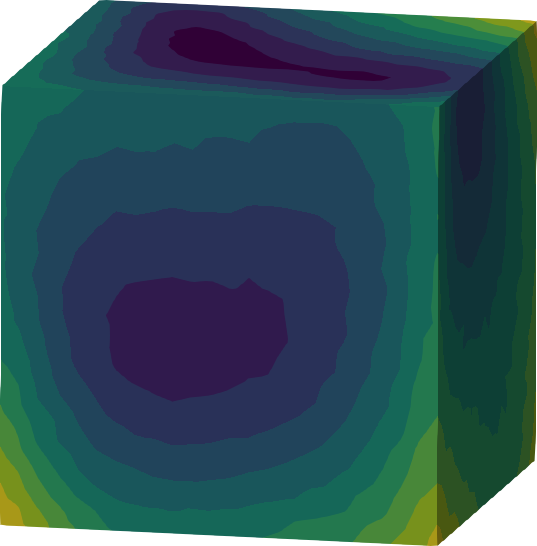}}\hspace{1mm}
    \adjustbox{valign=t}{\includegraphics[scale=0.13]{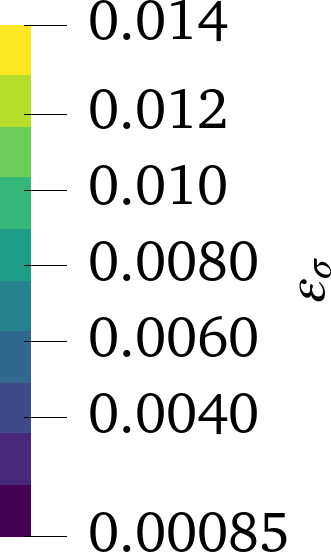}} \\
    \vspace*{6mm}
    \rotatebox{90}{\hspace*{-3cm}Displacement error}\hspace{1mm}
    \adjustbox{valign=t}{\includegraphics[scale=0.15]{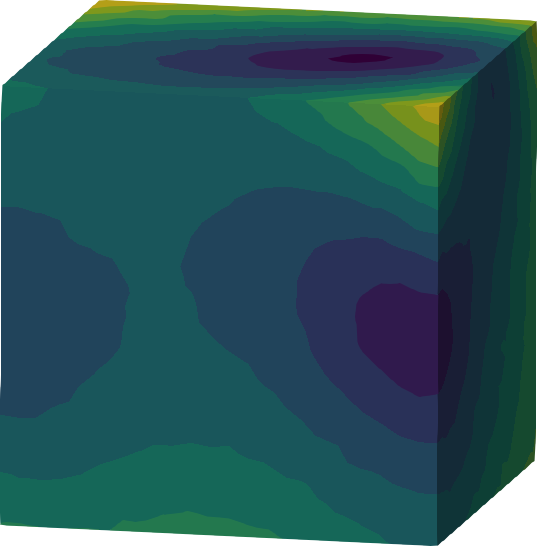}}\hspace{1mm}
    \adjustbox{valign=t}{\includegraphics[scale=0.13]{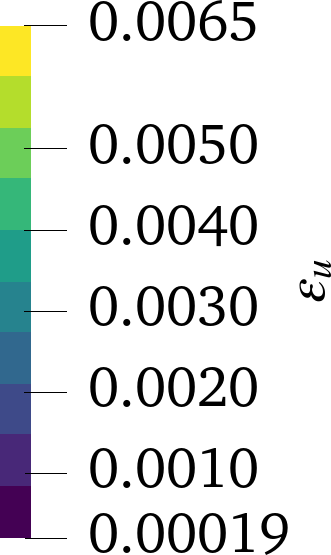}}\hspace{5mm}
    \adjustbox{valign=t}{\includegraphics[scale=0.15]{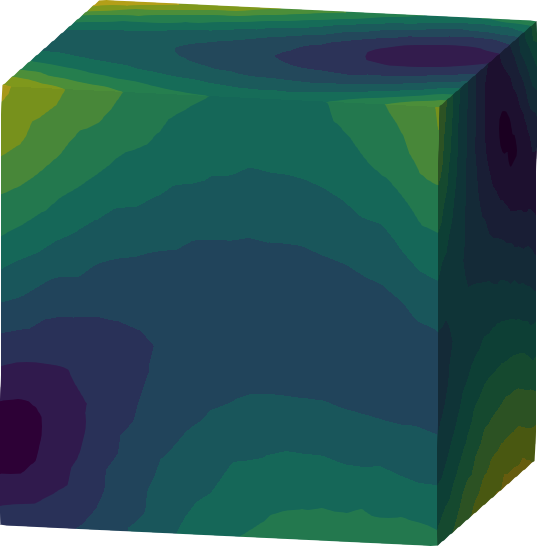}}\hspace{1mm}
    \adjustbox{valign=t}{\includegraphics[scale=0.13]{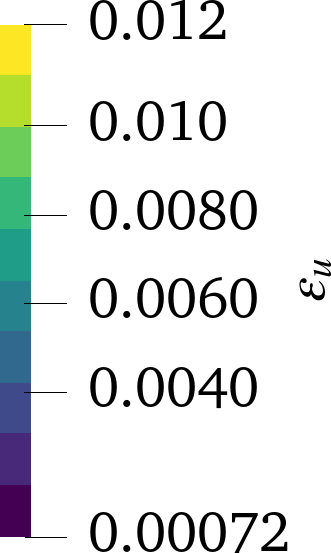}}\hspace{5mm}
    \adjustbox{valign=t}{\includegraphics[scale=0.15]{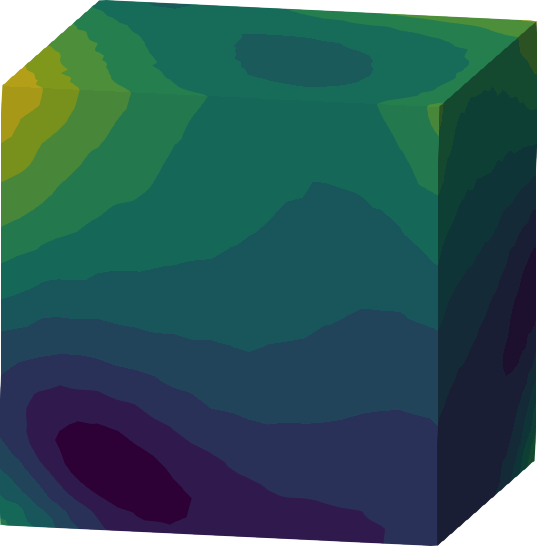}}\hspace{1mm}
    \adjustbox{valign=t}{\includegraphics[scale=0.13]{fig/test_1_error_displ_label.png}} \\
    \caption{First row: domain and boundary conditions for the three cases.
    	Second row: visualization of the deformed geometries obtained through the proposed method, together with the norm of the displacement. Third and fourth rows: displacement and stress errors. The maximum error across all the cases is $1.4\%$.}
    \label{fig:test_mech_uniform}
\end{figure}

This test is chosen because it illustrates the fundamental elastic responses, namely uniform uniaxial stretch, uniform biaxial stretch, and simple shear. Moreover, a wide range of boundary conditions is adopted to provide reliable tests of the proposed method. We consider a material with $\lambda = 1$ and $\mu = 1$. 
In each case, the domain is a cube of side length $L=1$, and different boundary conditions are applied to obtain the three cases mentioned above; see \Cref{fig:test_mech_uniform}. The same figure also shows the corresponding deformed configurations, which facilitate the physical interpretation of the elastic responses. 
In the following, we detail the boundary conditions and the corresponding analytical solution.
\begin{itemize}
	\item \textbf{Uniaxial stretch.} This case is displacement-driven: the three components of the displacement are prescribed on the cube faces as $\bfu = (0,0,\overline{w}z/L)$, where $\overline{w}$ is the vertical displacement prescribed at the top, given by $\overline{w} = 0.1$. The analytical expressions of the non-zero stress and displacement components are
	\begin{align*} 
		&\sigma_{xx} = \lambda \frac{\overline{w}}{L}, & u_z(z) = \frac{\overline{w}}{L} z, \\
		&\sigma_{yy} = \lambda \frac{\overline{w}}{L}, \\
		&\sigma_{zz} = \lambda \frac{\overline{w}}{L}.
	\end{align*}
	\item \textbf{Biaxial stretch.} In this case, zero shear traction boundary conditions are applied on all faces. Also, vanishing normal displacement is prescribed on all faces, except for the top and right faces where the normal displacement is set to $\overline{w} = 0.1$. The exact non-zero stress and displacement components are
	\begin{align*} 
		&\sigma_{xx} = 2\lambda\frac{\overline{w}}{L}, & u_y(z) = \frac{\overline{w}}{L} z, \\
		&\sigma_{yy} = 2(\lambda + \mu) \frac{\overline{w}}{L}, & u_z(z) = \frac{\overline{w}}{L} z, \\
		&\sigma_{zz} = 2(\lambda + \mu) \frac{\overline{w}}{L}.
	\end{align*}
	\item \textbf{Simple shear.} This case is stress-driven, with the bottom face fixed and traction applied on all the other faces. The applied shear stress is $\overline{\tau} = 0.5$. The exact non-zero stress and displacement components are
	\begin{align*}
		&\sigma_{yz} = \overline{\tau}, &\quad u_z = \overline{\tau}z/\mu.
	\end{align*}
\end{itemize}
%
%
In each case, we employ four semi-holomorphic neural networks, $N_{\phi_x}$, $N_{\phi_y}$, $N_{\phi_z}$, and $N_{\chi}$, all sharing the same architecture with layer sizes $(2, 32, 1)$. The value $m^\ell$, which defines the split between the holomorphic and non-holomorphic parts, is set to $m^\ell = n^\ell / 2$. An exponential activation function is used in each layer for the holomorphic component, while the cPReLU function is adopted for the non-holomorphic component. \\
We set, for the uniaxial and biaxial cases, $\zeta = \sin(\theta)x+\cos(\theta)y + i z$, and, for the shear case, $\zeta = \sin(\theta)y + \cos(\theta)z + i x$. In each case, due to the simplicity of the problems at hand, we adopt the composite midpoint rule \cite{Burden2011} with only $4$ intervals equally spaced in $[0, 2\pi)$. \\
In each case, the loss function used depends on the squared Euclidean norm of the residual of the boundary conditions:
\begin{align*}
	&\loss = \sum_{i=1}^{n_\text{train}} \loss_p, \\
    &\loss_p(\bfx) = \sum_{j=1}^{6} \alpha_j\omega_j(\bfx)\|\boundarymech[j](\overline{\bfu}(\bfx), \overline{\bfsigma}(\bfx))\|^2, 
\end{align*}
where $\boundarymech[j]$ is the boundary condition operator on the face $j$ of the six faces of the cube, as described above and illustrated in \Cref{fig:test_mech_uniform}; the weights are $\alpha_j = 1$, $\forall j$; $\overline{\bfu}$ and $\overline{\bfsigma}$ are, respectively, the displacement and the stress obtained through the neural networks. We remind that $\overline{\bfsigma}$ is computed through \eqref{eq:stress_tensor} where the derivatives of the displacement are computed using the automatic differentiation. Here, $n_\text{train} = 238$ while $n_\text{test} = 26$ points form the test dataset. Training lasts $\numprint{2000}$ epochs, until the loss function reaches a fairly stationary low value. Further details can be found in \Cref{sec:appendix_uniform}.

The relative errors defined in \eqref{eq:errs_mech_local_u}-\eqref{eq:errs_mech_local_sigma} are represented in \Cref{fig:test_mech_uniform}. We observe that the method accurately captures the elementary displacement states as well as the corresponding stress states. In each case, the error is indeed small, with a local maximum value across all cases of approximately $1.4\%$.

\subsubsection{3D device subjected to non-uniform stress}\label{sec:test_mech_non_uniform}

\begin{figure}[H]
    \subfloat[Domain and b.c.]{\includegraphics[width=0.26\linewidth]{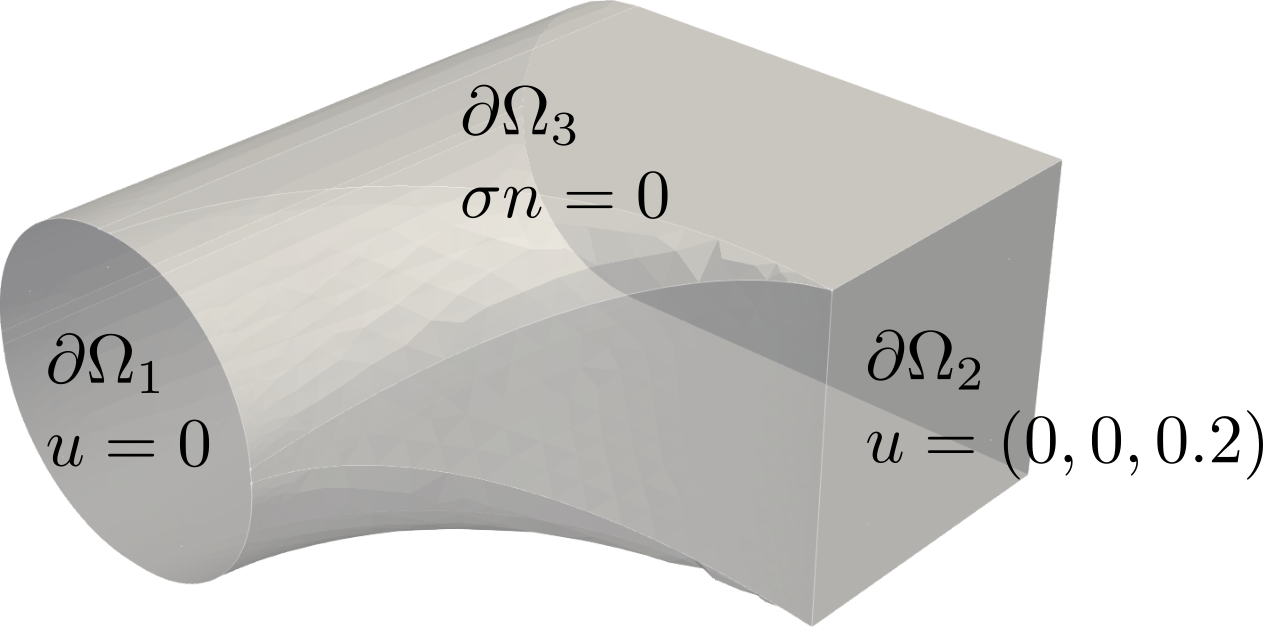}}\hspace{5mm}
    \subfloat[Traction - FEM]{\includegraphics[width=0.25\linewidth]{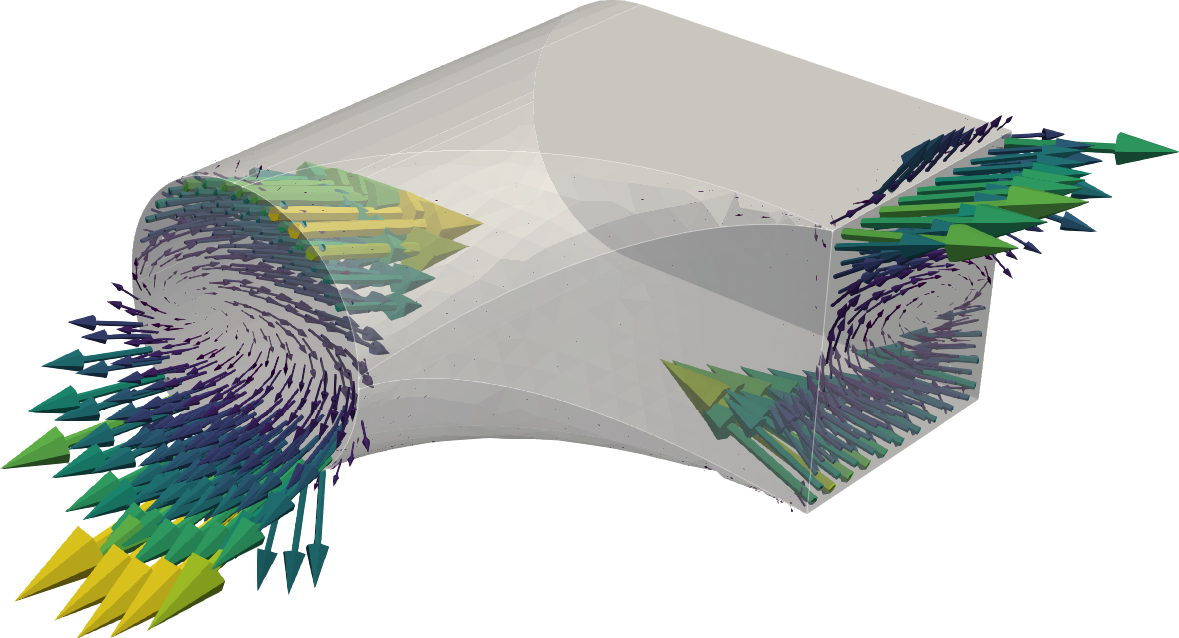}}\hspace{5mm}
    \subfloat[Traction - HOL]{\includegraphics[width=0.25\linewidth]{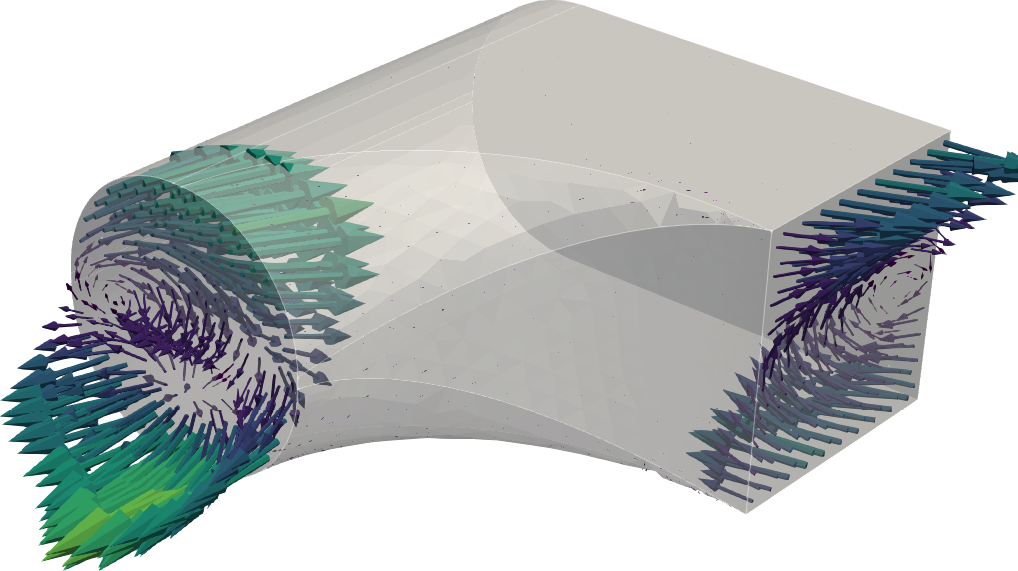}}\hspace{1mm}
    \includegraphics[width=0.08\linewidth]{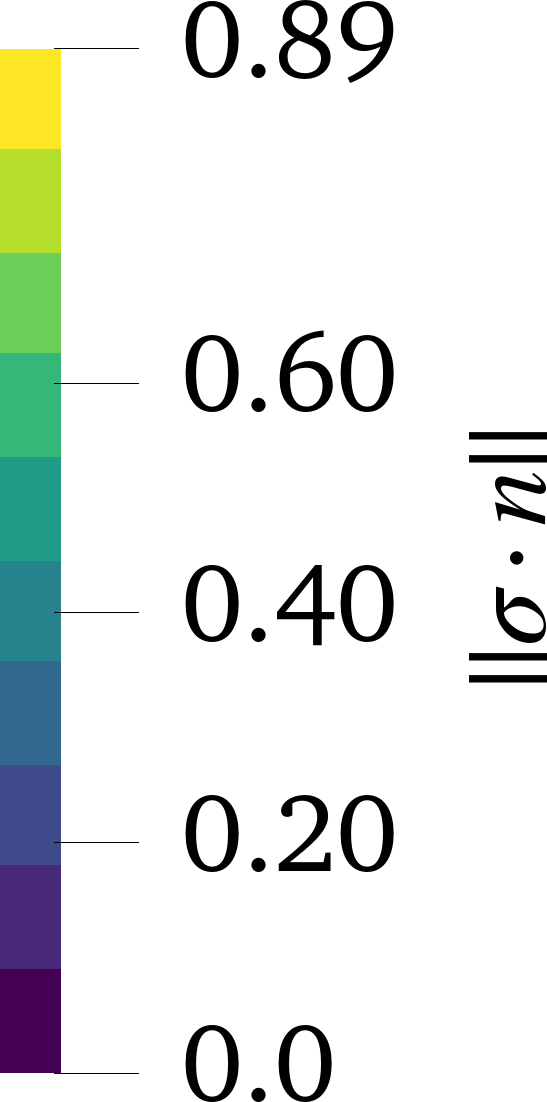} \\
    \vspace*{3mm}
    \includegraphics[scale=0.8]{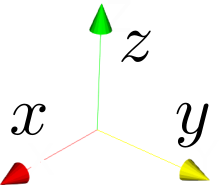} \\
    \vspace*{6mm}
    \subfloat[Def. domain - FEM]{\includegraphics[width=0.2\linewidth]{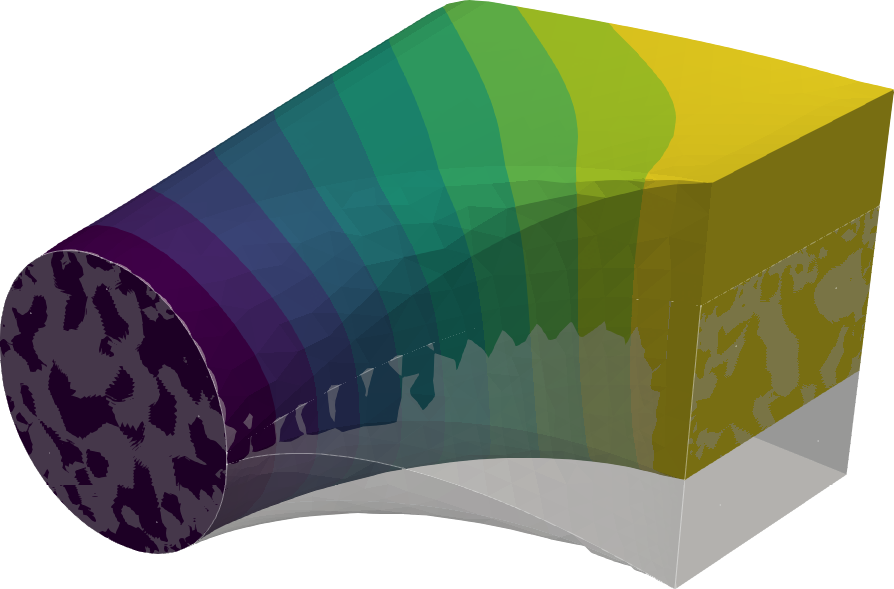}}\hspace{5mm}
    \subfloat[Def. domain - HOL]{\includegraphics[width=0.2\linewidth]{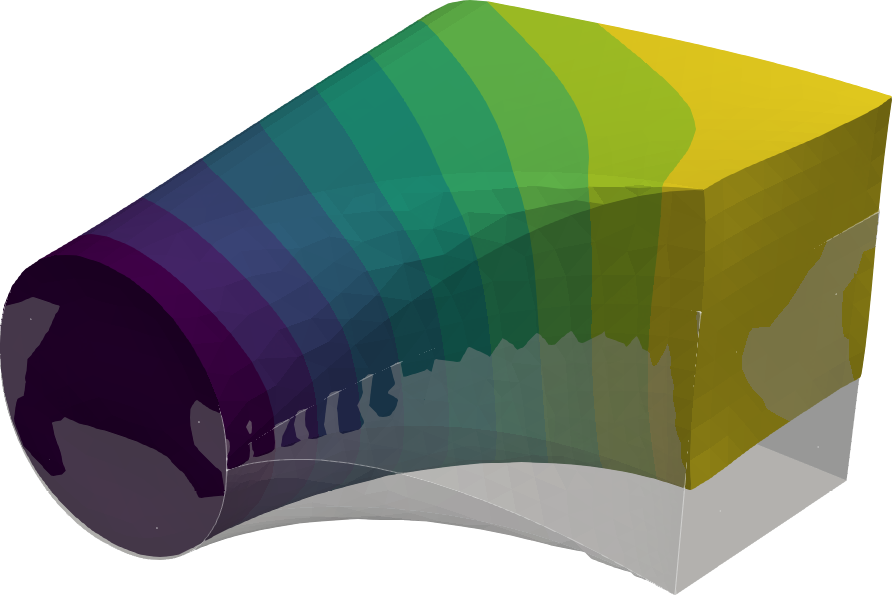}}\hspace{1mm}
    \includegraphics[width=0.08\linewidth]{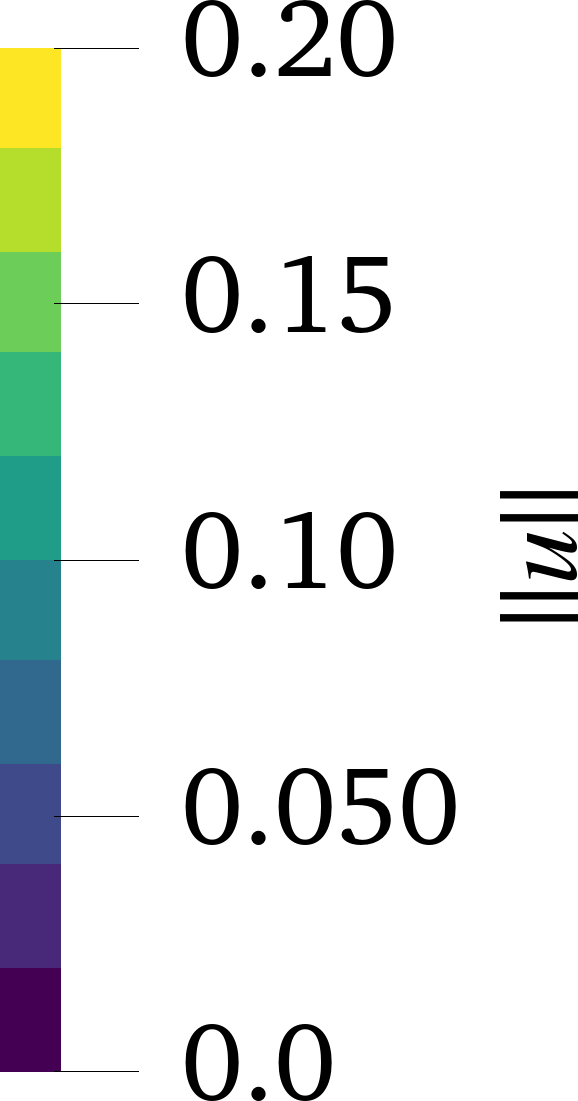} \hspace{5mm}
    \subfloat[Relative err. - Displacement]{\includegraphics[width=0.25\linewidth]{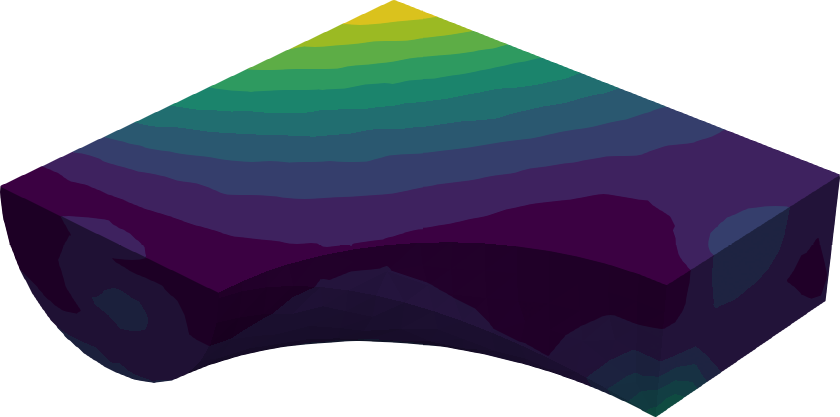}}\hspace{1mm}
    \includegraphics[width=0.08\linewidth]{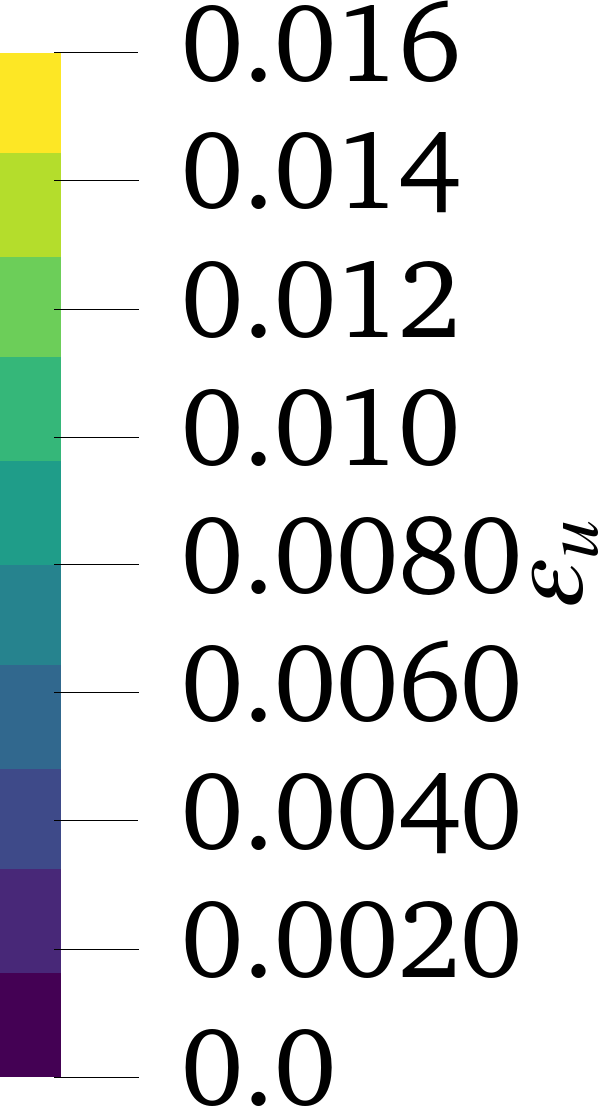} \\
    \vspace*{6mm}
    \subfloat[Stress - FEM]{\includegraphics[width=0.2\linewidth]{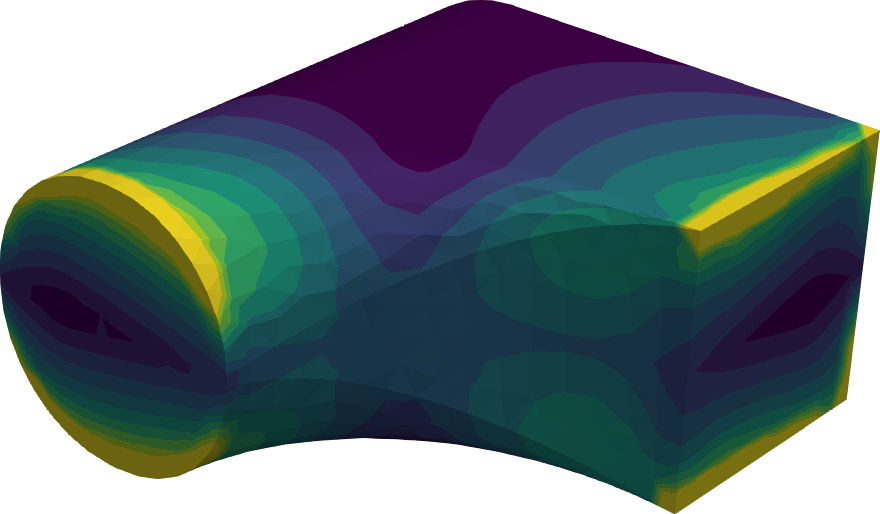}}\hspace{5mm}
    \subfloat[Stress - HOL]{\includegraphics[width=0.2\linewidth]{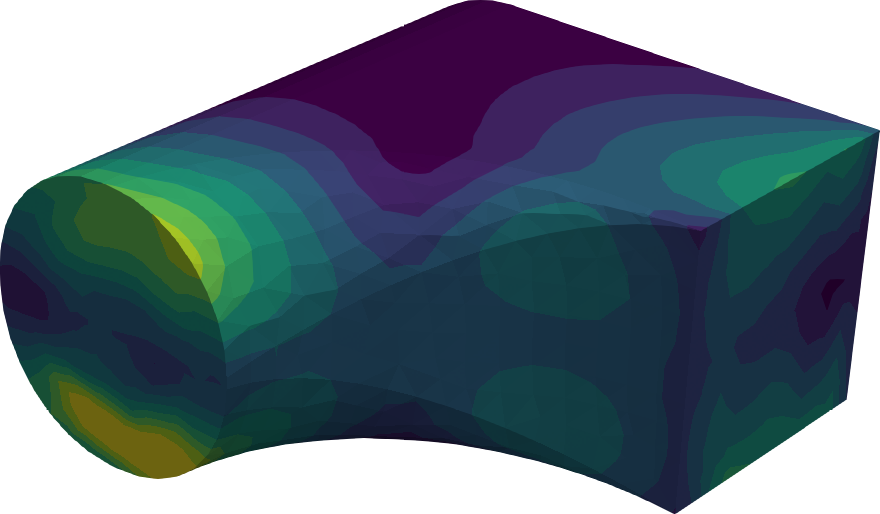}}\hspace{1mm}
    \includegraphics[width=0.08\linewidth]{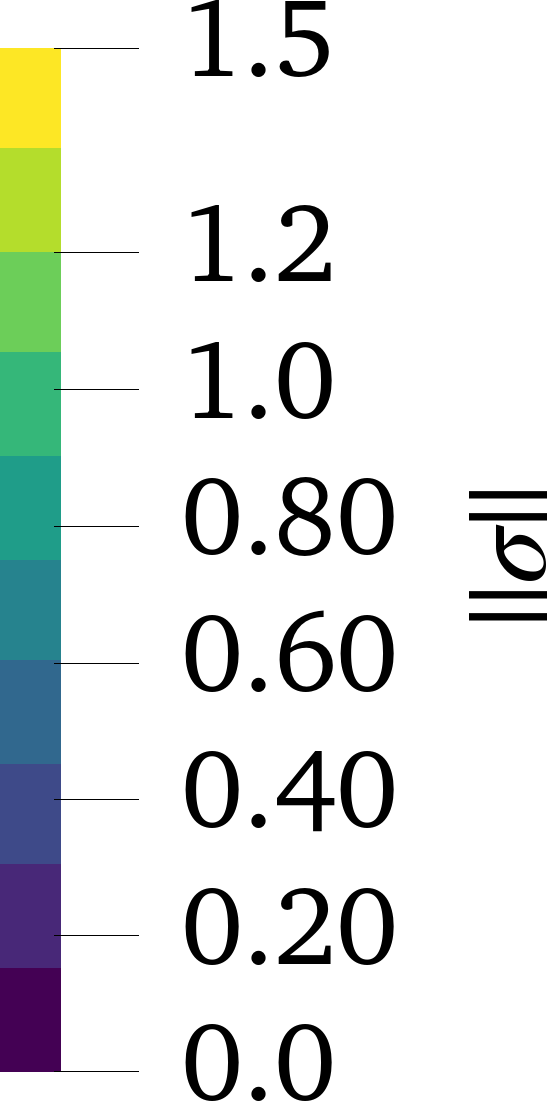} \hspace{5mm}
    \subfloat[Relative err. - Stress]{\includegraphics[width=0.25\linewidth]{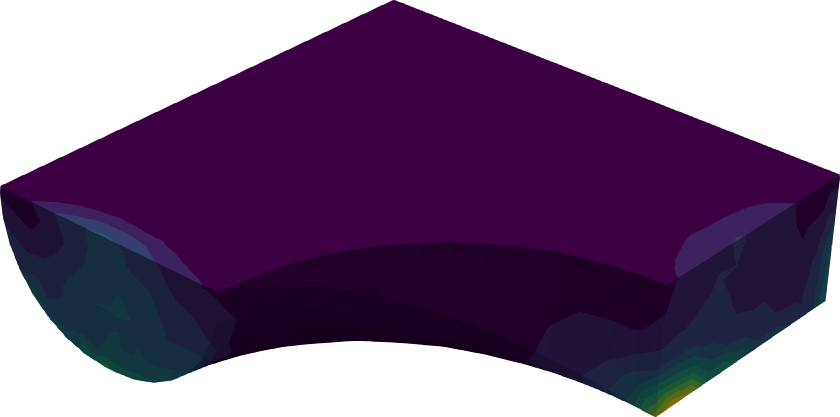}}\hspace{1mm}
    \includegraphics[width=0.08\linewidth]{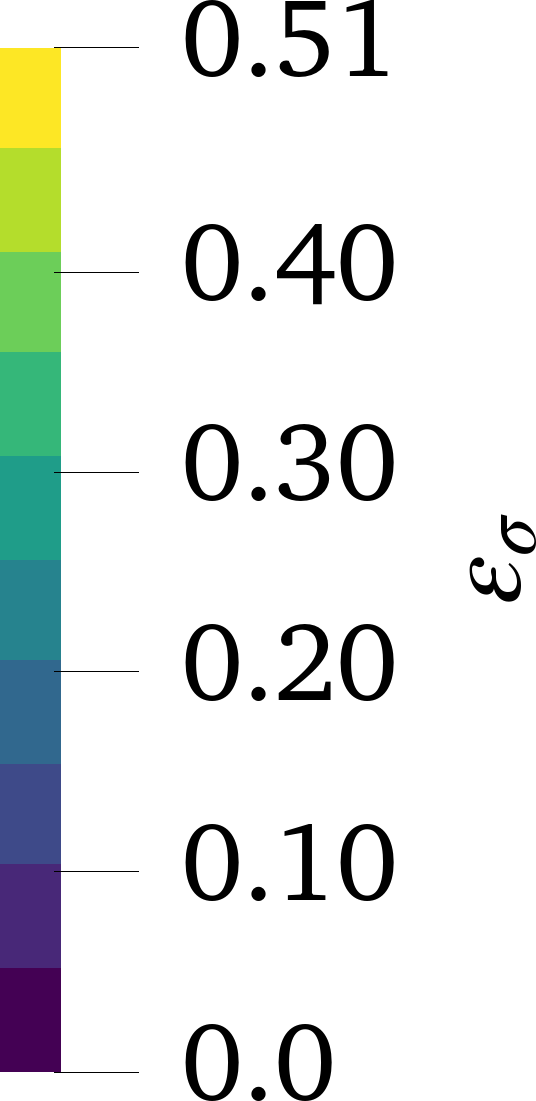} \\
    \caption{First row: domain with boundary conditions and tractions. The tractions obtained through the proposed method, denoted by HOL in the label, are similar to those from FEM. Note that interpolation of stresses introduces spurious tractions near the edges in the FEM case.
    Second row: deformed geometries and displacement error. The error is low and peaks on the boundary.
    Third row: stress and error. The stresses are well reproduced throughout the domain, except in small regions corresponding to stress concentrations.}
    \label{fig:test_mech_non_uniform}
\end{figure}

In this case, we adopt a generic three-dimensional geometry that cannot be obtained by extrusion of a two-dimensional domain; see \Cref{fig:test_mech_non_uniform}. The lengths along the $x$- and $y$-axes are equal to $1$, while the length along the $z$-axis is equal to $0.5$.
Moreover, this test is chosen because both the stress and displacement fields exhibit nonlinear variations along all three spatial directions. The elastic material is characterized by $\nu = 1.1$ and $\lambda = 3.5$. \\
Zero displacement is prescribed on the circular face, while a vertical displacement of magnitude $0.2$ is applied on the opposite square face, see \Cref{fig:test_mech_non_uniform}. Zero traction is imposed on the remaining faces. These boundary conditions are defined by the following boundary operators: 
\begin{align*}
	\begin{aligned}
		&\boundarymech[1](\bfu, \bfsigma) \coloneq \bfu,  && \on \partial{\Omega_1}, \\
		&\boundarymech[2](\bfu, \bfsigma) \coloneq  \bfu - (0, 0, 0.2), && \on \partial{\Omega_2}, \\
		&\boundarymech[3](\bfu, \bfsigma) \coloneq  \bfsigma \bfn,  && \on \partial{\Omega_3},
	\end{aligned}
\end{align*}
where the face indexing is depicted in \Cref{fig:test_mech_non_uniform}. \\
The neural networks share the same architecture, with layer sizes $(2, 64, 64, 64, 64, 1)$. The split value is set to $m^\ell = n^\ell / 2$. An exponential activation function is used in each layer for the holomorphic component, while the cPReLU function is adopted for the non-holomorphic component. Here, we set $\zeta = \sin(\theta)y + \cos(\theta)z + i x$. The integral in \eqref{eq:displ_phi_chi} is discretized using the composite two-point Gauss--Legendre quadrature rule \cite{Burden2011}, with 16 integration intervals equally spaced in $[0, 2\pi)$. Preliminary investigations show no improvement in accuracy when increasing the number of integration intervals. 
The loss function is given by
\begin{align}
	&\loss = \sum_{i=1}^{n_\text{train}} \loss_p(\bfx_i), \nonumber \\
	&\loss_p(\bfx) = \alpha_1 \omega_1(\bfx) \| \boundarymech[1](\overline{\bfu}(\bfx), \overline{\bfsigma}(\bfx)) \|^2 + \alpha_2 \omega_2(\bfx) \| \boundarymech[2](\overline{\bfu}(\bfx), \overline{\bfsigma}(\bfx)) \|^2 + \alpha_3 \omega_3(\bfx) \|\boundarymech[3](\overline{\bfu}(\bfx), \overline{\bfsigma}(\bfx))\|^2, \label{eq:pointwise_loss_hol}
\end{align}
where $\omega_i$ are defined in \eqref{eq:selection_fcn}. We set $\alpha_i = 1,  \ i=1,2,3$. The stress and displacement fields $\overline{\bfsigma}$ and $\overline{\bfu}$ are those obtained from the neural networks. The number of training points is $n_\text{train} = \numprint{7745}$ while $n_\text{test} = 861$ form the test dataset. Training lasts $\numprint{40000}$ epochs, until the loss function reaches a  stationary low value. More details can be found in \Cref{sec:appendix_non_uniform}. \\
Since an analytical solution is not available for this case, we consider as reference a numerical solution obtained with the FEM software Abaqus using 10-node, second order tetrahedral elements with mean size 0.04.

In \Cref{fig:test_mech_non_uniform}, we report the stress and deformation fields for this test. In particular, the second row shows a comparison of the deformed configurations, together with the corresponding relative error. The latter does not exceed approximately $1.6\%$, indicating good precision in the representation of the displacement field. \\
The third row of \Cref{fig:test_mech_non_uniform} presents the stress field and the associated relative error. The stress is well reproduced throughout most of the domain, except in small regions near the edges of the constrained faces, which correspond to stress concentration zones. \\
Furthermore, the first row of the same figure compares the boundary tractions obtained from the FEM solution and the proposed method. The tractions are computed as $\bfsigma \cdot \bfn$. We observe that the tractions are qualitatively well captured. 
It is worth noting that the FEM solution exhibits non-physical tractions near the edges; this is due to the interpolation of the stress tensor from integration points to face centers, where the normal vector needed to compute the traction is defined. By contrast, the proposed method does not exhibits this issue, as the solution can be evaluated directly at any point in the domain.

\begin{figure}[H]
	\centering
	\subfloat[Stress - PINN]{\includegraphics[width=0.22\linewidth]{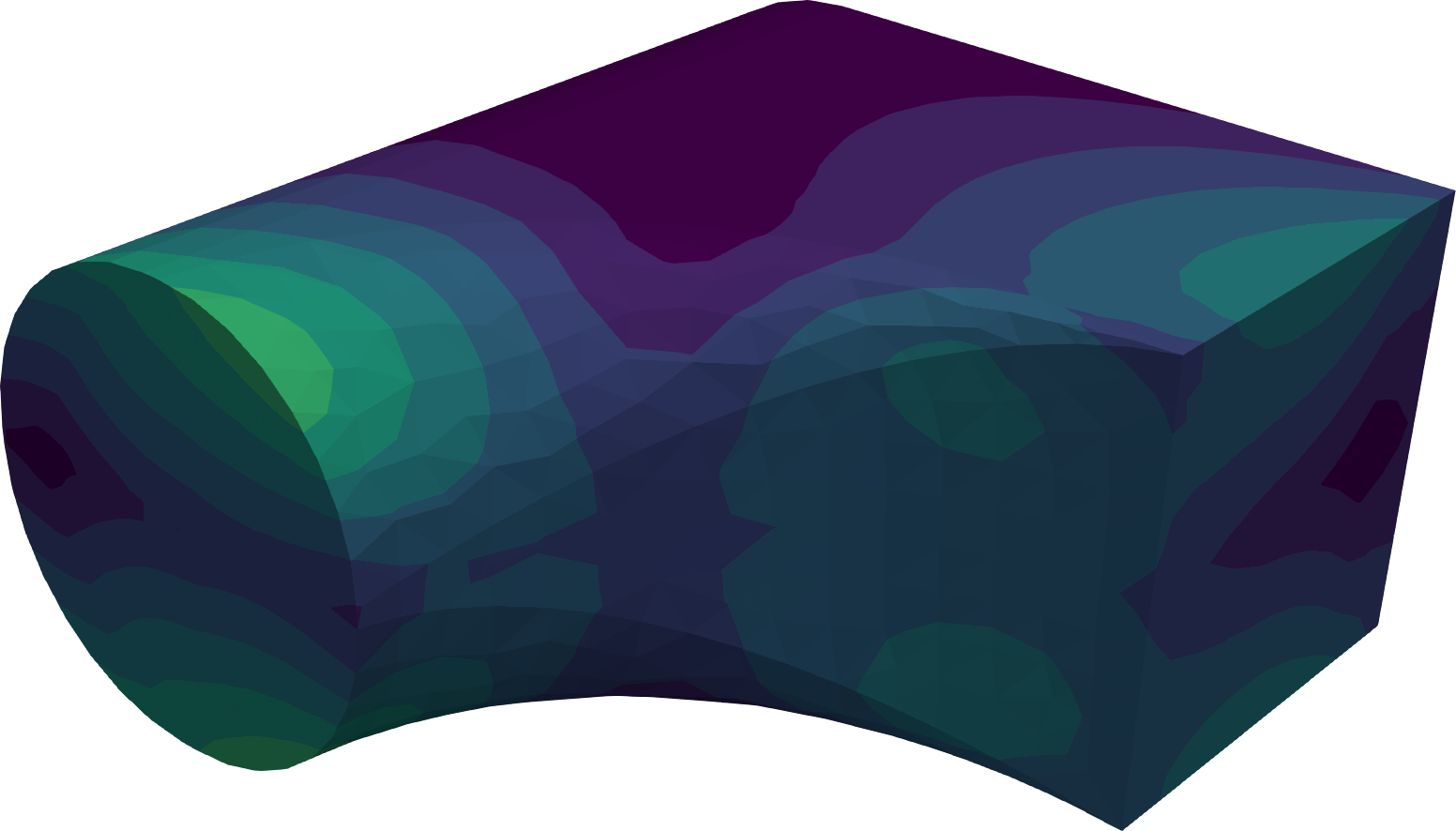}}\hspace{1mm}
	\includegraphics[width=0.08\linewidth]{fig/test_mech_device_stress_label.png}\hspace{5mm}
	\subfloat[Equilibrium - PINN]{\includegraphics[width=0.25\linewidth]{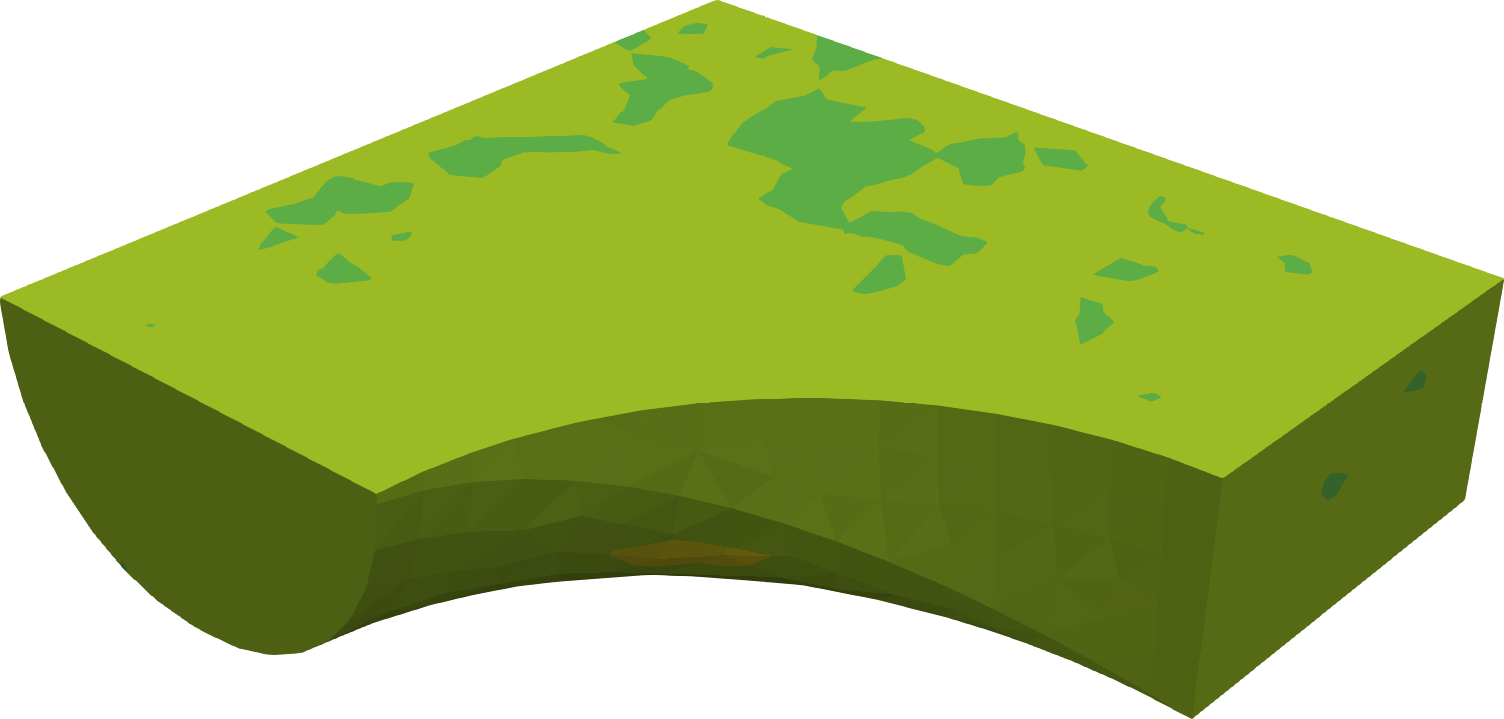}}
    \hspace{5mm}
	\subfloat[Equilibrium - HOL]{\includegraphics[width=0.25\linewidth]{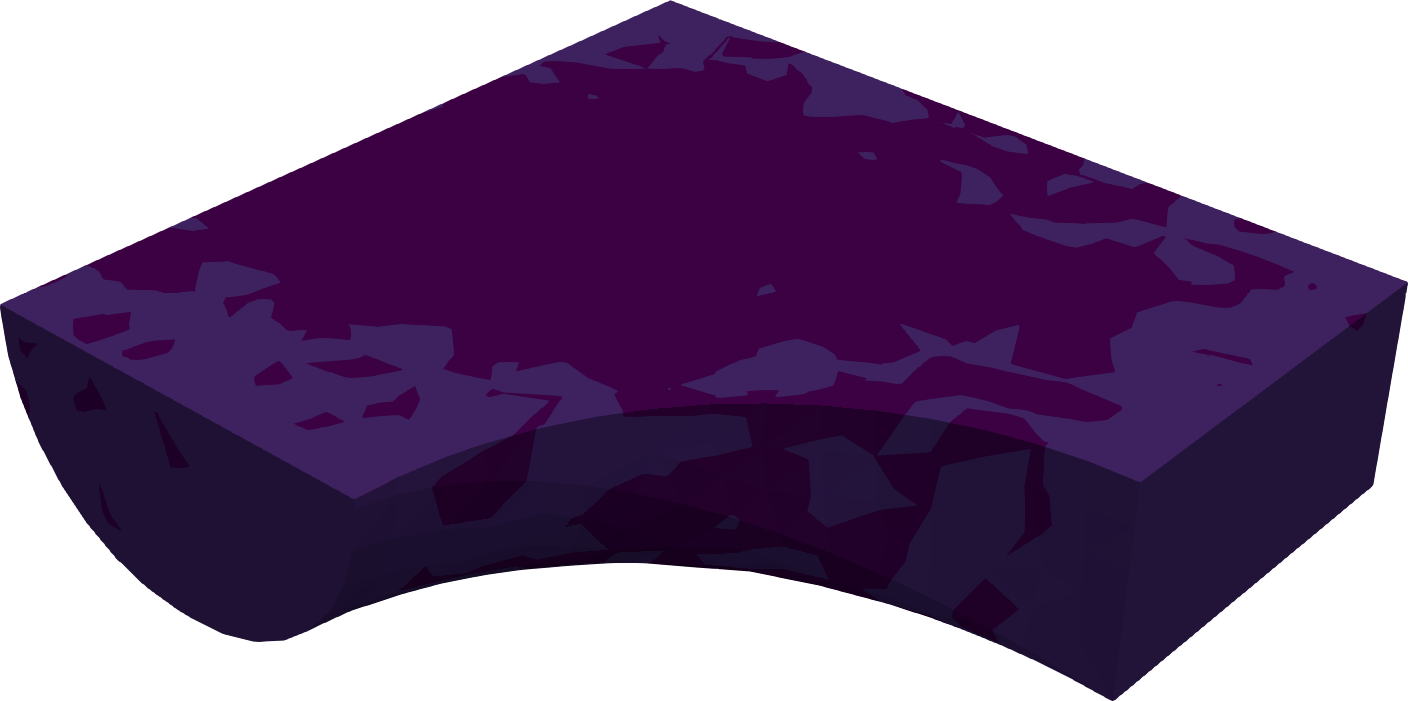}}\hspace{1mm}
	\includegraphics[width=0.08\linewidth]{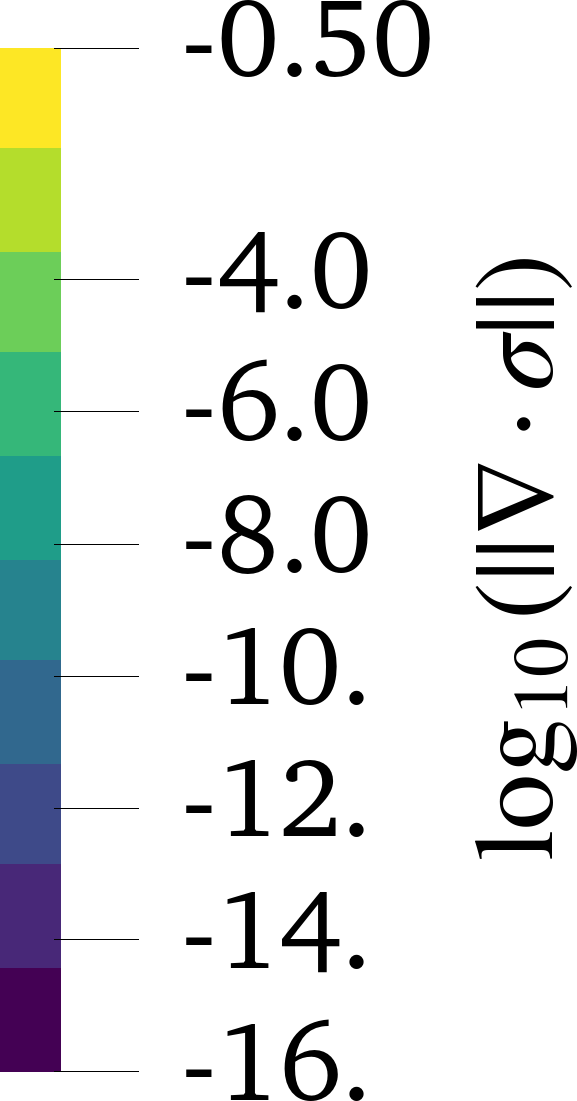} \\
	\caption{Leftmost panel: stress obtained with a standard PINN approach. Fairly accurate stress field is obtained using this approach. Central and right panels: visualization of the residual of the equilibrium equation, namely $\|\divergence \bfsigma - 0\| = \|\divergence \bfsigma\|$, in logarithmic scale. The proposed method, denoted by HOL, satisfies the equilibrium equation exactly as the residual is of the order of machine precision. For the standard PINN approach instead, the satisfaction of the equilibrium depends on the quality of the training.}
	\label{fig:test_mech_device_pinn}
\end{figure}

\paragraph{Evaluation against PINN.}
As done in \Cref{sec:test_lapl_fluid}, we present a comparison with a standard PINN formulation to better assess the properties of the proposed method. In this case, we use a single neural network to approximate the mapping $(x,y,z) \to (u_x,u_y,u_z)$.
We design the network to have both a comparable number of layers and total trainable parameters to $N_{\phi_x}$, $N_{\phi_y}$, $N_{\phi_z}$, and $N_\chi$, which is equal to $4\times \numprint{12737} = \numprint{50948}$. Accordingly, the network used for the PINN approach has layers of width $(3, 128, 128, 128, 128, 3)$, yielding a total number of trainable parameters equal to $\numprint{50435}$, as summarized in \Cref{tab:pinn_hol_mech}. 
The SiLU activation function is employed, ensuring sufficient smoothness of the displacement field. Preliminary experiments using other common smooth activation functions, e.g., softsign and tanh, produced ineffective results for this test case.
The PINN loss function is defined as
\begin{align*}
	&\loss_\text{PINN} = \sum_{i=1}^{n_\text{train}} \loss_p(\bfx_i)  + \sum_{i=1}^{n_{\text{train},\Omega}} \loss_\text{in}(\bfx_i), \\
    & \loss_\text{in}(\bfx) = \|
    \mu \laplacian \bfu + (\lambda + \mu)\nabla(\divergence \bfu)\|^2,
\end{align*}
where $\loss_p$ is defined in \eqref{eq:pointwise_loss_hol}, and $n_{\text{train},\Omega}$ denotes the number of training points in the domain, equal to $\numprint{37971}$. Note that, in this case $n_{\text{train}, \Omega} \approx 14\, n_\text{train}$ and in the proposed method, these training points in the interior domain are not required, significantly reducing the total number of training points, as summarized in \Cref{tab:pinn_hol_mech}. \\
%
%
In \Cref{fig:test_mech_device_pinn}, we present the stress field obtained with the PINN approach, which is similar to the FEM solution shown in \Cref{fig:test_mech_non_uniform}. Quantitative results are reported in \Cref{tab:pinn_hol_mech}, which presents the global errors on the displacement and stress fields for both the PINN and the proposed method. 
The two methods yield comparable overall accuracy. 
The PINN proves marginally more accurate in reproducing the displacement field 
($\overline{\eps}_{\bfu} = 0.04595$ for PINN versus $\overline{\eps}_{\bfu} = 0.04823$ for HOL), 
whereas it is less accurate in reproducing the stress field 
($\overline{\eps}_{\bfsigma} = 0.3117$ for PINN versus $\overline{\eps}_{\bfsigma} = 0.2527$ for HOL).
This discrepancy may be attributed to the fact that, in the PINN, the displacement is approximated 
directly as the output of a neural network, whereas in the proposed method it is recovered through 
a less direct computation, as shown in~\eqref{eq:displ_phi_chi}.
Concerning efficiency, \Cref{tab:pinn_hol_mech} shows that, contrary to the results seen for the Laplace problem of \Cref{sec:test_lapl_fluid}, the proposed method provides only a moderate speed-up (674
min for PINN versus 575 min for HOL, see \Cref{sec:appendix_non_uniform} for further details). A possible explanation is that, as noted in \Cref{sec:3d_elasticity}, the Papkovich-Neuber representation increases the number of unknown functions from three to four, which may confer an advantage to the PINN approach.\\
Finally, to highlight a key difference between the proposed method and the PINN approach, we consider the equilibrium equation:
\begin{equation}\label{eq:equilibrium}
    \divergence\bfsigma = 0.
\end{equation}
The satisfaction of \eqref{eq:mech_PN_P}--\eqref{eq:mech_PN_B} implies an equilibrated solution; therefore, \eqref{eq:equilibrium} is satisfied exactly \cite{Sadd2021}.
The residual of \eqref{eq:equilibrium} is shown in \Cref{fig:test_mech_device_pinn}, where we can observe that the equilibrium is exactly satisfied with the proposed method, whereas in the PINN approach the precision with which the governing equations are satisfied depends on the quality of the training.
\begin{table}[h]
\centering
\begin{tabular}{llllll}
     & $|\mathcal{W}_\text{tot}|$ & Tot. training points & Training \sib{\min} & $\overline{\varepsilon}_{\bfu}$ & $\overline{\varepsilon}_{\bfsigma}$  \\
     \hline
 PINN &  \numprint{50435} & \numprint{45716} &  674 & 0.04594 & 0.3117  \\
 HOL  & \numprint{50948} & \numprint{7745} & 575 & 0.04823 & 0.2527 \\ 
\end{tabular}
\caption{Comparison between the physics-informed neural network (PINN) approach and the proposed method (HOL), in terms of number of trainable parameters, number of training points, training time, and global relative errors defined in \eqref{eq:errs_mech_global_u}-\eqref{eq:errs_mech_global_sigma}.}
\label{tab:pinn_hol_mech}
\end{table}

\clearpage

\section{Conclusion}

This work introduced a neural-network-based framework for the solution of three-dimensional Laplace-type problems, as well as for a broader class of boundary value problems that admit a reformulation in terms of harmonic potentials. \\
From a methodological standpoint, the proposed framework can be regarded as a physics-informed neural network approach. As in standard PINNs, neural networks are employed as global ansatz spaces and the training procedure is driven by the governing equations. However, substantial differences compared to conventional PINN approaches exist. In the latter, the PDE is enforced approximately through the minimization of a residual evaluated at interior collocation points. In the proposed approach instead, the PDE residual vanishes identically by construction, independently of the optimization process, training accuracy, and numerical quadrature adopted for the integral representation. As a consequence, the obtained approximation satisfies exactly the underlying governing equations. Moreover, since no interior residual evaluation is required, training involves only collocation points located on the boundary. This avoids the computation of higher-order derivatives associated with the PDE residual and leads to a more favorable scaling of the number of training points with respect to the characteristic discretization length. \\
The starting point of the proposed framework is the integral representation of harmonic functions introduced by Whittaker. This representation leads to reformulations of the original boundary value problems, where the unknown fields are expressed in terms of functions that are holomorphic with respect to the first argument, while remaining generally non-holomorphic with respect to the second argument. This structure motivated the introduction of the \textit{semi-holomorphic} neural network architecture, specifically designed to preserve holomorphicity in the first input component, while allowing nonlinear interaction with the second, non-holomorphic, input.

The proposed methodology has been validated through several numerical experiments involving Laplace and linear elasticity problems under different geometries and boundary conditions. The results show that the method is capable of accurately reproducing both scalar and vector fields in complex three-dimensional settings. Comparison with standard PINN formulations highlights both strengths and limitations of the proposed approach. The exact satisfaction of the governing equations constitutes a significant advantage, particularly for solid mechanics applications in which equilibrium and compatibility are fundamental physical constraints. Moreover, the numerical results demonstrate a reduction in training time compared with the PINN method. The observed speedup is approximately a factor of 5 for the considered Laplace problem and about 15\% for the considered linear elasticity problem. The comparatively modest improvement in the latter case may be partly attributed to the adopted Papkovich-Neuber representation, which increases the number of unknown functions from three to four.

Some limitations and open questions remain. For example, the dependence of the solution on the angular variable $\theta$ may be discontinuous or exhibit sharp transitions, whereas neural networks are naturally better suited for approximating smooth functions. This issue may limit the approximation capabilities of the present architecture in specific cases. Additional developments are also required from the computational standpoint. Although the method avoids the evaluation of PDE residuals in the domain, the training process remains relatively slow, which is a common limitation of PINN-type approaches. Different optimization algorithms, pruning techniques \cite{Ballini2026}, or problem-tailored acceleration procedures may improve the overall computational performance. Along these lines, we note that, for the linear elasticity problem, the number of potentials, and thus the number of neural networks, can be reduced under suitable hypotheses \cite{TranCong1989}. This would shrink the search space and potentially improve training efficiency. Also, further investigations may benefit from alternative complex representations \cite{Piltner2001,MaresCarreno2021}, or from the use of different potential functions in both the stress and displacement formulations \cite{Sadd2021}. 

Overall, the present work shows that incorporating analytical structures directly into the neural network architecture provides an effective strategy for the solution of three-dimensional PDEs. More generally, the proposed framework illustrates how classical analytical representations and modern neural-network-based methodologies can be combined in a mathematically consistent way, preserving intrinsic properties of the governing equations while retaining the flexibility of meshless learning approaches. This can potentially open new research directions, both in the mathematical analysis of property-preserving neural network approximations and in the development of computational methods for PDEs.

\section*{Acknowledgement}
This work was supported by the Danish Research Council for Independent Research through the grant no. 2035-00142B ``Network-inspired models to predict the strength of heterogeneous materials''. The authors declare that they have no known competing financial interests or personal relationships that could have appeared to influence the work reported in this paper.

\section*{Data availability}
The code for all test cases is available at
\url{https://github.com/enricoballini/Holomorphic-NN-for-3D-BVP.git}.

\clearpage

\bibliographystyle{cas-model2-names-unsrt} 
\bibliography{all}

\clearpage

\appendix
\section{Computational framework and implementation details}\label{sec:hyperparameters_details}

We provide in this appendix information on practical aspects of the test cases presented in \Cref{sec:tests}.
The code for all test cases is available, cf. info in the Data availability section. It is written in Python 3 and relies on JAX \cite{jax2018github} neural network routines and uses CUDA \cite{cuda} for GPU acceleration.
The computation of the discrete form of $\overline{V}$, \eqref{eq:V_discrete}, and $\bfu$, \eqref{eq:displ_phi_chi}, is parallelized over the quadrature points and the training points within each minibatch. The functions related to the weight updates are compiled using \texttt{jax.jit}. We remark that optimizing the code to achieve the best performance is beyond the scope of this work as well as accurate timing and evaluation of the computational costs.

\subsection{Implementation details - Manufactured solution on complex geometry}\label{sec:appendix_support}
\begin{figure*}[h]
	\centering
	\includegraphics[width=0.48\linewidth]{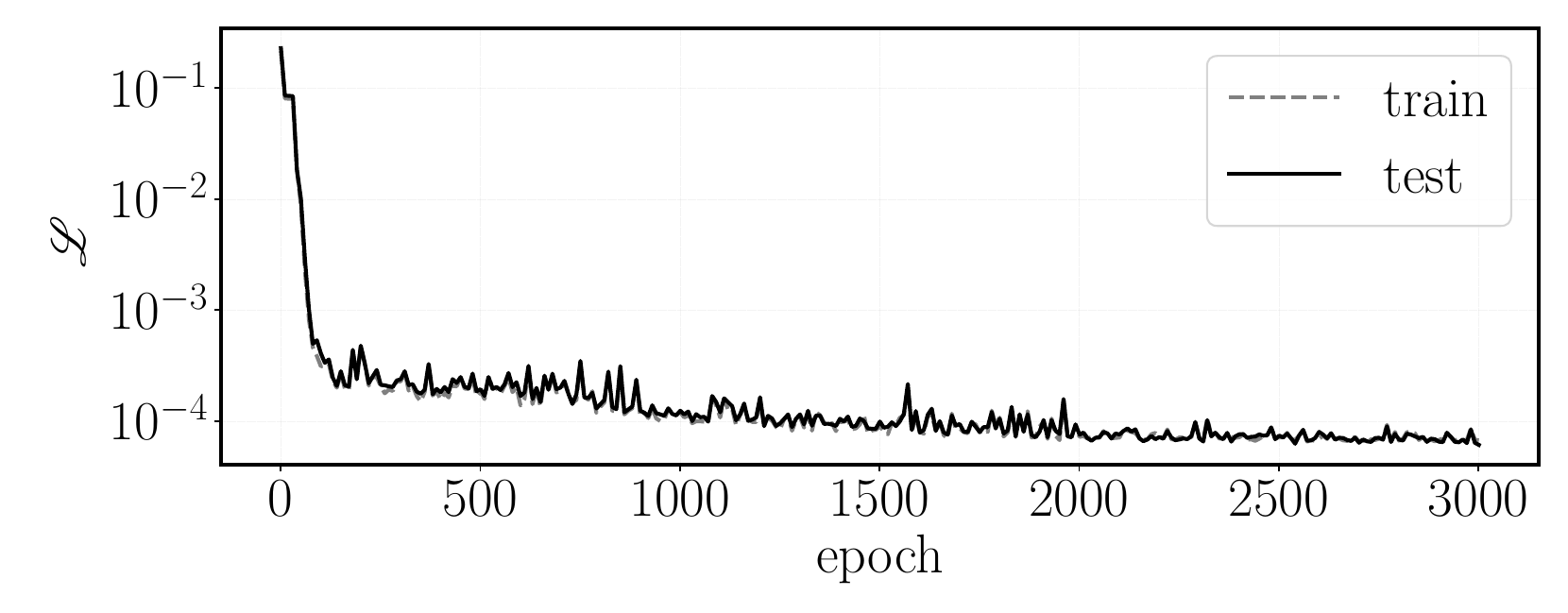}
	\caption{Value of the loss versus the epochs for the test case presented in \Cref{sec:test_lapl_complex}.}
	\label{fig:loss_lapl}
\end{figure*}

This test case was run on an NVIDIA GeForce RTX 4060. The training computations are performed in single precision, while the post-processing is carried out in double precision.  By single precision, we mean that both the real and imaginary parts of any complex number are represented as single-precision floating-point numbers. Similarly for double precision. \\
The characteristic dimension of the domain is approximately $300$ length units, which is not a suitable input scale for neural networks. The coordinates are therefore scaled by a factor of $300$ before being provided to the networks. \\
The weights are initialized following the strategy proposed in \cite{Calafa2024}.
During training with the ADAM algorithm, with hyperparameters $\beta_1=0.9$, $\beta_2=0.999$, and $\epsilon=10^{-8}$ (see \cite{Kingma2014}), the training data are shuffled at each epoch. The minibatch has size $32$. The training lasts a few minutes. Preliminary results shows that no meaningful advantages are provided by adopting different first-order or second-order optimization strategy, e.g. AMSgrad \cite{Sashank_amsgrad} or L-BFGS \cite{Liu1989_LBFGS}. \\
The learning rate, $\ell_r$, is varied according to the following schedule:
\begin{equation*}
	\ell_r = \begin{cases}
		10^{-4} & \text{if } \ n_\text{epoch} < 1000, \\
		5\times10^{-5} & \text{if } \ 1000 \leq n_\text{epoch} < 2000, \\
		10^{-5} & \text{if } \  2000 \leq n_\text{epoch},
	\end{cases}
\end{equation*}
where $n_\text{epoch}$ denotes the number of epochs. \\
In \Cref{fig:loss_lapl}, the loss function computed using the point on the training dataset and the same function but using the point of the test dataset are reported. The training is carried out for a large number of epochs to ensure a stable value of the loss in the final stages. No overfitting is observed.

\subsection{Implementation details - Neumann problem in a quarter torus}\label{sec:appendix_flow}
This test case was run on an NVIDIA GeForce RTX 4060. The training computations are performed in single precision, while the post-processing is carried out in double precision.
The weights of the network associated to the PINN method are initialized following the strategy proposed in \cite{Glorot2010} while the weights of the network associated to the proposed method are initialized following the strategy proposed in \cite{Calafa2024}. Different seeds are used for different networks. \\
The following hyperparameters are used for both the PINN and the proposed method. During training with the ADAM algorithm, with hyperparameters $\beta_1=0.9$, $\beta_2=0.999$, and $\epsilon=10^{-8}$ (see \cite{Kingma2014}), the training data are shuffled at each epoch. The minibatch has size $64$. \\
The learning rate, $\ell_r$, is varied according to the following schedule:
\begin{equation*}
	\ell_r = \begin{cases}
		5\times10^{-3} & \text{if } \ n_\text{epoch} < 1000, \\
		10^{-5} & \text{if } \ 1000 \leq n_\text{epoch}. \\
	\end{cases}
\end{equation*}
In \Cref{fig:loss_lapl}, the loss values are reported, both for the PINN strategy, in panel (a), and the proposed method, in panel (b). The training is carried out for a large number of epochs to ensure a stable value of the loss in the final stages. No overfitting is observed.
\begin{figure*}[h]
	\centering
	\subfloat[PINN]{\includegraphics[width=0.48\linewidth]{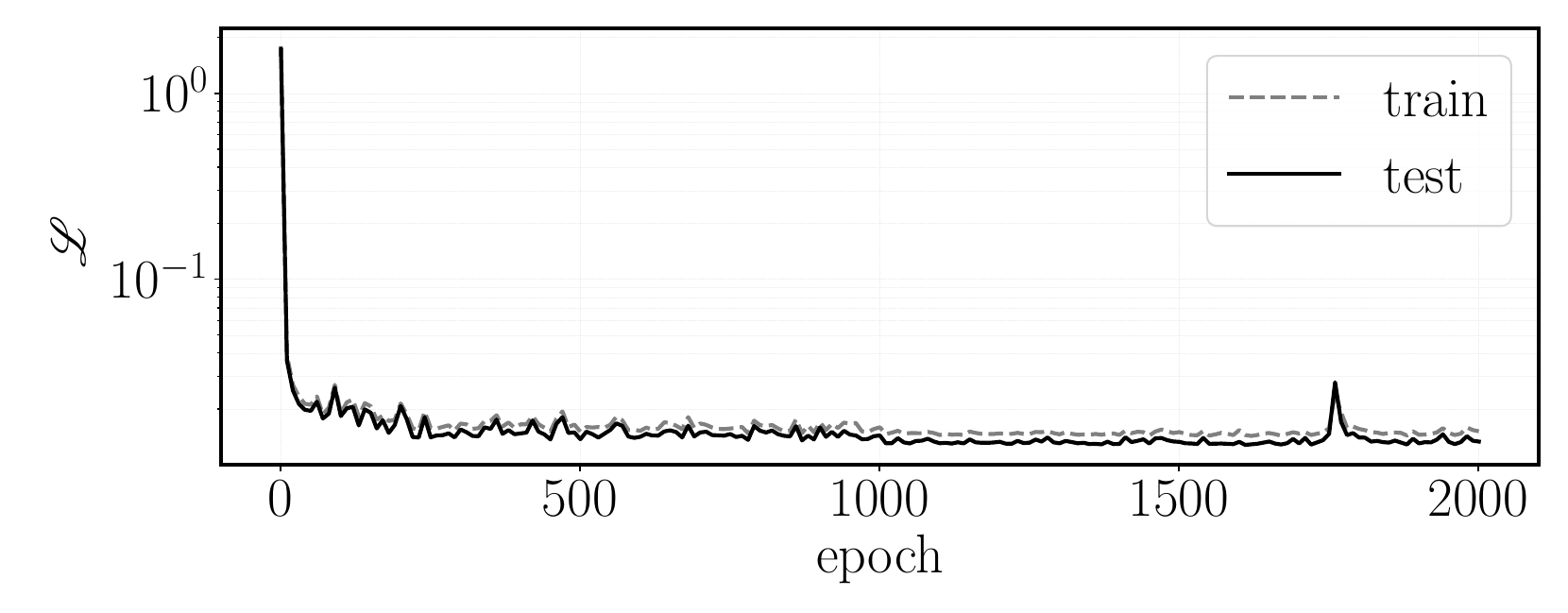}}\hspace*{4mm}
	\subfloat[HOL]{\includegraphics[width=0.48\linewidth]{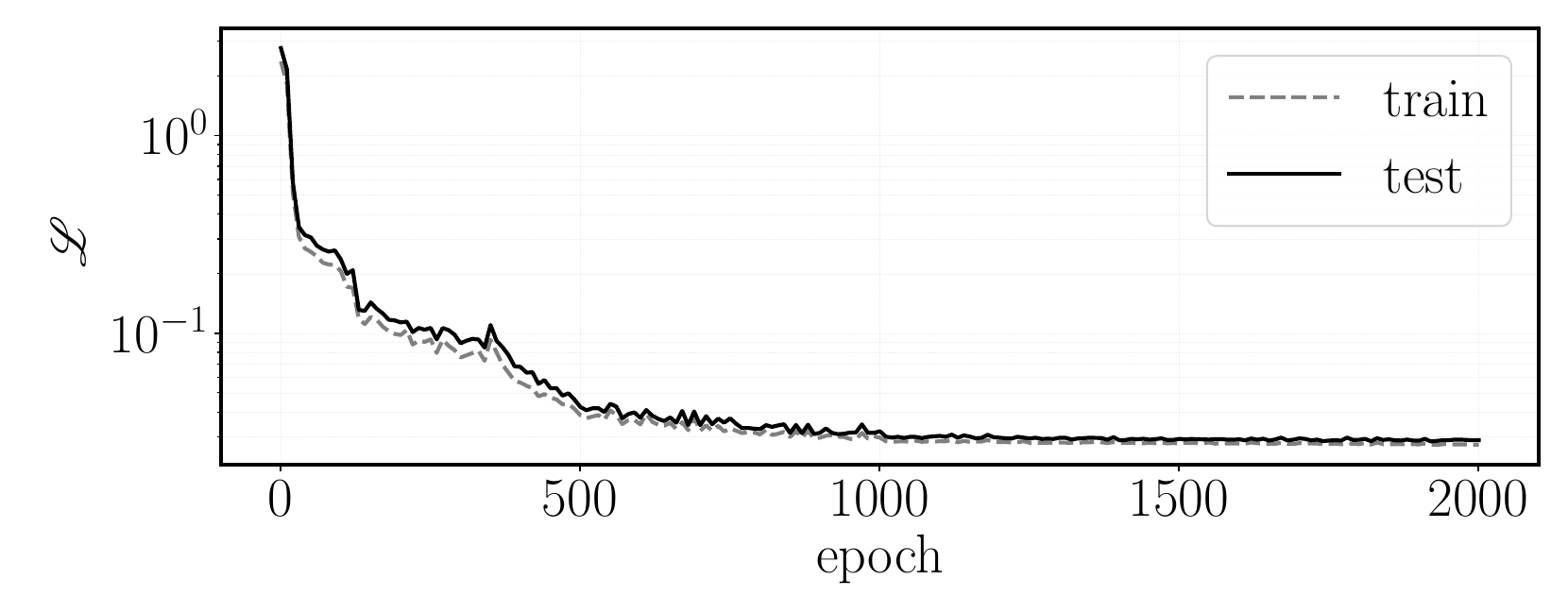}}
	\caption{Value of the loss versus the epochs for the test case presented in \Cref{sec:test_lapl_fluid}.}
	\label{fig:loss_lapl_flow}
\end{figure*}

\subsection{Implementation details - Elementary cases under uniform stress}\label{sec:appendix_uniform}
\begin{figure}[h]
	\subfloat[Loss]{\includegraphics[width=0.42\linewidth]{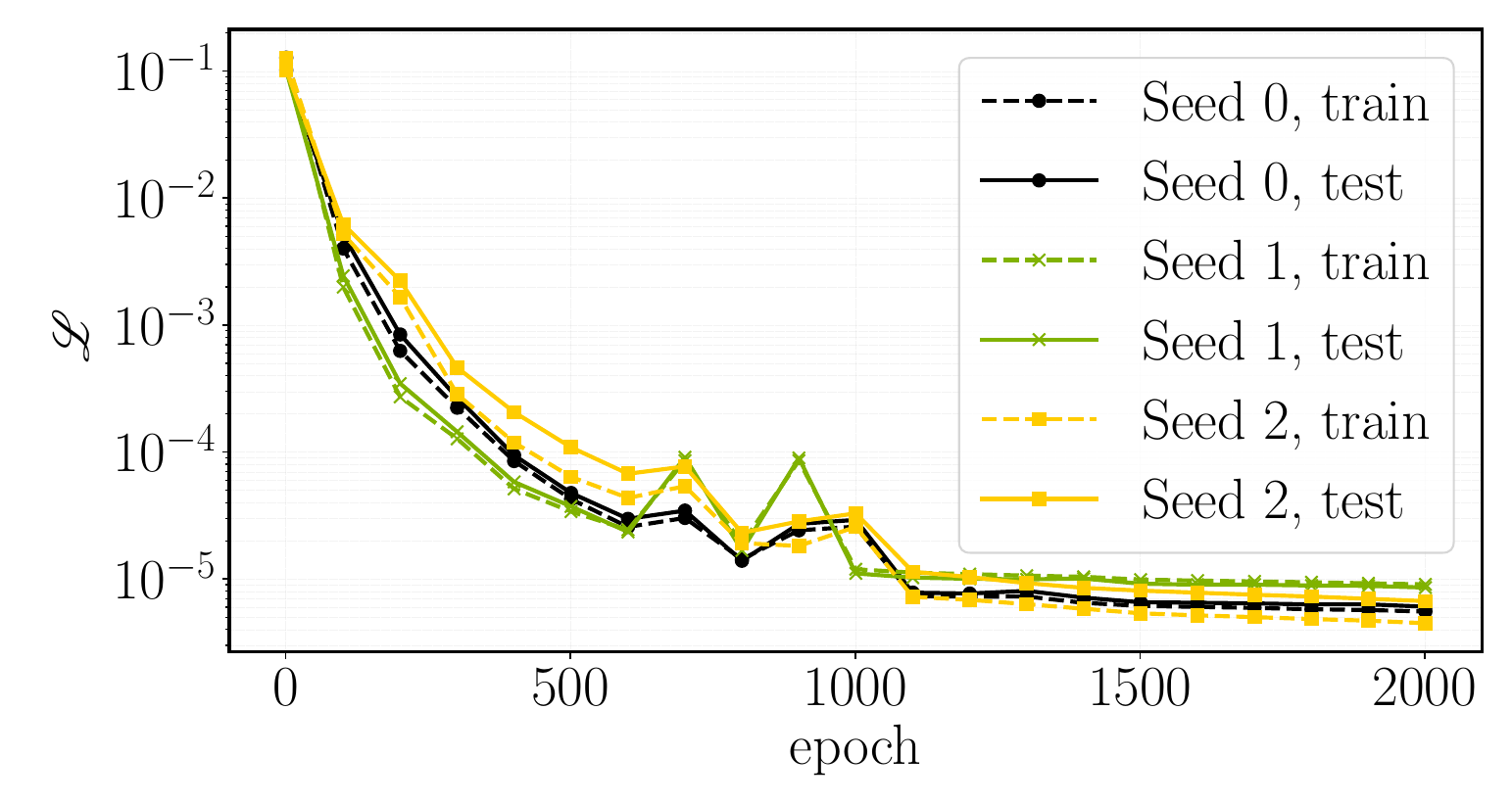}}\hspace*{5mm}
	\begin{minipage}{0.58\textwidth}
		\vspace*{-34mm}
		\subfloat[Stability of the layers]{\adjustbox{valign=t}{\includegraphics[width=0.7\linewidth]{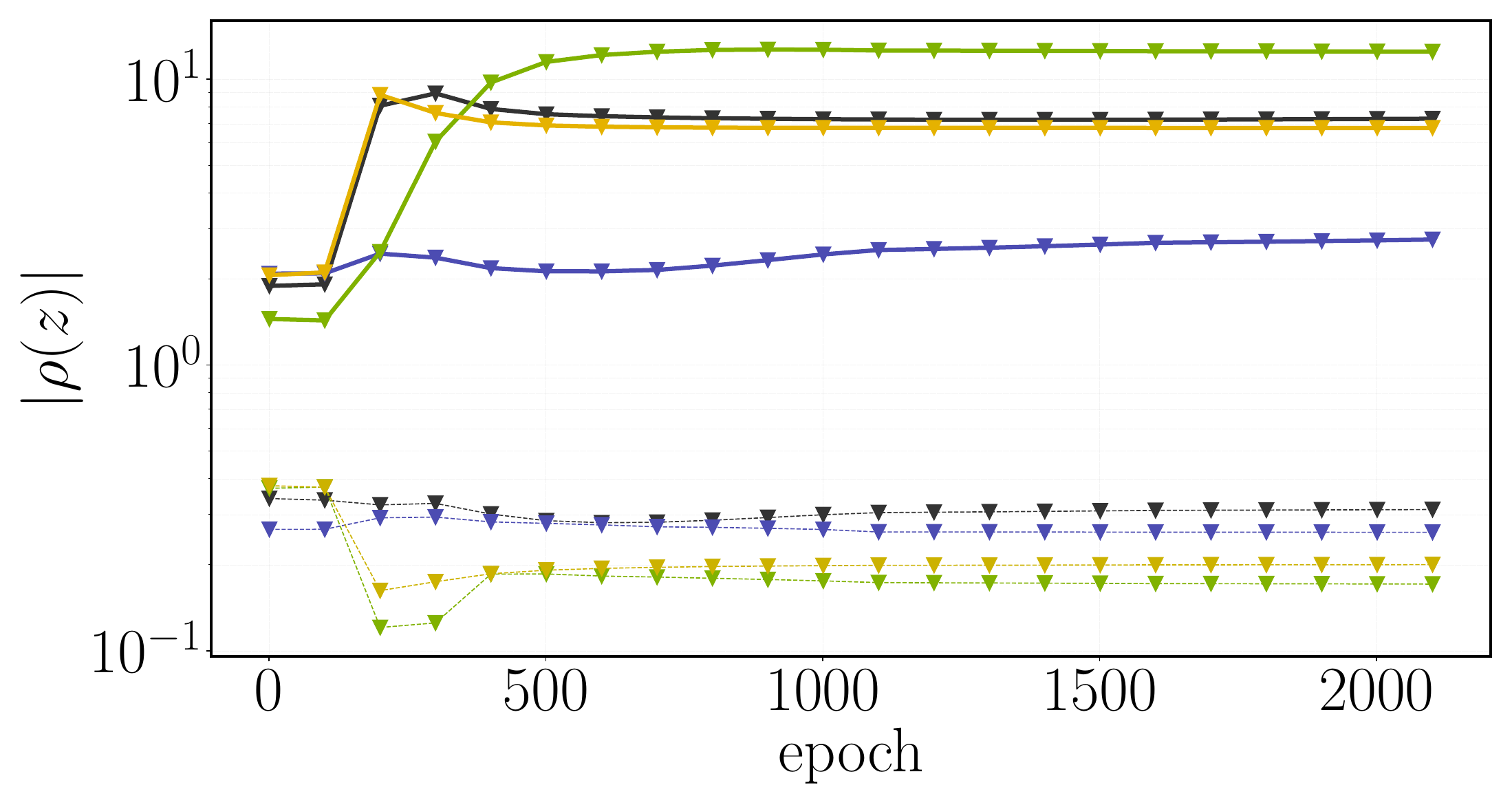}}
		\adjustbox{valign=t}{\includegraphics[width=0.25\linewidth]{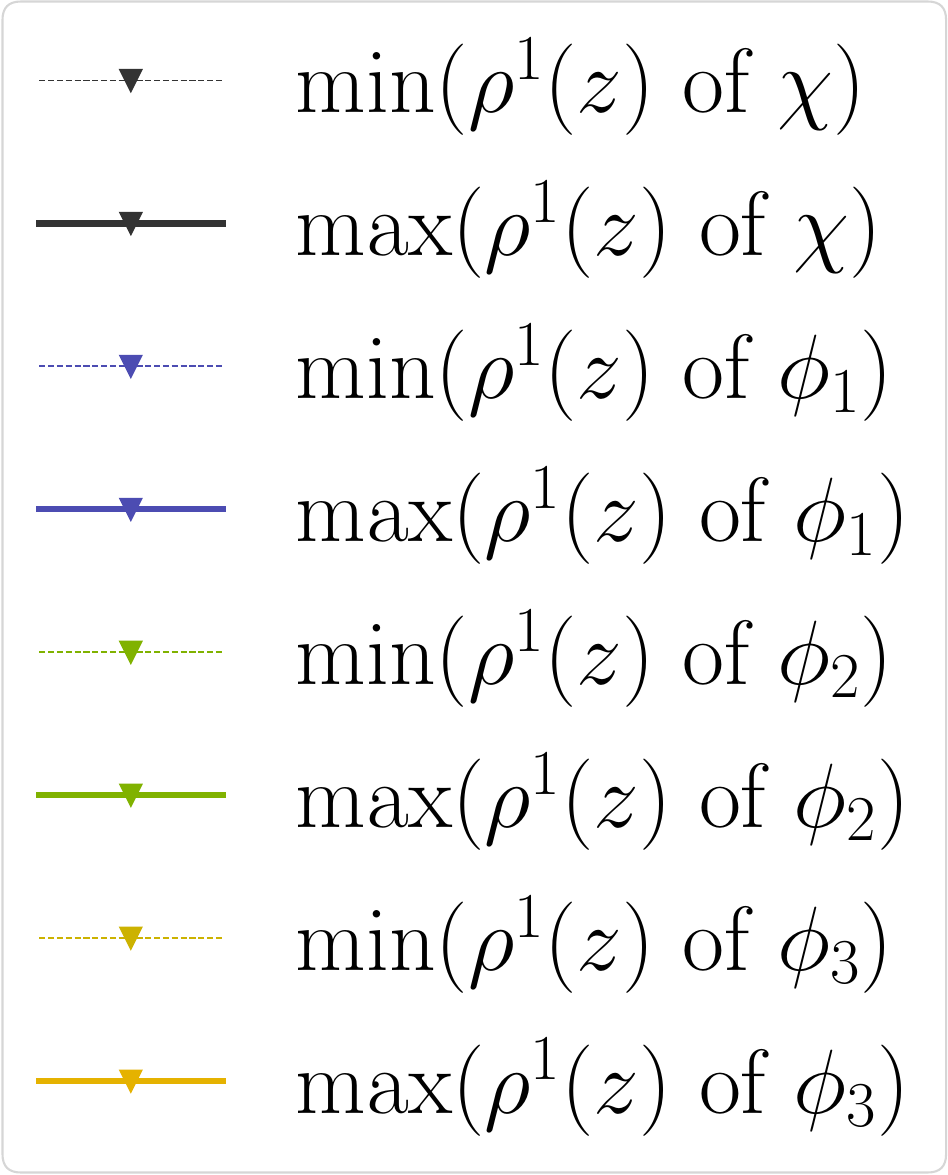}}}
	\end{minipage}
	\caption{(a) Loss versus epoch for different weight initialization. (b) maximum and minimum module of the output of the activation function for the hidden layer, i.e. layer number $1$.}
	\label{fig:loss_uniform}
\end{figure}

Each subcase of this test was run on an NVIDIA GeForce RTX 4060. The training computations are performed in single precision, while the post-processing is carried out in double precision.  \\
The same training hyperparameters described below are used for each subcase.
The weights are initialized following the strategy proposed in \cite{Calafa2024}. Different seeds are used for different networks. \\
During training with the ADAM algorithm, with hyperparameters $\beta_1=0.9$, $\beta_2=0.999$, and $\epsilon=10^{-8}$ (see \cite{Kingma2014}), the training data are shuffled at each epoch. The minibatch has size $32$. The training lasts a few seconds. \\
The learning rate, $\ell_r$, is varied according to the following schedule:
\begin{equation*}
	\ell_r = \begin{cases}
		10^{-3} & \text{if } \ n_\text{epoch} < 1000, \\
		2\times10^{-4} & \text{if } \ 1000 \leq n_\text{epoch} < 1500, \\
		5\times10^{-5} & \text{if } \  1500 \leq n_\text{epoch}.
	\end{cases}
\end{equation*}
In \Cref{fig:loss_uniform}, the loss values are reported. The training is carried out for a large number of epochs to ensure a stable value of the loss in the final stages. No overfitting is observed. Additionally, for this simple shear case (see \Cref{sec:test_mech_uniform}), the training was repeated $3$ times with different weight initialization. The loss value is fairly stable among the different initializations as visible in \Cref{fig:loss_uniform}.
To ensure stability of the updates, denoted by $\Delta W$, we apply a clipping operation to the weight updates, defined as follows:
\begin{align}\label{eq:clip_update}
	 \Delta W \leftarrow \operatorname{clip}(\real(\Delta W), t_h) + i\operatorname{clip}(\imag(\Delta W), t_h), 
\end{align} 
where the function clip is 
\begin{align*}
	\text{clip}(x, t_h) =
	\begin{cases}
		-t_h & \text{if } x < -t_h, \\
		x & \text{if } -t_h \le x \le t_h, \\
		t_h & \text{if } x > t_h.
	\end{cases}
\end{align*}
We set $t_h = 10^{-2}$. \Cref{fig:loss_uniform} also reports the minimum and maximum values of the modulus of the output of the activation function in the hidden layers of the neural network. These values are depicted for each epoch. They remain stable throughout training, with the maximum across the four networks being approximately $10$. This indicates overall stability of the training process, despite the use of an exponential activation function, which could potentially lead to uncontrolled large outputs.

\subsection{Implementation details - 3D device subjected to non-uniform sterss}\label{sec:appendix_non_uniform}
\begin{figure}[h]
	\centering
	\subfloat[Training and test points.]{\includegraphics[width=0.3\linewidth]{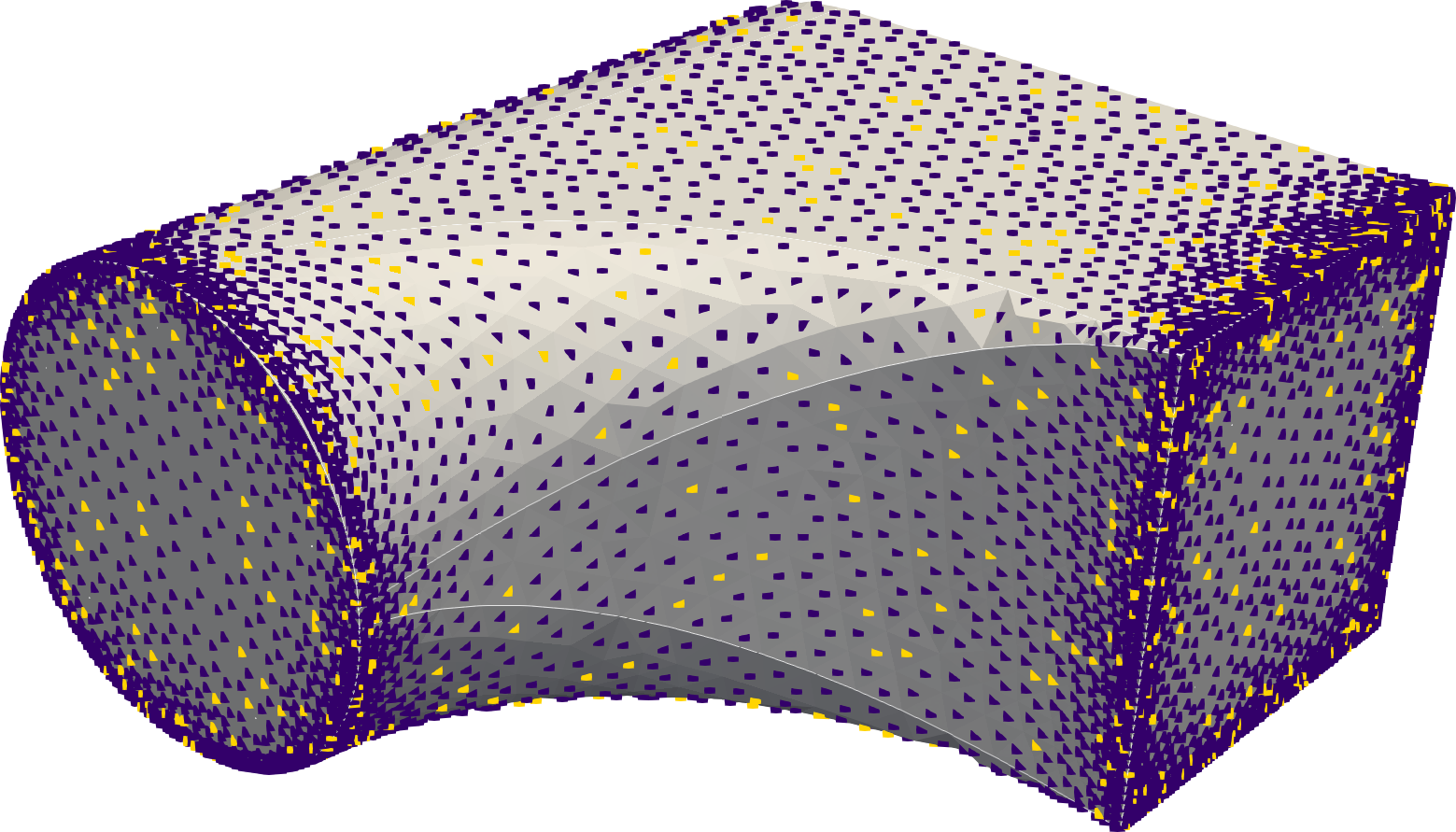}\hspace{5mm}} \\
    \subfloat[Loss - PINN]{\includegraphics[width=0.45\linewidth]{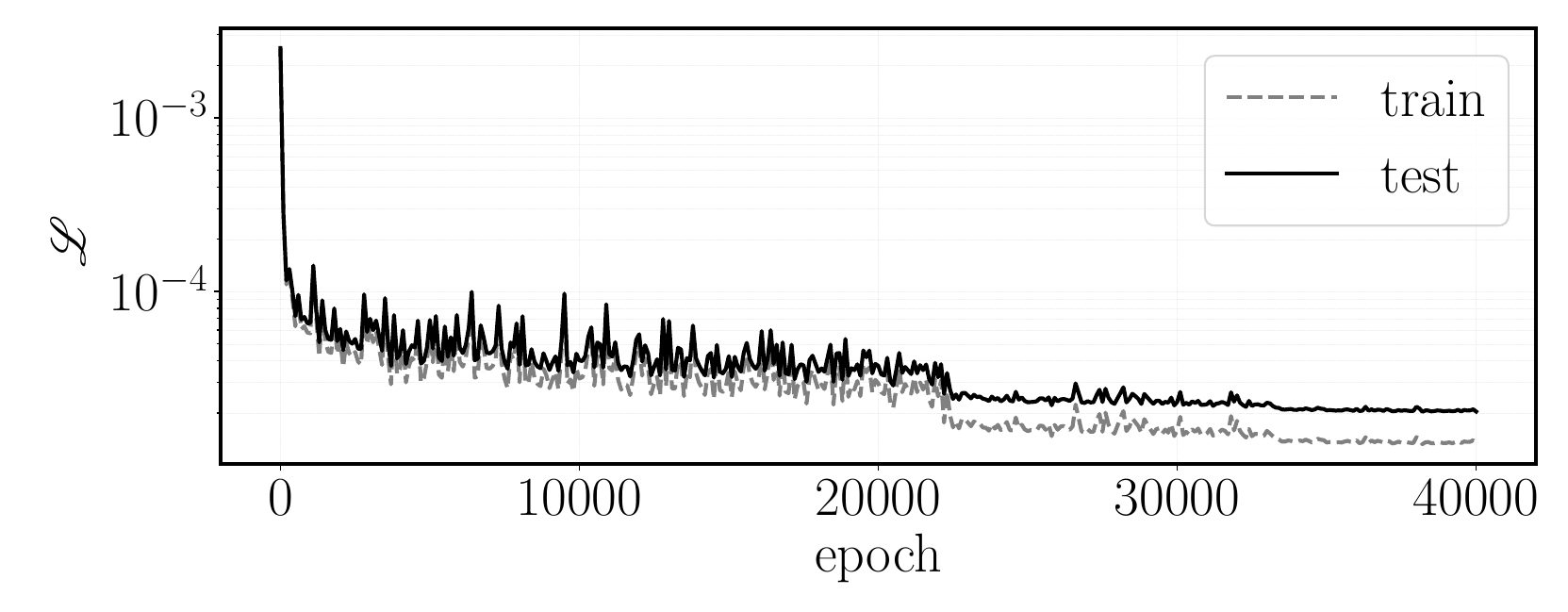}} 
	\subfloat[Loss - HOL]{\includegraphics[width=0.45\linewidth]{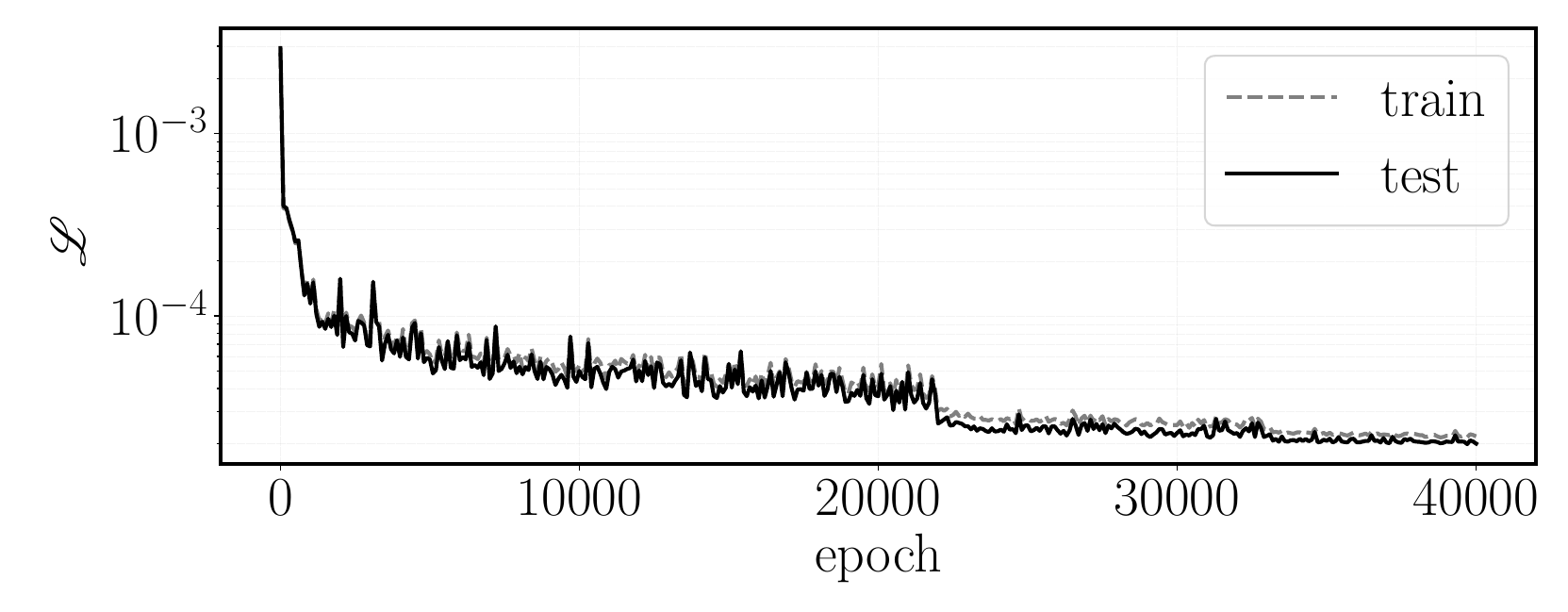}} \\
	\subfloat[Stability of the layers]{\adjustbox{valign=c}{\includegraphics[width=0.5\linewidth]{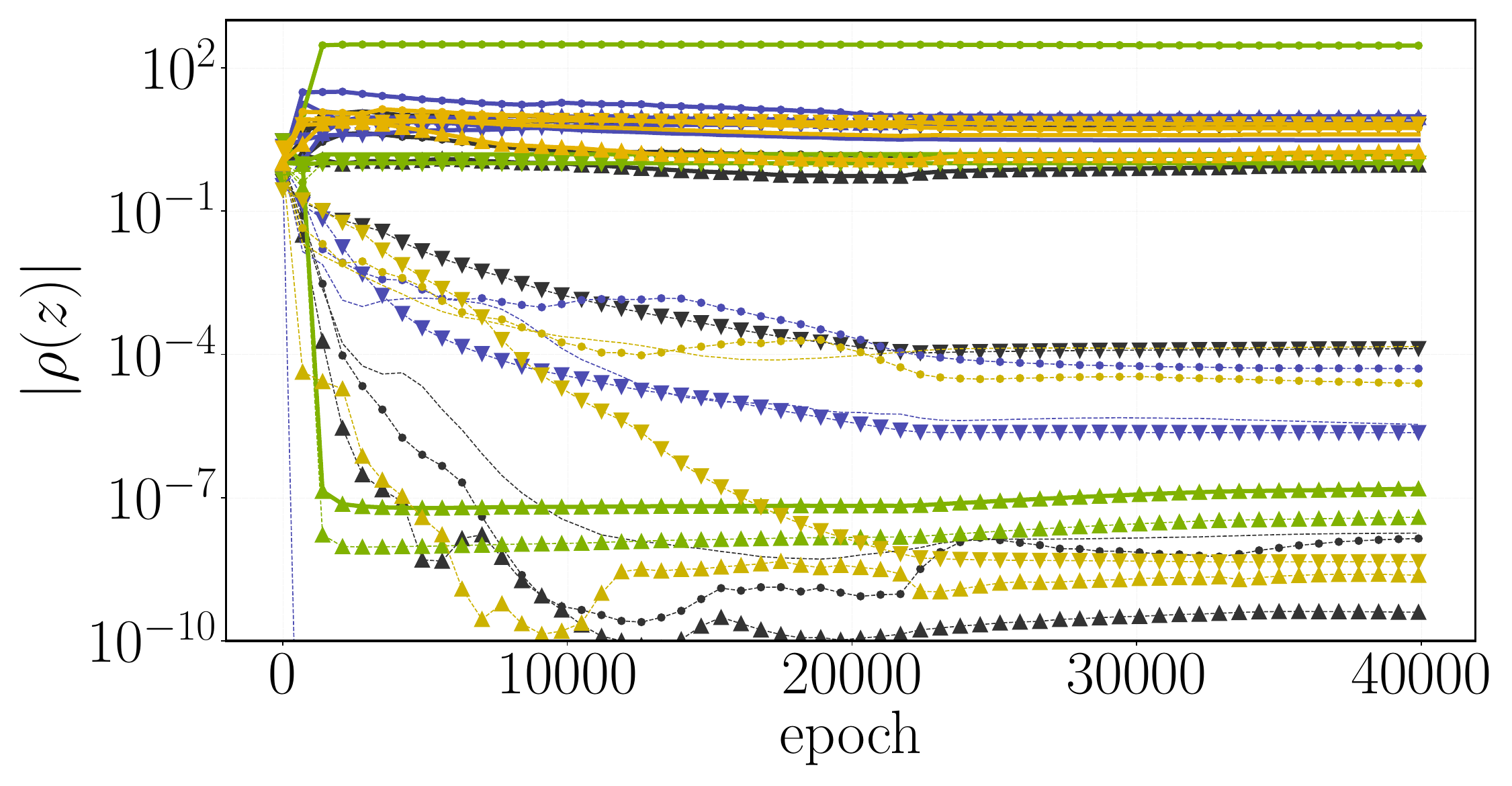}}
	\adjustbox{valign=c}{\includegraphics[width=0.30\linewidth]{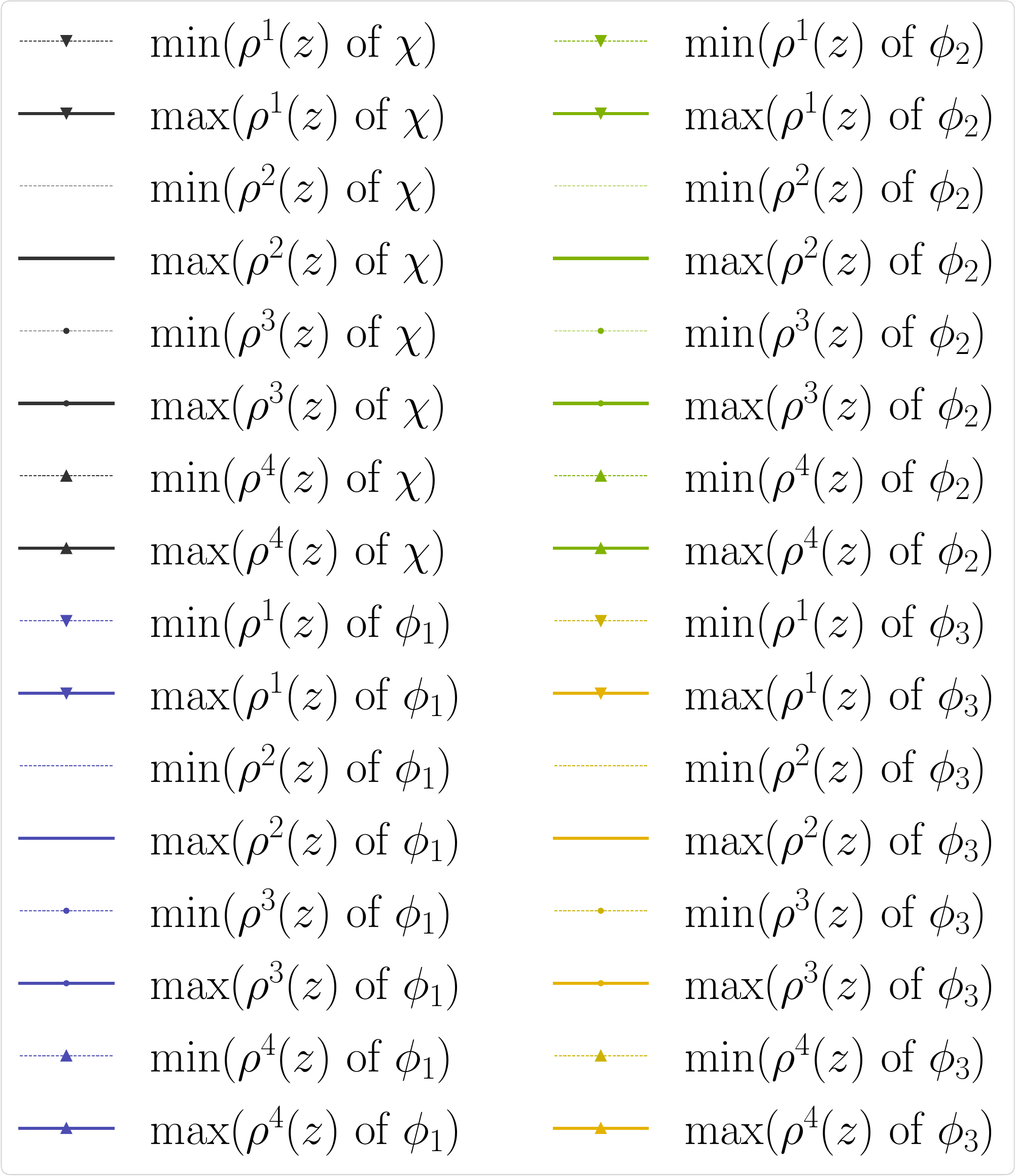}}}	
	\caption{(a) Training points in dark color, test points in light color. (b) Loss versus epoch for different weight initialization. (c) maximum and minimum module of the output of the activation function for the hidden layer, i.e. layer $1$.}
	\label{fig:loss_non_uniform}
\end{figure}

This test was run on an Nvidia Tesla A100 GPU. Both the training and the post-processing computations are performed in double precision. The training and test coordinates point are shown in panel (a) of \Cref{fig:loss_non_uniform}. A fine point distribution is applied on the sharp borders, where the error of the stress is higher than the error on the other part of the domain, see \Cref{fig:test_mech_non_uniform}.  
The weights of the network associated to the PINN method are initialized following the strategy proposed in \cite{He2015} while the weights of the networks associated to the proposed method are initialized following the strategy proposed in \cite{Calafa2024}. Different seeds are used for different networks. \\
The following hyperparameters are used for both the PINN and the proposed method.
During training with the ADAM algorithm, with hyperparameters $\beta_1=0.9$, $\beta_2=0.999$, and $\epsilon=10^{-8}$ (see \cite{Kingma2014}), the training data are shuffled at each epoch. The minibatch has size $128$. \\
The learning rate, $\ell_r$, is varied according to the following schedule:
\begin{equation*}
	\ell_r = \begin{cases}
		10^{-3} & \text{if } \ n_\text{epoch} < 22000, \\
		2\times10^{-4} & \text{if } \ 22000 \leq n_\text{epoch} < 33000, \\
		5\times10^{-5} & \text{if } \  33000 \leq n_\text{epoch}.
	\end{cases}
\end{equation*}
In \Cref{fig:loss_non_uniform}, the loss values are reported. The training is carried out for a large number of epochs to ensure a stable value of the loss in the final stages. However, only a slightly higher loss is observed if the training was stopped after a shorter number of epochs. This would lead to comparable results in terms of $\bfu$ and $\bfsigma$, while reducing the training time.
To ensure stability of the updates, denoted by $\Delta W$, we apply a clipping operation to the weight updates, defined in \eqref{eq:clip_update}, with $t_h = 10^{-5}$.
Panel (c) of \Cref{fig:loss_non_uniform} reports the minimum and maximum values of the modulus of the output of the activation function in the hidden layers of the neural network. The maximum values remain stable throughout training. This indicates overall stability of the training process, despite the use of an exponential activation function and 4 hidden layers which could potentially lead to uncontrolled large outputs.

\end{document}